\documentclass[leqno,12pt]{article}
\usepackage{a4}
\usepackage{amsthm,amssymb}
\usepackage{amsbsy,amsfonts,amsmath}

\usepackage{latexsym}

\numberwithin{equation}{section}

\newcommand{\hh }{\hat{h}}
\newcommand{\hps}{\hat{\psi }}

\newcommand{\bh }{\bar{h}}
\newcommand{\bpsi }{\bar{\psi }}

\newcommand{\cL }{\mathcal{L}}
\newcommand{\cM }{\mathcal{M}}
\newcommand{\cP }{\mathcal{P}}
\newcommand{\cT }{\mathcal{T}}
\newcommand{\cU }{\mathcal{U}}
\newcommand{\cZ }{\mathcal{Z}}
\newcommand{\cF }{\mathcal{F}}

\newcommand{\dg }{{\rm{deg}}}

\newcommand{\cUpm }{\cU ^{\prime <0}}
\newcommand{\cUpo }{\cU ^{\prime \,0}}
\newcommand{\cUpp }{\cU ^{\prime >0}}

\newcommand{\circledplus }{
\setlength{\unitlength}{1mm}
\begin{picture}(4.2,4)
	\put(1.98,1.4){\circle{3}}
	\put(1,.9){\text{\tiny $+$}}
\end{picture}
}
\newcommand{\End }{\mathrm{End}}
\newcommand{\id }{\mathrm{id}}
\newcommand{\isv }{I}
\newcommand{\isor }{\hat{\delta}}
\newcommand{\Klf }{\mathcal{K}_4}
\newcommand{\lact }{\boldsymbol{\cdot }}
\newcommand{\mfg }{\mathfrak{g}}
\newcommand{\mfh }{\mathfrak{h}}
\newcommand{\mfn }{\mathfrak{n}}
\newcommand{\mfp }{\mathfrak{p}}
\newcommand{\mfu }{\mathfrak{u}}

\newcommand{\ndC }{\mathbb{C}}
\newcommand{\ndN }{\mathbb{N}}
\newcommand{\ndR }{\mathbb{R}}
\newcommand{\ndZ }{\mathbb{Z}}
\newcommand{\ndX }{\mathbb{X}}
\newcommand{\ndY }{\mathbb{Y}}

\newcommand{\ndI }{\mathbb{I}}
\newcommand{\ndIodd }{\mathbb{I}^{\rm {odd}}}
\newcommand{\ndA }{\mathbb{A}}
\newcommand{\ndD }{\mathbb{D}}
\newcommand{\ndH }{\mathbb{H}}
\newcommand{\ndE }{\mathbb{E}}
\newcommand{\ndF }{\mathbb{F}}

\newcommand{\hndI }{{\hat {\mathbb{I}}}}
\newcommand{\hndIodd }{{\hat {\mathbb{I}}}^{\rm {odd}}}
\newcommand{\hndA }{{\hat {\mathbb{A}}}}
\newcommand{\hndD }{{\hat {\mathbb{D}}}}
\newcommand{\iisor }{{\hat{\delta}}}

\newcommand{\wmfgp }{{\widetilde{\mfg }}^\prime }
\newcommand{\wmfhp }{{\widetilde{\mfh }}^\prime}
\newcommand{\wmfnpu }{{\widetilde{\mfn }}_+}
\newcommand{\wmfnmi }{{\widetilde{\mfn }}_-}

\newcommand{\mfgp }{\mfg^\prime}

\newcommand{\psl }{\mathfrak{psl}}
\newcommand{\psu }{\mathfrak{psu}}
\newcommand{\tndX }{\tilde{\mathbb{X}}}

\newcommand{\per }{\gamma }
\newcommand{\rsign }{\epsilon }
\newcommand{\sbf }[2]{(#1,#2)}
\newcommand{\sd }{\mathcal{D}}
\newcommand{\tW }{\widetilde{W}}
\newcommand{\tWext }{\widetilde{W}^\mathrm{ext}}
\newcommand{\Upm }{U^{\prime <0}}
\newcommand{\Upo }{U^{\prime \,0}}
\newcommand{\Upp }{U^{\prime >0}}
\newcommand{\Wext }{W^\mathrm{ext}}

\newcommand{\mbb}[1]{\mathbb{#1}}

\newtheorem{theorem}{Theorem}[section]
\newtheorem{proposition}[theorem]{Proposition}

\newtheorem{lemma}[theorem]{Lemma}

\theoremstyle{definition}
\newtheorem{defin}[theorem]{Definition}
\newtheorem{notation}[theorem]{Notation}

\theoremstyle{remark}
\newtheorem{remar}[theorem]{Remark}

\newcommand{\al}{{\alpha}}
\newcommand{\tr}{\triangleright}

\newcommand{\wlbr}{{[\![}}
\newcommand{\wrbr}{{]\!]}}

\newcommand{\te}{{\tilde e}}
\newcommand{\tnu }{{\tilde \nu }}
\newcommand{\ts}{{\tilde s}}
\newcommand{\ttau }{{\tilde \tau }}
\newcommand{\tomega}{{\tilde \omega}}

\newcommand{\aorigin}{{
\setlength{\unitlength}{1mm}
\begin{picture}(80,55)(10,-8)

\put(20,50){$d=4$}

\put(10,5){\circle{6}}\put(3,3){$3$}
\put(08,3){\line(1,1){4}}\put(08,7){\line(1,-1){4}}
\put(70,5){\circle{6}}\put(76,3){$1$}
\put(68,3){\line(1,1){4}}\put(68,7){\line(1,-1){4}}
\put(40,50){\circle{6}}\put(39,54.5){$2$}
\put(38,48){\line(1,1){4}}\put(38,52){\line(1,-1){4}}
\put(40,20){\circle{6}}\put(39,13){$0$}
\put(38,18){\line(1,1){4}}\put(38,22){\line(1,-1){4}}

\put(40, 23){\line(0,1){24}}\put(41, 30){$x+1$}
\put(13, 5){\line(1,0){54}}\put(35, 1){$x+1$}
\put(11, 8){\line(2,3){27}}\put(18, 29){$-x$}
\put(13, 6){\line(2,1){25}}\put(22, 15){$-1$}
\put(69, 8){\line(-2,3){27}}\put(57, 29){$-1$}
\put(67, 6){\line(-2,1){25}}\put(50, 15){$-x$}

\end{picture}}}

\newcommand{\azero}{{
\setlength{\unitlength}{1mm}
\begin{picture}(80,55)(10,-8)

\put(20,50){$d=0$}

\put(10,5){\circle{6}}\put(3,3){$3$}

\put(70,5){\circle{6}}\put(76,3){$1$}

\put(40,50){\circle{6}}
\put(39,54.5){$2$}

\put(40,20){\circle{6}}\put(39,13){$0$}
\put(38,18){\line(1,1){4}}\put(38,22){\line(1,-1){4}}

\put(40, 23){\line(0,1){24}}\put(42, 33){$-x-1$}

\put(13, 6){\line(2,1){25}}\put(22, 14){$1$}

\put(67, 6){\line(-2,1){25}}\put(55, 14){$x$}

\end{picture}
}}

\newcommand{\aone}{{
\setlength{\unitlength}{1mm}
\begin{picture}(80,55)(10,-8)

\put(20,50){$d=1$}

\put(10,5){\circle{6}}\put(3,3){$2$}

\put(70,5){\circle{6}}\put(76,3){$0$}

\put(40,50){\circle{6}}
\put(39,54.5){$3$}

\put(40,20){\circle{6}}\put(39,13){$1$}
\put(38,18){\line(1,1){4}}\put(38,22){\line(1,-1){4}}

\put(40, 23){\line(0,1){24}}\put(42, 33){$-x-1$}

\put(13, 6){\line(2,1){25}}\put(22, 14){$1$}

\put(67, 6){\line(-2,1){25}}\put(55, 14){$x$}

\end{picture}
}}

\newcommand{\atwo}{{
\setlength{\unitlength}{1mm}
\begin{picture}(80,55)(10,-8)

\put(20,50){$d=2$}

\put(10,5){\circle{6}}\put(3,3){$1$}

\put(70,5){\circle{6}}\put(76,3){$3$}

\put(40,50){\circle{6}}
\put(39,54.5){$0$}

\put(40,20){\circle{6}}\put(39,13){$2$}
\put(38,18){\line(1,1){4}}\put(38,22){\line(1,-1){4}}

\put(40, 23){\line(0,1){24}}\put(42, 33){$-x-1$}

\put(13, 6){\line(2,1){25}}\put(22, 14){$1$}

\put(67, 6){\line(-2,1){25}}\put(55, 14){$x$}

\end{picture}
}}

\newcommand{\athree}{{
\setlength{\unitlength}{1mm}
\begin{picture}(80,55)(10,-8)

\put(20,50){$d=3$}

\put(10,5){\circle{6}}\put(3,3){$0$}

\put(70,5){\circle{6}}\put(76,3){$2$}

\put(40,50){\circle{6}}
\put(39,54.5){$1$}

\put(40,20){\circle{6}}\put(39,13){$3$}
\put(38,18){\line(1,1){4}}\put(38,22){\line(1,-1){4}}

\put(40, 23){\line(0,1){24}}\put(42, 33){$-x-1$}

\put(13, 6){\line(2,1){25}}\put(22, 14){$1$}

\put(67, 6){\line(-2,1){25}}\put(55, 14){$x$}

\end{picture}
}}

\newcommand{\Dynkindiagrams}{{
\setlength{\unitlength}{1mm}
\begin{picture}(210,155)(0,0)

\put(67,163){$\tau_{f_2,2}$}
\put(60,160){\vector(1,0){20}}
\put(80,158){\vector(-1,0){20}}
\put(67,154){$\tau_{f_2,0}$}

\put(67,43){$\tau_{f_2,1}$}
\put(60,40){\vector(1,0){20}}
\put(80,38){\vector(-1,0){20}}
\put(67,34){$\tau_{f_2,3}$}

\put(-12,93){$\tau_{f_3,1}$}
\put(0,85){\vector(0,1){20}}
\put(2,105){\vector(0,-1){20}}
\put(9,93){$\tau_{f_3,2}$}

\put(128,93){$\tau_{f_3,3}$}
\put(140,85){\vector(0,1){20}}
\put(142,105){\vector(0,-1){20}}
\put(149,93){$\tau_{f_3,0}$}

\put(0,120){\atwo}\put(80,120){\azero}

\put(40,100){\vector(-1,1){20}}\put(20,118){\vector(1,-1){20}}
\put(27,115){$s_{2,4}$}\put(25,101){$s_{2,2}$}
\put(40,95){\vector(-1,-1){20}}\put(20,73){\vector(1,1){20}}
\put(29,78){$s_{1,1}$}\put(25,92){$s_{1,4}$}

\put(100,100){\vector(1,1){20}}\put(120,118){\vector(-1,-1){20}}
\put(104,115){$s_{0,4}$}\put(106,101){$s_{0,0}$}
\put(100,95){\vector(1,-1){20}}\put(120,73){\vector(-1,1){20}}
\put(106,92){$s_{3,4}$}\put(106,78){$s_{3,3}$}

\put(40,60){\aorigin}
\put(0,0){\aone}\put(80,0){\athree}
\end{picture}
}}

\newcommand{\finitedimDynkindiagrams}{{
\setlength{\unitlength}{1mm}
\begin{picture}(70,15)(0,0)

\put(10,5){\circle{6}}\put(9,10){$1$}
\put(40,5){\circle{6}}\put(39,10){$2$}
\put(70,5){\circle{6}}\put(69,10){$3$}

\put(38,3){\line(1,1){4}}\put(38,7){\line(1,-1){4}}

\put(13, 5){\line(1,0){24}}\put(24, 8){$1$}

\put(43, 5){\line(1,0){24}}\put(54, 8){$x$}

\end{picture}
}}
 
\title{Drinfeld second realization of the\\
quantum affine  superalgebras of $D^{(1)}(2,1;x)$\\
via the Weyl groupoid
}

\author{Istv\'an Heckenberger,\quad Fabian Spill, \\
Alessandro Torrielli\quad and \quad Hiroyuki Yamane}
\providecommand{\bysame}{\leavevmode\hbox to3em{\hrulefill}\thinspace}
\providecommand{\MR}{\relax\ifhmode\unskip\space\fi MR }
\providecommand{\href}[2]{#2}

\date{}

\begin{document}

\maketitle

%In this note, we use an argument 
%similar to that in  \cite{beck1}~hep-th/9404165.

%Addresses: 1: Mathematisches Institut,
%Universit\"at Leipzig,
%Pf 100920, D-04009 Leipzig,
%Germany

%2:
%Humboldt--Universit\"at zu Berlin,
%Institut f\"ur Physik,
%Newtonstra\ss e 15,
%D-12489 Berlin, Germany

%3:
%Center for Theoretical Physics,
%Laboratory for Nuclear Sciences,
%and,
%Department of Physics,
%Massachusetts Institute of Technology,
%Cambridge, Massachusetts 02139, USA

%4:
%Department of Pure and
%Applied Mathematics,
%Graduate School of 
%Information Science and Technology,
%Osaka University,
%Toyonaka, Osaka 
%560-0043, Japan
%\newline\par
%E-mails: 1: Istvan.Heckenberger@math.uni-leipzig.de \par
%2: spill@physik.hu-berlin.de \par
%3: torriell@mit.edu \par
%4: yamane@ist.osaka-u.ac.jp
%\newline\par
%Preprint number(s): MIT-CTP 3835 

\begin{abstract}
We obtain Drinfeld second realization
of the quantum affine superalgebras associated with
the affine Lie superalgebra $D^{(1)}(2,1;x)$. 
%We use a Weyl groupoid, 
%inspired by a method developed by Beck 
%using the (extended) affine
%Weyl groups.
Our results are analogous to those obtained by
Beck for the quantum affine algebras. Beck's analysis uses heavily
the (extended) affine Weyl groups of the affine Lie algebras.
In our approach the structures are based on a Weyl groupoid. 
\end{abstract}

\centerline{\small Preprint 
numbers: MIT-CTP 3835, 
%(Dept. of Physics, MIT),
HU-EP-07/15 }

\section{Introduction}

In this paper we study the quantum deformation of the affine
Lie superalgebra $D^{(1)}(2,1;x)$, where $x\in\ndC\setminus\{0,-1\}$.
%Our first main result is
%Theorem~\ref{th:HopfTri}, in which
In Definition~\ref{definition:defU'd},
for any $q\in\ndC$ 
%which is not a root of $1$,
%we give the defining relations of 
such that $q(q^n-1)(q^{nx}-1)(q^{n(x+1)}-1)\ne 0$ for all $n\in\ndN$,
we define
the quantized enveloping algebras
of $D^{(1)}(2,1;x)$ by the defining relations (cf. \cite{yamane1}
\footnote{They were originally given in \cite[Remark~7.1.1]{yamane1}
(or Prop.6.3.1(vii),(viii) in q-alg/9603015).}) in terms of the Chevalley-Serre generators.
In Theorem~\ref{th:exisT}
we attach to any simple (even \textit{and} odd) reflection 
(cf. \cite{S}) a
Lusztig type isomorphism between two such algebras.
These isomorphisms satisfy Coxeter type relations, see Theorem~\ref{th:Tcoxrel}.
In Definition~\ref{definition:defUD} we give Drinfeld second realization of
quantum $D^{(1)}(2,1;x)$. Our main result is 
Theorem~\ref{theorem:iso}, see also
Theorems~\ref{theorem:prDs} and \ref{theorem:fin}, where
we show that the two realizations are isomorphic as algebras.
See \cite{Drinfeld} for the original Drinfeld second realization
of the quantum affine algebras. The argument in this paper was inspired by 
Beck's work \cite{beck1} and we utilize the Weyl groupoid
instead of the Weyl group.
Khoroshkin and Tolstoy \cite{KohTol} obtained results concerning
quantum affine superalgebras
relevant to 
this paper. %without using the Weyl groupoids.

Our work was motivated by recent results in Hopf algebra theory and
in theoretical physics, in particular the AdS/CFT correspondence.
We sketch those aspects of these developments which are relevant
for our work.

The Lie superalgebra $D(2,1;x)$ has a very interesting relation to
$A(1,1)=\psl(2|2)$, which is the only classical basic Lie
superalgebra allowing for a nontrivial universal central extension
\cite{IK} with three central elements
(see Section~\ref{sec:sectiontwo}
for $D(2,1;x)$ and $A(1,1)$).
One can obtain this centrally extended algebra
$\psl(2|2)\circledplus \mbb{C}^3$ by a contraction from $D(2,1;x)$
in the limit $x\rightarrow - 1$.
The Lie superalgebra $\psl(2|2)$ and its central extensions
have recently become important in the context of the
AdS/CFT correspondence \cite{Maldacena}
(for comprehensive reviews, the reader is referred to \cite{rev}).
This conjecture relates the maximal supersymmetric Yang-Mills theory
in four dimensions to string theory formulated on $\mathrm{AdS}_5 \times S_5$.
On the gauge theory side of this correspondence one can think
of a certain class of operators as integrable spin chains,
and apply the Bethe ansatz technique to calculate their energy
spectrum \cite{M,BS}. The symmetry algebra of the gauge theory,
which is the superconformal algebra $\psu(2,2|4)$,
is reduced to $\mfu(1)\circledplus (\psu(2|2)\times \psu(2|2)) \circledplus \mfu(1)$
upon choosing an appropriate vacuum for the spin chain.
Excitations transform under
$\mfu(1)\circledplus (\psu(2|2)\times \psu(2|2)) \circledplus \mfu(1)$,
and the $S$-matrix which intertwines two modules is physically
interpreted as the scattering matrix of those excitations
(detailed descriptions are contained in \cite{BB}).
Interestingly, the $S$-matrix is already fixed,
up to a scalar prefactor, by vanishing of its commutators
with the generators of the centrally extended
$(\psu(2|2)\times \psu(2|2)) \circledplus \mbb{C}^3$ algebra \cite{B1,B2},
when one twists the universal enveloping algebra with an additional
braiding element \cite{P}. 
The complete symmetry algebra has been recently related to a
twisted Yangian \cite{B3}. The spectral parameter of the Yangian,
the eigenvalues of the central charges and the braiding are all
linked on the fundamental evaluation representation.

Due to its close relation to $\psl (2|2)$ it is very promising
to study the affine Lie superalgebra $D^{(1)}(2,1;x)$
(see Section~\ref{sec:sectiontwo}
for $D^{(1)}(2,1;x)$). Since one
can obtain Yangians from quantum affine algebras one can consider
physical models with quantum $D^{(1)}(2,1;x)$ symmetry as deformations
of models with Yangian $\psl (2|2)$ symmetry.
In this paper we do the first steps by deriving Drinfeld's second
realization of quantum $D^{(1)}(2,1;x)$, which we need for
further investigations of finite dimensional representations
and studies of the universal $R$-matrix. Our key tool 
is the Weyl groupoid of (quantum) $D^{(1)}(2,1;x)$. The notion
of the Weyl groupoids
was initiated and 
has intensively been studied by the first author \cite{a-Heck06a}
in order to classify
Nichols algebras of diagonal type with a finite set of
Poincar\'e-Birkhoff-Witt generators. The interest in Nichols algebras
arose with a fundamental paper of Andruskiewitsch and Schneider
\cite{a-AndrSchn98} where they developed a method to classify pointed
Hopf algebras. The results of many papers culminated in a fairly
general classification result \cite{a-AndrSchn05p}
on finite dimensional pointed Hopf algebras with abelian coradical
over the complex numbers. In the heart of the theory the Weyl groupoid
seems to play one of the fundamental roles. Guided by this observation
the first and fourth authors started to investigate
the Weyl groupoids in more detail, and obtained a Matsumoto-type theorem
\cite{HecYam} for them.
The fourth author \cite{yamane1} essentially used the Weyl groupoids
to get Serre-type defining relations of the quantum affine 
superalgebras
and the Drinfeld second realization of the quantum 
%$\psl(m|n)^{(1)}$. 
$A^{(1)}(m,n)$ (see Remark~\ref{remar:typeasuper}
for the notation $A^{(1)}(m,n)$).
The fourth author \cite{yamane1}
utilized the quantum deformation of the universal central extension
of 
%$\psl (2|2)^{(1)}$ 
%$A^{(1)}(1,1)$
$[A^{(1)}(1,1),A^{(1)}(1,1)]$
to get a new $R$-matrix.

In this paper we use the following notation. 
Let $\ndZ $ and $\ndN $ denote the sets of integers
and positive integers, respectively, and
let $\ndR $ and $\ndC $ denote 
the fields of real and complex numbers, respectively.
The symbol $\delta_{ij}$, or $\delta_{i,j}$, denotes Kronecker's $\delta $, that is,
$\delta_{ij}=1$ if $i=j$, and $\delta_{ij}=0$ otherwise.

\section{The simple Lie superalgebra $D(2,1;x)$
and the affine Lie superalgebra 
$D^{(1)}(2,1;x)$ ($x\ne 0,-1$)}
\label{sec:sectiontwo}

As for the terminology concerning
affine Lie superalgebras, we refer to
\cite{kac1}, or to
\cite{IK}, \cite{vandeLeur}.

Let $\mathfrak{v}=\mathfrak{v}(0)\oplus\mathfrak{v}(1)$
be a $\ndZ /2\ndZ $-graded
$\ndC $-linear space. If $i\in \{0,1\}$ and $j\in \ndZ $ such that
$j-i\in 2\ndZ $ then let $\mathfrak{v}(j)=\mathfrak{v}(i)$.
If $X\in\mathfrak{v}(0)$ (resp.\ $X\in\mathfrak{v}(1)$)
then we write 
\begin{equation}\label{eqn:pxzeroone}
\mbox{$\dg (X)=0$ (resp.\ $\dg (X)=1$)}
\end{equation}
and
we say that $X$ is an {\it{even}} (resp.\ {\it{odd}})
element.
If $X\in\mathfrak{v}(0)\cup\mathfrak{v}(1)$, then we say
that $X$ is a
{\it {homogeneous}} element and that
$\dg (X)$ is the {\it{parity}} (or {\it{degree}}) of $X$.
If $\mathfrak{w}\subset \mathfrak{v}$ is a subspace and
$\mathfrak{w}=(\mathfrak{w}\cap\mathfrak{v}(0))\oplus
(\mathfrak{w}\cap\mathfrak{v}(1))$
(resp. $\mathfrak{w}\subset\mathfrak{v}(0)$,
resp. $\mathfrak{w}\subset\mathfrak{v}(1)$), then we say that
$\mathfrak{w}$ is a {\it{graded}}
(resp. {\it{even}}, resp. {\it{odd}})
subspace.

Let $\mathfrak{a}=\mathfrak{a}(0)\oplus\mathfrak{a}(1)$ be a
$\ndZ /2\ndZ $-graded
$\ndC $-linear space equipped with a bilinear
map $[\,,\,]:\mathfrak{a}\times\mathfrak{a}\rightarrow\mathfrak{a} $
such that $[\mathfrak{a}(i),\mathfrak{a}(j)]\subset \mathfrak{a}(i+j)$
($i$, $j\in\ndZ $); we
recall from the above paragraph that
\begin{equation}\label{eqn:defohai}
\mathfrak{a}(i)=\{X\in\mathfrak{a}\,|\,\dg (X)=i\}. 
\end{equation}
We say that $\mathfrak{a}=(\mathfrak{a},[\,,\,])$ is a
($\ndC $-){\it{Lie superalgebra}} if
for all homogeneous elements $X$, $Y$, $Z$ of $\mathfrak{a}$ the following equations
hold.
\begin{align*}
	[Y,X]=&\,-(-1)^{\dg (X)\dg (Y)}[X,Y],& &\text{(skew-symmetry)} \\
	[X,[Y,Z]]=&\,[[X,Y],Z]+(-1)^{\dg (X)\dg (Y)}[Y,[X,Z]]. & &
	\text{(Jacobi identity)}
\end{align*} 
%For $X\in\mathfrak{a}$, we define ${\rm{ad}}X\in{\rm{End}}_{\ndC }(\mathfrak{a})$
%by $({\rm{ad}}X)(Y)=[X,Y]$.

Let $\mathfrak{a}$ be a Lie superalgebra.
We say that a bilinear form
$(\cdot ,\cdot ):\mathfrak{a}\times\mathfrak{a}
\rightarrow\ndC $ is a
{\it{supersymmetric invariant form}} on $\mathfrak{a}$
if for all homogeneous elements $X$, $Y$, $Z$ of $\mathfrak{a}$ one has
\begin{equation}\nonumber
(Y,X)=(-1)^{\dg (X)\dg (Y)}(X,Y)\,\,\mbox{and}\,\,
(X,[Y,Z])=([X,Y],Z).
\end{equation}
A graded subspace $\mathfrak{i}$ of $\mathfrak{a}$
is called an {\it{ideal}} if one has 
$[X,Y]\in\mathfrak{i}$
for all homogeneous elements $X$ of $\mathfrak{a}$ 
and all homogeneous elements $Y$ of $\mathfrak{i}$.

Let $\ndI$ be a finite set.
Let $\ndA =(\ndA _{ij})_{i,j\in\ndI }$
be an $|\ndI|\times|\ndI|$ matrix with coefficients in $\ndC $.
Suppose that we are given an $|\ndI |\times|\ndI |$
diagonal matrix $\ndD 
=(\delta_{ij}\ndD _i)_{i,j\in\ndI }$
with $\ndD _i\in\ndC \setminus\{0\}$
satisfying the condition that 
$\ndD ^{-1}\ndA $ is a symmetric matrix,
that is, $^t(\ndD ^{-1}\ndA )
=\ndD ^{-1}\ndA $.
Let $\ndIodd$ be a subset of $\ndI$.  
Let $\wmfgp=\wmfgp(\ndA ,\ndIodd)$ be the 
$\ndC$-Lie superalgebra generated by the (homogeneous) 
elements
$\ndH _i$, $\ndE _i$, $\ndF _i$
($i\in \ndI$) with
\begin{align*}
\dg (\ndH _i)=&\,0 \qquad (i\in \ndI ),\\
\dg (\ndE _j)=&\,\dg (\ndF _j)=0 \qquad
(j\in \ndI\setminus\ndIodd),\\
\dg (\ndE _j)=&\,\dg (\ndF _j)=1 \qquad
(j\in \ndIodd).
\end{align*}
and defined by the relations 
\begin{align*}
[\ndH _i,\ndH _j]=&0,&
[\ndH _i,\ndE _j]=&\ndA _{ij}\ndE _j,&
[\ndH _i,\ndF _j]=&-\ndA _{ij}\ndF _j,&
[\ndE _i,\ndF _j]=&\delta _{ij}\ndH _i&
(i,j\in \ndI ). 
\end{align*}
Let $\wmfhp=\wmfhp (\ndA, \ndIodd )$, $\wmfnpu$ and $\wmfnmi$ be the Lie subsuperalgebras
generated by the sets $\{\ndH _i|i\in\ndI\}$,
$\{\ndE _i|i\in\ndI \}$ and 
$\{\ndF _i|i\in\ndI\}$, respectively.
Then $\{\ndH _i|i\in\ndI\}$
is a $\ndC$-basis of $\wmfhp$ and hence one has
$\dim\wmfhp=|\ndI|$. Further, one obtains the decomposition
$\wmfgp= \wmfnpu\oplus \wmfhp\oplus \wmfnmi$
as a $\ndC$-vector space. The Lie superalgebras
$\wmfnpu$ and $\wmfnmi$ are free Lie superalgebras generated by the sets
$\{\ndE _i|i\in\ndI\}$
and $\{\ndF _i|i\in\ndI\}$,
respectively.
Let ${\mathfrak{r}}_+$ (resp. ${\mathfrak{r}}_-$) be the largest ideal of $\wmfgp$ 
which is contained in $\wmfnpu$
(resp. $\wmfnmi$).
Let $\mfgp=\mfgp(\ndA ,\ndIodd)$
be the quotient Lie superalgebra $\wmfgp/({\mathfrak{r}}_+\oplus{\mathfrak{r}}_-)$.
%Let ${\mathbb{P}}:\wmfgp\to\mfgp$ be the canonical map.
%Denote ${\mathbb{P}}(\ndH _i)$, ${\mathbb{P}}(\ndE _i)$
%and ${\mathbb{P}}(\ndF _i)$ again by
%$\ndH _i$, $\ndE _i$ and $\ndF _i$,
%respectively. Let 
%${\mfh }^\prime:={\mathbb{P}}(\wmfhp)$,
%${\mfn }_+:={\mathbb{P}}(\wmfnpu)$
%and ${\mfn }_-:={\mathbb{P}}(\wmfnpu)$.
Let $\ndH _i$, $\ndE _i$, $\ndF _i$, 
${\mfh }^\prime={\mfh }^\prime (\ndA, \ndIodd )$, ${\mfn }_+$, and ${\mfn }_-$
be the images of $\ndH _i$, $\ndE _i$, $\ndF _i$, 
$\wmfhp=\wmfhp (\ndA, \ndIodd )$, $\wmfnpu$, and $\wmfnmi$, respectively,
under the canonical projection $\wmfgp\to\mfgp$. 
Then $\mfgp={\mfn }_+\oplus{\mfh }^\prime\oplus{\mfn }_-$.
Further, 
there exists a (unique) Lie superalgebra
$\mfg=\mfg(\ndA ,\ndIodd )=\mfg(\ndA ,\ndD ,\ndI ,\ndIodd )$
with the following properties.
\begin{itemize}
	\item[(i)] $\mfg$ includes $\mfgp$ as a Lie subsuperalgebra.
	\item[(ii)] There exists an even subspace ${\mfh }^{\prime\prime}=
        {\mfh }^{\prime\prime}(\ndA , \ndIodd )$ of $\mfg$
		such that $\mfg={\mfh }^{\prime\prime}\oplus\mfgp$,
		$\dim{\mfh }^{\prime\prime}=|\ndI|-{\rm{rank}}\ndA $,
%		$\dg (\ndH ^{\prime\prime})=0$ for all $\ndH ^{\prime\prime}\in
%		{\mfh }^{\prime\prime}$, 
        and
%		$\dim[{\mfh }^{\prime\prime},{\mfh }^{\prime\prime}]
%		=\dim[{\mfh }^\prime,{\mfh }^{\prime\prime}]=0$.
		$[{\mfh }^{\prime\prime},{\mfh }^{\prime\prime}]=
		[{\mfh }^\prime,{\mfh }^{\prime\prime}]=\{0\}$.
\item[(iii)] Let $\mfh =\mfh (\ndA ,\ndIodd ):={\mfh }^\prime\oplus{\mfh }^{\prime\prime}$,
so $\mfg ={\mfn }_+\oplus{\mfh }\oplus{\mfn }_-$
as a $\ndC$-vector space.
Then for each $i\in\ndI$
there exists $\al_i\in{\mfh }^*$
such that $[\ndH ,\ndE _i]
=\al_i(\ndH )\ndE _i$ and
$[\ndH ,\ndF _i]
=-\al_i(\ndH )\ndF _i$
for all $\ndH \in{\mfh }$. Further,
$\al_i$ ($i\in \ndI $) are linearly independent
elements of ${\mfh }^*$.
\end{itemize}
%We call $\ndI$ the {\it {index set}} of $\mfg$.

For $\beta\in{\mfh }^*$, 
let $\mfg_\beta =\mfg(\ndA ,\ndIodd )_\beta
:=\{X\in \mfg|[\ndH ,X]=\beta(\ndH)X\,\mbox{for all $\ndH\in\mfh $}\}$.
Let $\Phi=\Phi(\ndA ,\ndIodd ):=\{\beta\in{\mfh }^*\setminus\{0\}|\dim \mfg_\beta \ne 0\}$.
The set $\Phi $ is called the {\it {root system}} of $\mfg $ and
the elements of $\Phi$ are called {\it {roots}}.
For $\al\in\Phi$ the space $\mfg_\al $ is called the {\it {root space}}
of $\al$.
Note that one obtains the decomposition
$\mfg = \mfh \oplus (\oplus_{\al\in\Phi}\mfg_\al )$
as a $\ndC$-vector space.

Note that $\mfgp=[\mfg,\mfg]$.

It is well-known that
there exists a (unique) nondegenerate supersymmetric invariant form
$(\cdot|\cdot)$ on $\mfg$ 
such that 
\begin{align*}
	(\ndH _i|\ndH )=\ndD _i\al_i(\ndH )\quad &
	\text{for all $i\in\ndI$, $\ndH \in\mfh $,}\\
	(\ndH ^{\prime\prime}_1|\ndH ^{\prime\prime}_2)=0\quad &
	\text{for all $\ndH ^{\prime\prime}_1,\ndH ^{\prime\prime}_2
	\in{\mfh }^{\prime\prime}$,}\\
	(\ndE _i|\ndF _j)=\delta_{ij}\ndD _i\quad &
	\text{for all $i,j\in\ndI$.}
\end{align*}

Assume that $\mfg$ is a finite dimensional simple Lie superalgebra.
It is well-known that $\mfg=\mfgp$ (i.e., $\mfh={\mfh }^\prime$),
and $\dim\mfg_\al=1$ for all
$\al\in\Phi$. The (non-twisted) {\it{affine Lie superalgebra}}
${\hat{\mfg }}={\hat{\mfg }}(\ndA ,\ndIodd )$ is the $\ndZ
/2\ndZ $-graded vector space
$$\mfg\otimes\ndC [t,t^{-1}]\oplus\ndC H_\iisor\oplus\ndC
H_{\Lambda_0}$$
such that $\dg (X\otimes t^m)=\dg (X)$ for all 
homogeneous $X\in \mfg $ and $m\in \ndZ $ and
$\dg (H_\iisor )=\dg (H_{\Lambda_0})=0$
(that is, $\hat{\mfg }(0)= \mfg (0)\otimes\ndC [t,t^{-1}]\oplus\ndC H_\iisor\oplus\ndC
H_{\Lambda_0}$ and $\hat{\mfg }(1)= \mfg (1)\otimes\ndC [t,t^{-1}]$), 
together with the super-bracket
\begin{eqnarray}\nonumber
\lefteqn{[X\otimes t^m+a_1H_\iisor+b_1H_{\Lambda_0},
Y\otimes t^n+a_2H_\iisor +b_2H_{\Lambda_0}]}
\\
& &=[X,Y]\otimes t^{m+n}+m\delta_{m+n,0}(X|Y)H_\iisor +
b_1nY\otimes t^n-b_2mX\otimes t^m \nonumber
\end{eqnarray}
for all $m,n\in\ndZ$, $a_1,a_2,b_1,b_2\in\ndC $
and homogeneous elements $X$, $Y$ of $\mfg$.
Note that $[{\hat{\mfg }},{\hat{\mfg }}]=\mfg\otimes\ndC
[t,t^{-1}]\oplus\ndC H_{{\hat{\delta}}}$.

The affine Lie superalgebra ${\hat{\mfg }}$ is
identified with an infinite dimensional contragredient Lie superalgebra
$\mfg (\hndA, \hndIodd )=\mfg (\hndA ,\hndD ,\hndI ,\hndIodd )$ 
%with an index set $\hndI$ 
in the following way. 
Let $\theta$ be the (unique) highest element of
$\Phi (\ndA, \ndIodd )$, that is,
$\theta\in\Phi (\ndA, \ndIodd )$ and
$\theta+\al _i\notin \Phi (\ndA, \ndIodd )$
for all $i\in\ndI $.
Let $\ndE _\theta\in\mfg (\ndA, \ndIodd ) _\theta\setminus\{0\}$
and $\ndE _{-\theta}\in\mfg (\ndA, \ndIodd ) _{-\theta}\setminus\{0\}$.
Then $\ndE _{\pm\theta}$ are homogeneous elements of $\mfg (\ndA, \ndIodd )$.
Further, one has
$\dg (\ndE _\theta)=\dg (\ndE _{-\theta})$,
$[\ndE _{-\theta},\ndE _\theta ]\in
\mfh (\ndA, \ndIodd )$,
and $(\ndE _{-\theta}|\ndE _\theta)\ne 0$.
Let $\ndH _o^\prime:=[\ndE _{-\theta},\ndE _\theta ]$.
Then $\mfg (\hndA, \hndIodd )=\mfg (\hndA ,\hndD ,\hndI ,\hndIodd )$ 
is the contragredient 
Lie superalgebra
defined with $\hndA$, $\hndD$, $\hndI$, and $\hndIodd$ below.
\begin{itemize}
	\item[(i)] 
    $\hndI$ is a set given by adding an element $o$ to 
    $\ndI$, that is, $\hndI=\ndI \cup\{o\}$
    and $|\hndI|=|\ndI| +1$.
    \item[(ii)] $\hndIodd$ is the subset of $\hndI$ defined as follows.
    If $\dg (\ndE _\theta)=0$, then let $\hndIodd:=\ndIodd$.
    If $\dg (\ndE _\theta)=1$, then let $\hndIodd:=\ndIodd \cup\{o\}$.
    \item[(iii)] $\hndD=(\delta_{ij}\hndD_i)_{i.j\in\hndI}$
    is the $|\hndI|\times |\hndI|$ diagonal matrix defined by
    $\hndD _i=\ndD _i$ ($i\in\ndI $) and 
    $\hndD _o=(\ndE _{-\theta}|\ndE _\theta)$.
    \item[(iv)] $\hndA=(\hndA_{ij})_{i.j\in\hndI}$
    is the $|\hndI|\times |\hndI|$ matrix defined by
    $\hndA _{ij}=\ndA _{ij}$,
    $\hndA _{oj}=\hndD _j^{-1}(\ndH _o^\prime |\ndH _j)$,
    $\hndA _{io}=\hndD _o^{-1}(\ndH _i |\ndH _o^\prime )$,
    and
    $\hndA _{oo}=\hndD _o^{-1}(\ndH _o^\prime |\ndH _o^\prime)$
    for $i,j\in\ndI $.
    (Note that $^t(\hndD ^{-1}\hndA )=\hndD ^{-1}\hndA $.)
\end{itemize}

More precisely, there exists an isomorphism
$\varphi:\mfg (\hndA, \hndIodd )\to\hat{\mfg }$
such that 
$\varphi (\ndH _i)=\ndH _i\otimes 1$,
$\varphi (\ndE _i)=\ndE _i\otimes 1$,
$\varphi (\ndF _i)=\ndF _i\otimes 1$ ($i\in\ndI $),
$\varphi (\ndH _o)=\ndH _o^\prime\otimes 1+\hndD _o H_\iisor $,
$\varphi (\ndE _o)=\ndE _{-\theta}\otimes t$,
$\varphi (\ndF _o)=\ndE _\theta\otimes t ^{-1}$,
and 
%$\varphi ({\mfh }^{\prime\prime}(\hndA , \hndIodd ))
%=\ndC H _{\Lambda _0}$
$\varphi ({\mfh }^{\prime\prime}(\hndA , \hndIodd ))
=\ndC (H _{\Lambda _0}-{\frac 1 2}(X|X)H_\iisor +X)$
for some $X\in{\mfh }^\prime (\ndA, \ndIodd )$.

Let ${\dot{\theta}}$, $\iisor $, and $\Lambda _0$
be
the elements
of
$\mfh (\hndA ,\hndIodd )^*$ defined by
${\dot{\theta}}(\varphi^{-1} ( H_\iisor ))
={\dot{\theta}}(\varphi^{-1} ( H_{\Lambda_0} )) =0$,
$\iisor (\varphi^{-1} ( H_\iisor ))
=\Lambda _0(\varphi^{-1} ( H_{\Lambda _0} ))=0$,
$\iisor (\varphi^{-1} ( H_{\Lambda _0} ))
=\Lambda _0(\varphi^{-1} ( H_\iisor ))$ $=1$,
${\dot{\theta}}(\varphi^{-1} (\ndH\otimes 1))=\theta(\ndH )$, and
$\iisor (\varphi^{-1} (\ndH\otimes 1))=\Lambda _0(\varphi^{-1} (\ndH\otimes 1))=0$
for all $\ndH\in\mfh (\ndA ,\ndIodd )$.
Then one has $\al_o=\iisor - {\dot{\theta}}$
and, moreover, 
$( \varphi^{-1} ( H_\iisor ) | \ndH ) 
=\iisor (\ndH)$ and 
$( \varphi^{-1} ( H_{\Lambda_0} ) | \ndH ) 
=\Lambda_0 (\ndH)$
for all $\ndH\in\mfh (\hndA ,\hndIodd )$.

Now we define the Lie superalgebra
$D(2,1;x)$ and the affine Lie superalgebra
$D^{(1)}(2,1;x)$.
Let $\ndI=\{1,2,3\}\subset\ndN $ and $\ndIodd =\{2\}$.
Let $x\in\ndC\setminus\{0\}$ and
\begin{equation}\label{eq:thdimAD}
\ndA =
\left(
\begin{array}{ccc}
2 & -1 & 0 \\
1 & 0 & x \\
0 & -1 & 2 
\end{array}
\right),\quad
\ndD =
\left(
\begin{array}{ccc}
-1 & 0 & 0 \\
0 & 1 & 0 \\
0 & 0 & -x^{-1}
\end{array}
\right).
\end{equation}
First consider the case $x\not=-1$.
Then $\mfg=\mfgp$ and $\dim\mfg=17$. and, moreover,
$\Phi=\{\pm\al _1, \pm\al _2, \pm\al _3, 
\pm(\al _1+\al _2), \pm(\al _2+\al _3), 
\pm(\al _1+\al _2+\al _3), 
\pm(\al _1+2\al _2+\al _3) \}$. 
Further, $\mfg$ is a finite dimensional
simple Lie superalgebra and is called $D(2,1;x)$.
The affine Lie superalgebra ${\hat{\mfg }}$ is denoted
by $D^{(1)}(2,1;x)$. As mentioned above, $D^{(1)}(2,1;x)$
is identified with a contragredient Lie superalgebra. Its
Dynkin diagram is given in Figure~\ref{fig:Dynkin}
(see the one labeled $d
=2$ in Figure~\ref{fig:Dynkin}
especially).

\begin{figure}
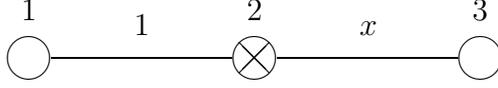

\begin{center}
\finitedimDynkindiagrams
\end{center}
\caption{(Standard) Dynkin diagram of $D(2,1;x)$
($x\ne 0,-1$)}
\label{fig:finDynkin}
\end{figure}

Now assume that $x=-1$. Then $\dim\mfg=16$, $\dim \mfg^\prime =15$ and 
$\mfg ,\mfg ^\prime $ are called $\mathfrak{gl}(2|2)$ and
${\mathfrak{sl}}(2|2)$, respectively. 
Further, $\Phi=\{\pm\al _1, \pm\al _2, \pm\al _3, 
\pm(\al _1+\al _2), \pm(\al _2+\al _3), 
\pm(\al _1+\al _2+\al _3) \}$. 
However
${\mathfrak{sl}}(2|2)$ is not simple, and
${\psl(2|2)}:=
A(1,1)
%A(1|1)
:=
{\mathfrak{sl}}(2|2) /\ndC(\ndH _1 +2\ndH _2
+\ndH _3)$ is a $14$-dimensional simple Lie superalgebra.
We obtain a 17-dimensional 
Lie superalgebra from  $D(2,1;x)$ 
by performing a specialization of $x$ at $-1$.
It is a universal central extension of ${\psl(2|2)}$
and ${\mathfrak{sl}}(2|2)$.
Similarly, we obtain 
a universal central extension of ${\psl(2|2)}\otimes\ndC [t,t^{-1}]$
and ${\mathfrak{sl}}(2|2)\otimes\ndC [t,t^{-1}]$
from $[D^{(1)}(2,1;x),D^{(1)}(2,1;x)]$
by performing a specialization of $x$ at $-1$, see \cite{IK}.

\begin{remar}\label{remar:typeasuper} {\it{The Lie superalgebras $\mathfrak{gl}(m+1|n+1)$,
$\mathfrak{sl}(m+1|n+1)$, $\psl(n+1|n+1)$, $A(m,n)$
and $A^{(1)}(m,n)$.}}
Let $m$ and $n$ be non-negative integers such that $m+n\geq 1$.
For $i,j\in\{1,\ldots,m+n+2\}$, let ${\bf{E}}_{i,j}$ denote the $(m+n+2)\times (m+n+2)$ matrix
having 1 in $(i,j)$ position and 0 otherwise, that is,
the $(i,j)$-matrix unit. 
Let ${\bf{E}}_{m+n+2}$ denote the $(m+n+2)\times (m+n+2)$ unit matrix, that is,
$\sum_{i=1}^{m+n+2}{\bf{E}}_{i,i}$.
Denote by ${\rm{M}}_{m+n+2}(\ndC )$ 
the $\ndC$-linear space of the $(m+n+2)\times (m+n+2)$-matrices,
that is, $\oplus_{i,j=1}^{m+n+2}\ndC {\bf{E}}_{i,j}$.
The  
Lie superalgebra $\mathfrak{gl}(m+1|n+1)$ is defined by
$\mathfrak{gl}(m+1|n+1)={\rm{M}}_{m+n+2}(\ndC )$
(as a $\ndC$-linear space),
$\mathfrak{gl}(m+1|n+1)(0)=(\oplus_{i,j=1}^{m+1}\ndC {\bf{E}}_{i,j})\oplus(\oplus_{i,j=m+2}^{m+n+2}\ndC {\bf{E}}_{i,j})$,
$\mathfrak{gl}(m+1|n+1)(1)=(\oplus_{i=1}^{m+1}\oplus_{j=m+2}^{m+n+2}\ndC {\bf{E}}_{i,j})\oplus(
\oplus_{j=1}^{m+1}\oplus_{i=m+2}^{m+n+2}\ndC {\bf{E}}_{i,j})$
and 
$[X,Y]=XY-(-1)^{\dg (X)\dg (Y)}YX$
for all $X$, 
$Y$ $\in\mathfrak{gl}(m+1|n+1)(0)\cup\mathfrak{gl}(m+1|n+1)(1)$,
where $XY$ and $YX$ mean the matrix product, that is,
${\bf{E}}_{i,j}{\bf{E}}_{k,l}=\delta_{j,k}{\bf{E}}_{i,l}$.
Define the $\ndC$-linear map ${\rm{str}}:\mathfrak{gl}(m+1|n+1)\to\ndC$
by ${\rm{str}}({\bf{E}}_{i,j})=\delta_{i,j}
(\sum_{k=1}^{m+1}\delta_{i,k}-\sum_{l=m+2}^{m+n+2}\delta_{i,l})$.
The Lie subsuperalgebra $\{X\in \mathfrak{gl}(m+1|n+1)\,|\,{\rm{str}}(X)=0\}$
of $\mathfrak{gl}(m+1|n+1)$ is denoted as $\mathfrak{sl}(m+1|n+1)$.
The finite dimensional simple Lie superalgebra $A(m,n)$ (cf.~\cite{kac1}) 
is defined as follows.
Let $\mathfrak{z}$ be the one dimensional ideal 
$\ndC {\bf{E}}_{m+n+2}$ of $\mathfrak{gl}(m+1|n+1)$.
If $m\ne n$, then $A(m,n)$ means $\mathfrak{sl}(m+1|n+1)$.
On the other hand, $A(n,n)$ means $\mathfrak{sl}(n+1|n+1)/\mathfrak{z}$, and
is also denoted as
$\psl(n+1|n+1)$. 
%(cf.~\cite{GP}). 

Recall the 
Lie superalgebras
$\mfg= \mfg (
\ndA ,
\ndD , \ndI, \ndIodd)$
and 
${\hat{\mfg }}={\hat{\mfg }}(\ndA ,\ndIodd )=
\mfg\otimes\ndC [t,t^{-1}]\oplus\ndC H_\iisor\oplus\ndC
H_{\Lambda_0}$, introduced above.

Define the supersymmetric invariant form 
$(\,,\,)$ 
on $\mathfrak{gl}(m+1|n+1)$
by 
$(X,Y)={\rm{str}}(XY)$.
%Denote by $(\,,\,)_{A(n,n)}$ the supersymmetric invariant form on $A(n,n)$
%induced from $(\,,\,)$. 
Then the infinite dimensional
Lie superalgebra 
%$A^{(1)}(n,n)$ 
%(resp. $\mathfrak{sl}(m+1|n+1)^{(1)}$, 
%resp. $\mathfrak{gl}(m+1|n+1)^{(1)}$)
$\mathfrak{gl}(m+1|n+1)^{(1)}$ 
is defined 
in the same way as that for ${\hat{\mfg }}$
with 
%$A(n,n)$ 
%(resp. $\mathfrak{sl}(m+1|n+1)$,
%resp. $\mathfrak{gl}(m+1|n+1)$) 
$\mathfrak{gl}(m+1|n+1)$
and 
%$(\,,\,)_{A(n,n)}$ 
%(resp. $(\,,\,)$,
%resp. $(\,,\,)$) 
$(\,,\,)$
in place of 
${\mfg }$ and $(\,|\,)$ respectively.
% (cf.~\cite{kac1}, see also \cite{IK}).
Further, $\mathfrak{sl}(m+1|n+1)^{(1)}$ means
the Lie subsuperalgebra
$\mathfrak{sl}(m+1|n+1)\otimes\ndC [t,t^{-1}]
\oplus\ndC H_\iisor\oplus\ndC 
H_{\Lambda_0}$
of $\mathfrak{gl}(m+1|n+1)^{(1)}$.
If $m\ne n$, then $A^{(1)}(m,n)$ means $\mathfrak{sl}(m+1|n+1)^{(1)}$.
On the other hand, $A^{(1)}(n,n)$ means 
%$\mathfrak{sl}(n+1|n+1)^{(1)}/(\oplus_{r\in\ndZ}\mathfrak{z}\otimes t^r)$
$\mathfrak{sl}(n+1|n+1)^{(1)}/(\mathfrak{z}\otimes \ndC [t,t^{-1}])$.
(See also \cite{kac1}, or \cite{IK}, for these notation.)

Assume
$\ndA$, $\ndD$, $\ndI$, and $\ndIodd$
to be
the $(m+n+1)\times (m+n+1)$ matrix 
$(-\delta_{i-j,1}
+2(1-\delta_{i,m+1})\delta_{i,j}
-(-1)^{\delta_{i,m+1}}\delta_{i-j,-1}
)_{1\leq i,j\leq m+n+1}$, 
the diagonal $(m+n+1)\times (m+n+1)$ matrix
$(
\delta_{i,j}(\sum_{k=1}^{m+1}\delta_{i,k}-\sum_{l=m+2}^{m+n+1}\delta_{i,l})
)_{1\leq i,j\leq m+n+1}$, 
$\{1,\ldots,m+n+1\}$,
and 
$\{m+1\}$ respectively.
%Recall the 
%contragredient Lie superalgebra
%$\mfg= \mfg (
%\ndA ,
%\ndD , \ndI, \ndIodd)$, introduced above.

Assume that $m\ne n$.
Then we identify $A(m,n)$ with $\mfg$,
since
there exists a unique isomorphism
$\varpi :\mfg \to A(m,n)$
such that $\varpi(\ndE _i)={\bf{E}}_{i,i+1}$
and
$\varpi(\ndF _i)={\bf{E}}_{i+1,i}$.
Further, we identify  $A^{(1)}(m,n)$ with 
the affine Lie superalgebra $\hat{\mfg }$,
since 
$(\varpi (X),\varpi (Y))=(X|Y)$
for all $X$, $Y\in \mfg$.

Assume that $m=n$.
Then $\mfg$ is isomorphic to 
$\mathfrak{gl}(n+1|n+1)$,
and we identify them.
Note that 
$\mfg$ is not simple since 
$[\mfg,\mfg]\ne \mfg$.
Nevertheless, we define 
$\hndA$,
$\hndD$, $\hndI$, and $\hndIodd$ in the same way as above,
and we let ${\bar{\mfg }}:=
\mfg(\hndA , \hndD , \hndI, \hndIodd)$.
Let ${\overline{\mathfrak{sl}}}(n+1|n+1)$
be the Lie subsuperalgebra 
$\mathfrak{sl}(n+1|n+1)^{(1)}
\oplus
\ndC {\bf{E}}_{1,1}$
of $\mathfrak{gl}(n+1|n+1)^{(1)}$.
Then
${\bar{\mfg }}$ 
is isomorphic to 
$
%(\mathfrak{sl}(n+1|n+1)^{(1)}
%   \rtimes 
%+
%\ndC {\bf{E}}_{1,1})
{\overline{\mathfrak{sl}}}(n+1|n+1)
/(\oplus_{r\in\ndZ\setminus\{0\}}\mathfrak{z}\otimes t^r)$
(cf.~\cite[Section~1.5]{yamane1}).

\end{remar}

\section{Semigroups and braid semigroups}
\subsection{Semigroups}

In this section we fix notations and terminology concerning semigroups.
This will be helpful for the definition of semigroups
by generators and relations.

Let $K$ be a non-empty set. We call $K$ a 
{\it{semigroup}} if it is equipped with a product 
$K\times K\rightarrow K$, $(x,y)\mapsto xy$, 
satisfying the associativity law, that is, $(xy)z=x(yz)$
for $x$, $y$, $z\in K$. If
$K$ is a semigroup, we call it a {\it{monoid}}
if there exists a unit $1\in K$, that is, 
$1x=x1=x$ for $x\in K$. 
If
$K$ is a semigroup and does not
have a unit, let ${\widehat K}$ denote the monoid
obtained from $K$ by
adding a unit.
\par
Let $H$ be a non-empty set and $L(H)$ a set of all the finite
sequences of elements of $H$,
so,
$L(H)=\{(h_1,\ldots,h_n)|n\in\ndN , h_i\in H\}$.
We regard $L(H)$
as a semigroup whose product is defined by
$(h_1,\ldots,h_m)(h_{m+1},\ldots, h_{m+n})
=(h_1,\ldots, h_m,h_{m+1},\ldots, h_{m+n})$.
Then $L(H)$ is called a {\it{free semigroup}} on $H$
and ${\widehat {L(H)}}$ a {\it{free monoid}} on $H$.
Let $\{(x_j,y_j)|j\in J\}$ be a subset of $L(H)\times L(H)$,
where $J$ is an index set.
As usual,
for at most two elements
$g$, $g^\prime$ of $L(H)$,
we let the notation 
$\{g,g^\prime\}$ mean
the subset of $L(H)$
consisting of $g$ and $g^\prime$.
(Hence the cardinality of $\{g,g^\prime\}$
is $2-\delta_{g,g^\prime}$.
As for $\delta_{g,g^\prime}$,
see the last paragraph in Introduction.)
%two if $g^\prime_1\ne g^\prime_2$
%and one otherwise.
For $g_1$, $g_2\in L(H)$,
%For $g_1$, $g_2\in L(H)$,
we write $g_1 \sim_1 g_2$ if the equation
$\{g_1,g_2\}=\{z_1x_jz_2, z_1y_jz_2\}$
(equality of subsets of $L(H)$) %consisting of at most two elements)
holds for some $j\in J$
and some $z_1$, $z_2\in{\widehat {L(H)}}$.
For $g$, $g^\prime\in L(H)$,
we write $g \sim g^\prime$ if $g = g^\prime$ or
there exist finitely many elements $g_1,\ldots,g_r$ of $L(H)$ such that
$g_1=g$, $g_r=g^\prime$ and $g_i \sim_1 g_{i+1}$. Then $\sim$ 
is an equivalence relation in $L(H)$. Let $L(H)/\!\!\sim$ be the set 
of the equivalence classes in $L(H)$ 
with respect to $\sim$.
For $g\in L(H)$, 
let $[g]$ be the 
element of $L(H)/\!\!\sim$ containing $g$. We regard
$L(H)/\!\!\sim$ as a semigroup so that the map
$L(H)\rightarrow L(H)/\!\!\sim$ defined by $g\mapsto [g]$ is a homomorphism.
We call $L(H)/\!\!\sim$ {\it{the semigroup generated by $H$ 
and defined by the relations $x_j=y_j$}} 
($j\in J$). 
When there is not fear of misunderstanding,
we also denote $[g]$ simply by $g$.

\subsection{The Weyl groupoid of $\boldsymbol{D^{(1)}(2,1;x)}$}

For the presentation of contragredient Lie superalgebras $\mfg $
one can use different Dynkin diagrams. This fact leads to the definition
of the Weyl groupoid as a symmetry object of $\mfg $.
% which strictly contains the Weyl group of the even part of $\mfg $.
General properties of such groupoids were investigated in \cite{HecYam}.
In this section we introduce the Weyl groupoid and the
extended Weyl groupoid
of the affine Lie superalgebra $D^{(1)}(2,1;x)$.

%Let $\sd$ be the subset of $\ndZ$ defined by
%\begin{gather*}
%\sd:=\{0,1,2,3,4\}. 
%\end{gather*} 
Let $\sd:=\{0,1,2,3,4\}$.
Let $\tr $: $S_5\times \sd \to \sd $
denote the usual (left) action of the symmetric group
$S_5$ on $\sd $ by permutations. The elements of the set $\sd $ will be used to
label different Dynkin diagrams for the affine Lie superalgebra $D^{(1)}(2,1;x)$.

\begin{figure}
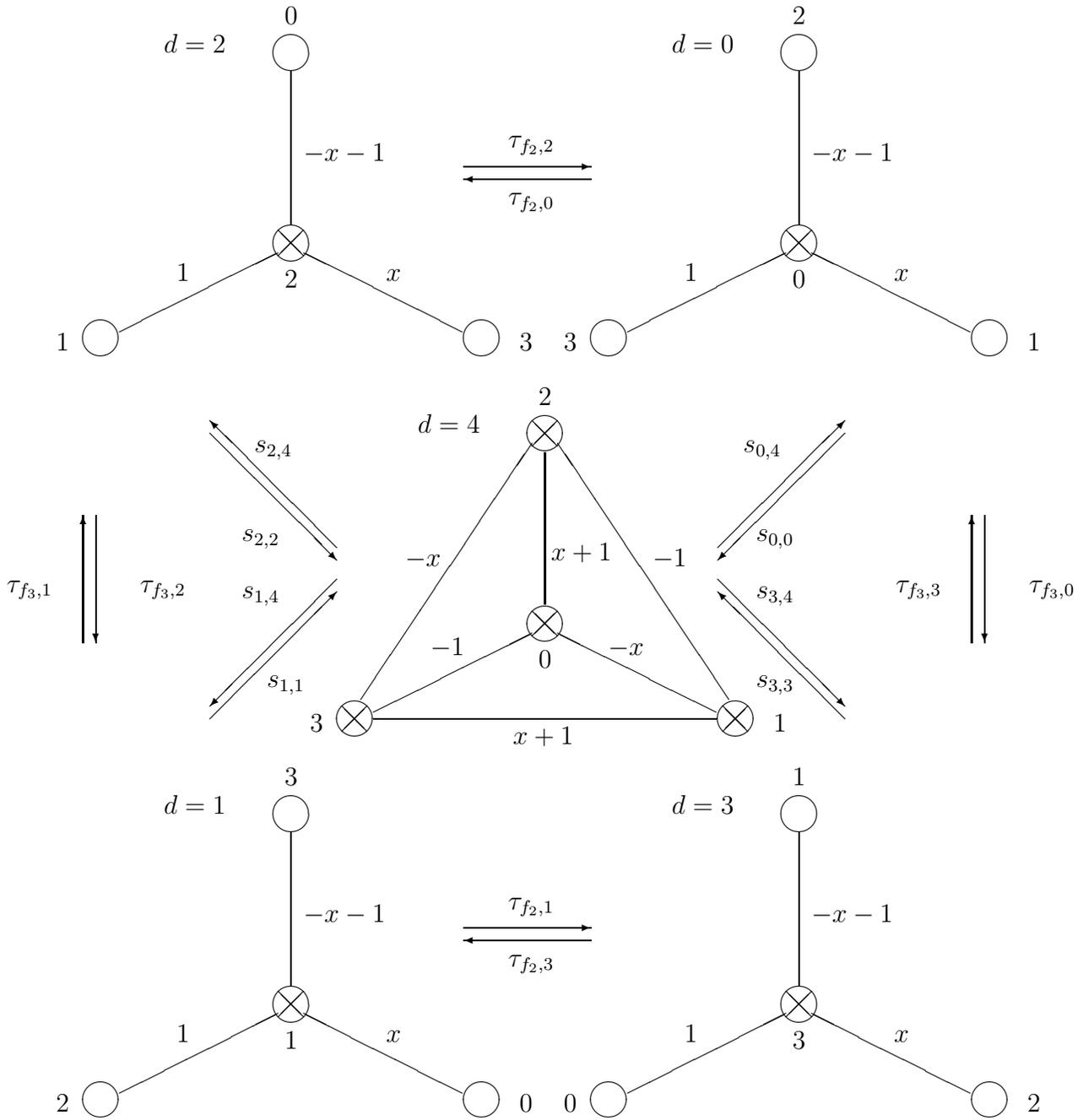

\begin{center}
\Dynkindiagrams
\end{center}
\caption{Dynkin diagrams of $D^{(1)}(2,1;x)$ ($x\ne 0,-1$)}
\label{fig:Dynkin}
\end{figure}

%Let $\isv$ be the subset of $\ndZ$ defined by
%\begin{gather*}
%\isv :=\{0,1,2,3\}.
%\end{gather*} 
Let $\isv :=\{0,1,2,3\}$.
Note that $|\isv |=4$ is the rank of $D^{(1)}(2,1;x)$.
This set will be used to label the vertices in a given Dynkin diagram.
In order to define the Weyl groupoid
we will need the following structure constants.
For $d\in \sd \setminus \{4\}$ and $i,j\in \isv $ with $i\ne j$
let 
\begin{align}
m_{i,j;d}:=&
\begin{cases}
2 & \mbox{if $i\ne d\ne j$,}  \\
3 & \mbox{otherwise},
\end{cases}&
m_{i,j;4}:=&3.
\label{eq:mijd}
\end{align}
The index $d=4$ is distinguished, see Figure~\ref{fig:Dynkin}.
Note that
$$
m_{i,j;d}=m_{j,i;d}=m_{i,j;n_i\tr d}=m_{i,j;n_j\tr d}
$$
for all $d\in \sd $ and $i,j\in \isv $,
where $n_i=(i\,4)$ as an element in $S_5$.
Further one has
\begin{align*}
n_i\tr d=&d\quad
\mbox{if $m_{i,j;d}=2$,}&
n_in_jn_i=&n_jn_in_j =(i\,j).
\end{align*}

The extended Weyl groupoid defined below contains even and odd reflections
and elements corresponding to permutations of vertices of Dynkin diagrams.
Note that any permutation of vertices of a given diagram can be identified
with a permutation $f$ of $I$. In our setting only the Klein four-group
\begin{align}\label{eq:A4}
	\Klf =\{f_0:=\id ,f_1:=(01)(23),f_2:=(02)(13),f_3:=(03)(12)\}
\end{align}
will be needed.
Further, any $f\in \Klf $ induces a permutation $\gamma (f)$
of Dynkin diagrams $d\in \sd $. Thus one obtains a group homomorphism
$$\per :\Klf \to S_5=\mathrm{Perm}\,\sd ,$$
defined by the following formula, see Figure~\ref{fig:Dynkin}.
\begin{align*}
	\per (f)\tr d=
	\begin{cases}
		f(d) & \text{if $d\in \{0,1,2,3\}$,}\\
		4 & \text{if $d=4$}
	\end{cases}
\end{align*}
for all $f\in \Klf $. This formula is accidentally true in our setting,
and can not be generalized to arbitrary contragredient Lie superalgebras.
By abuse of notation we will also write $f\tr d$
instead of $\gamma (f)\tr d$.

Let $\Wext $ be the semigroup generated by
\begin{align}
\{0\}\cup\{e_d\,|\,d\in \sd \}
\cup\{s_{i,d}\,|\,i\in \isv ,d\in \sd \}
\cup\{\tau _{f,d}\,|\,f\in \Klf , d\in \sd \}
	\label{eq:Wextgen}
\end{align}
and defined by the following relations~\eqref{eq:rel1}--\eqref{eq:rel9}:
\begin{align}
	0=&0u=u0\,\,\makebox[0pt][l]{for all elements $u$ in \eqref{eq:Wextgen},}	
	\label{eq:rel1}\\
	e_de_d=&e_d,\quad \makebox[0pt][l]{$e_de_{d'}=0$ for $d\not=d'$,}
	\label{eq:rel2}\\
	e_{n_i\tr d}s_{i,d}=&s_{i,d},\quad
	\makebox[0pt][l]{$s_{i,d}e_d=s_{i,d}$,}
	\label{eq:rel3}\\
	s_{i,n_i\tr d}s_{i,d}=&e_d,
	\label{eq:rel4}\\
	s_{i,d}s_{j,d}=&s_{j,d}s_{i,d}&
	&\text{if $m_{i,j;d}=2$,}
	\label{eq:rel5}\\
	s_{i,n_jn_i\tr d}s_{j,n_i\tr d}s_{i,d}=&
	s_{j,n_in_j\tr d}s_{i,n_j\tr d}s_{j,d}&  &\text{if $m_{i,j;d}=3$,}
	\label{eq:rel6}\\
	e_{f\tr d}\tau _{f,d}=&\tau _{f,d},\quad
	\makebox[0pt][l]{$\tau _{f,d}e_d=\tau _{f,d}$,
	\quad $\tau _{f_0,d}=e_d$,}
	\label{eq:rel7}\\
	\tau _{f,f'\tr d}\tau _{f',d}=&\tau _{ff',d}\quad
	\makebox[0pt][l]{for $f,f'\in \Klf $,}
	\label{eq:rel8}\\
	\tau _{f,n_i\tr d}s_{i,d}=&s_{f(i),f\tr d}\tau _{f,d}.
	\label{eq:rel9}
\end{align}

\begin{defin}
	The semigroup $\Wext $ is called the \textit{extended Weyl
	groupoid} of the affine Lie superalgebra
	$D^{(1)}(2,1;x)$.
	\footnote{Here we use the less standard
	terminology concerning groupoids, in which the multiplication is
	globally defined, but may be zero (instead of nondefined).
	Then all groupoids are semigroups. Removing $0$ from the semigroup
	gives a groupoid in the standard sense.}
	The subgroupoid of $\Wext $ generated by the set
\begin{align}
\{0\}\cup\{e_d\,|\,d\in \sd \}
\cup\{s_{i,d}\,|\,i\in \isv ,d\in \sd \}
	\label{eq:Wgen}
\end{align}
        is called the \textit{Weyl groupoid}
	of $D^{(1)}(2,1;x)$ and will be denoted by $W$.
\end{defin}

%For all $i$, $j\in \isv $ there exists a unique
%$f\in \Klf $ such that $f\tr i=j$ and $f\tr j=i$. We will also write
%$$\tau_{\{i,j\},d}:=\tau _{f,d}$$
%in this case.
%In particular one has $\tau_{\{i,i\},d}=e_d$ for all $i\in \isv $.
For an element
$w=s_{i_r,d_r}\cdots s_{i_2,d_2}s_{i_1,d_1}$ of $W$,
where $d_u:=n_{i_{u-1}}\cdots n_{i_1}\tr d_1$ for $1\le u\le r$, we
also use the abbreviations
\begin{align}
	w=&s_{i_r}\cdots s_{i_2}s_{i_1,d_1},&
	\tau_{f,d_r}w=&\tau _fs_{i_r}\cdots s_{i_2}s_{i_1,d_1},&
	w\tau _{f,f\tr d}=&s_{i_r}\cdots s_{i_2}s_{i_1}\tau _{f,f\tr d}.
	\label{eq:abbr}
\end{align}
If $r=0$, let $s_{i_r}\cdots s_{i_2}s_{i_1,d}:=e_d$,
$\tau_fs_{i_r}\cdots s_{i_2}s_{i_1,d}:=\tau_{f,d}$
and 
$s_{i_r}\cdots s_{i_1}\tau_{f,d}:=\tau_{f,d}$.
Note that
\begin{equation}
	\begin{aligned}
	\Wext =\{0\}\cup\{\tau_fs_{i_r}\cdots s_{i_2}s_{i_1,d}\,|\,&d\in \sd ,
	f\in \Klf ,\,r\in \ndN _0,\,i_1,\ldots ,i_r\in \isv ,\\
	&\mbox{$i_u \ne i_{u+1}$ for $1\leq u\leq r-1$}\}.
        \end{aligned}
	\label{eq:Welements}
\end{equation}

% \subsection{The map ${\widetilde {\rm {sgn}}}$}

Now we prove that the elements $\tau_fs_{i_r}\cdots s_{i_2}s_{i_1,d}$
in Eq.~\eqref{eq:Welements} are nonzero.
Let
$\ndR ^\sd $ be an $\ndR $-vector space of dimension $|\sd |$,
and let $\{v_d\,|\,d\in \sd \}$ be a fixed basis of $\ndR ^\sd $.
Then there is a unique semigroup homomorphism
$${\widetilde {\rm {sgn}}}:\Wext \rightarrow{\rm{End}}(\ndR ^\sd )$$ 
such that
\begin{equation}
\begin{aligned}
	{\widetilde {\rm {sgn}}}(0)(v_{d'})&=0,&
	{\widetilde {\rm {sgn}}}(e_d)(v_{d'})=&\delta_{dd'}v_{d'},\\
	{\widetilde {\rm {sgn}}}(s_{i,d})(v_{d'})=&(-1)\delta_{dd'}v_{n_i\tr d},&
	{\widetilde {\rm {sgn}}}(\tau _{f,d})(v_{d'})=&\delta_{dd'}v_{f\tr d'}
\end{aligned}
	\label{eq:sgntilde}
\end{equation}
for all $d,d'\in \sd $, $i\in \isv $, $f\in \Klf $.
In particular $\tau _fs_{i_r}\cdots s_{i_2}s_{i_1,d}\ne 0$ for all
$d\in \sd $, $f\in \Klf $, $r\in \ndN _0$, and $i_1,\ldots ,i_r\in \isv $.
Let
\begin{align*}
	e_d^{-1}:=&e_d,&
	(\tau _fs_{i_r}\cdots s_{i_2}s_{i_1,d})^{-1}:=s_{i_1}\cdots
	s_{i_{r-1}}s_{i_r}\tau _{f,fn_{i_r}\cdots n_{i_2}n_{i_1}\tr d}.
\end{align*}
One says that an expression $w=\tau _fs_{i_r}\cdots s_{i_2}s_{i_1,d}$ is
\textit{reduced} if $w=\tau _{f'}s_{j_s}\cdots s_{j_2}s_{j_1,d}$ implies that
$s\ge r$. In this case one defines $\ell (w):=r$.

\subsection{The braid semigroup}
\label{sec:brsgroup}

Lusztig \cite{b-Lusztig93} defined automorphisms of quantized enveloping algebras
of Kac-Moody Lie algebras
attached to all simple reflections of the corresponding Weyl group. These automorphisms
are not involutions, but nevertheless they 
satisfy some Coxeter relations. Analogously there exist isomorphisms
between different realizations of quantized $D^{(1)}(2,1;x)$, see
Section~\ref{sec:brgraction}, which also satisfy Coxeter relations.
At this place we introduce the abstract semigroup which forms a bridge between
the aforementioned isomorphisms and the Weyl groupoid of $D^{(1)}(2,1;x)$.

Let $\tWext $ be the semigroup generated by
\begin{align}
	\{0\}\cup\{\te _d\,|\,d\in \sd \}\cup\{\ts _{i,d}\,|\,i\in \isv ,d\in \sd \}
	\cup \{\ttau _{f,d}\,|\,f\in \Klf ,d\in \sd \}
	\label{eq:tWextgen}
\end{align}
and defined by the relations analogous to \eqref{eq:rel1}-\eqref{eq:rel3},
\eqref{eq:rel5}-\eqref{eq:rel9}.
Let $\tW $ be the subsemigroup of $\tWext $ generated by
\begin{align}
	\{0\}\cup\{\te _d\,|\,d\in \sd \}\cup\{\ts _{i,d}\,|\,i\in \isv ,d\in \sd \}.
	\label{eq:tWgen}
\end{align}

\begin{notation}\label{no:tabbr}
	For elements of $\tW $ and $\tWext $ we use a notation analogous to the one in
	Eq.~\eqref{eq:abbr}.
\end{notation}
Similarly to Eq.~\eqref{eq:Welements} one has
\begin{align}
	\tWext =\{0\}\cup\{\ttau_f\ts_{i_r}\cdots \ts_{i_2}\ts_{i_1,d}\,|\,&d\in \sd ,
	f\in \Klf ,\,r\in \ndN _0,\,i_1,\ldots ,i_r\in \isv \}.
	\label{eq:tWelements}
\end{align}
Note that there exists a canonical semigroup homomorphism
$\mfp :\tWext \rightarrow \Wext $ such that
\begin{align*}
	\mfp(0)=&0, & 
	\mfp(\ttau_f\ts_{i_r}\cdots \ts_{i_2}\ts_{i_1,d})=&
	\tau_fs_{i_r}\cdots s_{i_2}s_{i_1,d}
\end{align*}
for all $d\in \sd $, $f\in \Klf $, $r\in \ndN _0$, and $i_1,\ldots ,i_r\in \isv $.

\subsection{Special elements of $\boldsymbol{\tWext }$}
\label{ss:specelem}

The extended Weyl group of an affine Lie algebra $\hat{\mfg }$,
which is the group generated by the Weyl group and by diagram automorphisms
of $\hat{\mfg }$,
can be written as the semidirect
product of the (finite) Weyl group of $\mfg $ and a free abelian group,
and the latter can be identified with the weight lattice corresponding
to $\mfg $.
We expect that
a similar decomposition holds for the extended Weyl groupoid of any
affine Lie superalgebra. Since we are mainly interested in the (quantum)
affine Lie superalgebra $D^{(1)}(2,1;x)$, we will not work out here the
details of such a decomposition, but concentrate on those formulas
which are necessary to obtain
Theorems~\ref{th:exisT}, \ref{th:Tcoxrel}, and \ref{th:tzerotone}.
Nevertheless it may be helpful
to think about the elements $\omega ^\vee _{i,d}$ introduced below
as the generators of the weight lattice in the extended
Weyl groupoid $\Wext $.
 Further, the analog of the Weyl group of $\mfg $ will be
the subgroupoid of $\Wext $ generated by the reflections $s_{i,d}$,
where $i\in \isv \setminus \{0\}$, $d\in \sd \setminus \{0\}$.

%$\tomega^\vee_{i,a(j)}=\tomega^\vee_{\al_{i,a(j)}}
%=\tomega^\vee (i:a(j)\rightarrow a(j))\in\tW$ 
%with $i\in\{1,2,3\}$ and $j\in\{1,2,3,4\}$
  
Define $\omega^\vee_{i,d}\in\Wext $, where
  $i\in \isv \setminus \{0\}$ and $d\in \sd \setminus \{0\}$, as follows.
Let $i$, $j$, $k\in \isv \setminus \{0\}$ be such that
$\{i,j,k\}= \{1,2,3\}(=\isv \setminus \{0\})$.
Let
\begin{equation}
\begin{aligned}
	\omega^\vee_{i,4}:=&\tau _{f_i,4}s_{i,i}s_{k,i}s_{j,i}s_{i,4}
	=\tau _{f_i}s_is_ks_js_{i,4}
	=\tau _{f_i}s_is_js_ks_{i,4},\\
	\omega^\vee_{j,i}:=&\tau _{f_j,k}s_{j,k}s_{k,4}s_{i,i}s_{j,i}
	=\tau _{f_j}s_js_ks_is_{j,i},\\
	\omega^\vee_{i,i}:=&s_{0,i}s_{i,4}s_{j,j}s_{k,j}s_{j,4}s_{i,i}
	=s_0s_is_js_ks_js_{i,i}=s_0s_is_ks_js_ks_{i,i}. 
\end{aligned}
	\label{eq:omega}
\end{equation}
It is easy to check \cite{HecYam}
that all of the above expressions are reduced.
Define also $\tomega^\vee_{i,d}\in\tWext $, where
$i\in \isv \setminus \{0\}$ and $d\in \sd \setminus \{0\}$ by the
following formulas:
\begin{equation}
\begin{aligned}
	\tomega^\vee_{i,4}:=&\ttau _{f_i,4}\ts_{i,i}\ts_{k,i}\ts_{j,i}\ts_{i,4}
	=\ttau _{f_i}\ts_i\ts_k\ts_j\ts_{i,4}
	=\ttau _{f_i}\ts_i\ts_j\ts_k\ts_{i,4},\\
	\tomega^\vee_{j,i}:=&\ttau _{f_j,k}\ts_{j,k}\ts_{k,4}\ts_{i,i}\ts_{j,i}
	=\ttau _{f_j}\ts_j\ts_k\ts_i\ts_{j,i},\\
	\tomega^\vee_{i,i}:=&\ts_{0,i}\ts_{i,4}\ts_{j,j}\ts_{k,j}\ts_{j,4}\ts_{i,i}
	=\ts_0\ts_i\ts_j\ts_k\ts_j\ts_{i,i}=
	\ts_0\ts_i\ts_k\ts_j\ts_k\ts_{i,i}. 
\end{aligned}
	\label{eq:tomega}
\end{equation}
Note that $\mfp (\tomega ^\vee _{i,d})=\omega ^\vee _{i,d}$ for all
$i\in \isv \setminus \{0\}$ and $d\in \sd \setminus \{0\}$.

\begin{remar}
In order to understand the above definitions it is important to note that
for all $i\in I$ the vertex $i$ is playing a special role in the Dynkin diagram
labeled by $i$.
\end{remar}

	The elements of $\tWext $ defined above satisfy the equations
	$$ \te_d\tomega^\vee_{i,d}=\tomega^\vee_{i,d}\te_d=\tomega^\vee_{i,d}$$
	for all $i\in \isv \setminus \{0\}$ and $d\in \sd \setminus \{0\}$.
Further, for $i$, $j$, $k$ as above let
\begin{align}
	\tnu _{i,i}:=&\ttau _{f_i}\ts_i\ts_k\ts_{j,i},&
	\tnu _{j,i}:=&\ttau _{f_j}\ts_j\ts_k\ts_{i,i},&
	\tnu _{i,4}:=\ts_0\ts_i\ts_j\ts_k\ts_{j,4}.
	\label{eq:tnu}
\end{align}

\begin{lemma} One has
	\begin{align*}
		\tnu _{i,n_i\tr d}\ts_{i,d}=\tomega^\vee_{i,d}
	\end{align*}
	for all $i\in \isv \setminus \{0\}$ and $d\in \sd \setminus \{0\}$.
	\label{le:tnu}
\end{lemma}

\begin{proof} This follows immediately from Eq.s~\eqref{eq:tnu} and \eqref{eq:tomega}.
\end{proof}

\subsection{Some commutation relations in $\boldsymbol{\tWext }$}

In \cite[Lemma\,2.7]{a-Lusztig89} Lusztig studies the braid group
of the extended Weyl group of an affine Lie algebra.
In our setting the following related formulas
are valid.

\begin{theorem}\label{th:main}
\textrm{(1)} For all $i,j\in \isv \setminus \{0\}$ and
$d\in \sd \setminus \{0\}$ one has
\begin{align}
\tomega^\vee_{i,d}\tomega^\vee_{j,d}
=\tomega^\vee_{j,d}\tomega^\vee_{i,d}
\label{eq:crtomega1}
\end{align}

\textrm{(2)} Assume that $\{i,j,k\}=\{1,2,3\}$ and $d\in \sd \setminus \{0\}$.
Then one has
\begin{align}
\tnu _{i,d}\tomega^\vee_{i,d}
=&\ts_{i,d}(\tomega^\vee_{j,d})^{m_{i,j;d}-2}
(\tomega^\vee_{k,d})^{m_{i,k;d}-2},
\label{eq:crtomega2}\\
\tomega^\vee_{i,n_j\tr d}\ts_{j,d}
=&\ts_{j,d}\tomega^\vee_{i,d}.
\label{eq:crtomega3}
\end{align}
\end{theorem}

\begin{proof}
Suppose that $i\not=j$ and $d=4$.
One calculates
\begin{align*}
  \tomega^\vee_{i,d}\tomega^\vee_{j,d}=&
  \ttau _{f_i,4}\underline{\ts_{i,i}\ts_{k,i}\ts_{j,i}\ts_{i,4}
  \ttau _{f_j,4}}\ts_{j,j}\ts_{k,j}\ts_{i,j}\ts_{j,4}
%\\=&
=\underline{\ttau _{f_i,4}\ttau _{f_j,4}}\ts_{k,k}
    \underline{\ts_{i,k}\ts_{0,k}}\ts_{k,4}
    \ts_{j,j}\ts_{k,j}\ts_{i,j}\ts_{j,4} 
\\=&
\ttau _{f_k,4}\ts_{k,k}\ts_{0,k}\ts_{i,k}\underline{\ts_{k,4}
  \ts_{j,j}\ts_{k,j}}\ts_{i,j}\ts_{j,4}
%\\=&
=\ttau _{f_k,4}\ts_{k,k}\ts_{0,k}
    \underline{\ts_{i,k}\ts_{j,k}}
    \ts_{k,4}\underline{\ts_{j,j}\ts_{i,j}\ts_{j,4}}\\
  =&\ttau _{f_k,4}\ts_{k,k}\ts_{0,k}
    \ts_{j,k}\ts_{i,k}\ts_{k,4}\ts_{i,i}\ts_{j,i}\ts_{i,4},\\
  \tomega^\vee_{j,d}\tomega^\vee_{i,d}=&
  \ttau _{f_j,4}\underline{\ts_{j,j}\ts_{k,j}\ts_{i,j}\ts_{j,4}
  \ttau _{f_i,4}}\ts_{i,i}\ts_{k,i}\ts_{j,i}\ts_{i,4}
%\\=&
=\underline{\ttau _{f_j,4}\ttau _{f_i,4}}\ts_{k,k}
    \underline{\ts_{j,k}\ts_{0,k}}\ts_{k,4}
    \ts_{i,i}\ts_{k,i}\ts_{j,i}\ts_{i,4}\\
  =&\ttau _{f_k,4}\ts_{k,k}\ts_{0,k}\ts_{j,k}\underline{\ts_{k,4}
  \ts_{i,i}\ts_{k,i}}\ts_{j,i}\ts_{i,4}
%\\=&
=\ttau _{f_k,4}\ts_{k,k}\ts_{0,k}
    \ts_{j,k}\ts_{i,k}\ts_{k,4}\ts_{i,i}\ts_{j,i}\ts_{i,4}.
\end{align*}
The statement of part $(1)$ for $d\in \{1,2,3\}$ and part $(2)$
can be obtained analogously.
\end{proof}

\subsection{Symmetric bilinear forms}
\label{sec:sbf}

The affine Lie superalgebra $D^{(1)}(2,1;x)$ can be described with help of
different Dynkin diagrams. In this section we define symmetric bilinear forms
associated to all of these diagrams.
 
For any $d\in \sd $
let $V_d$ be a four dimensional
$\ndC $-vector space, and let
$\Pi_d=\{\al_{i,d}\,|\,i\in \isv \}$
be a basis of $V_d$.
Let $x\in\mathbb{C}\setminus\{0,-1\}$.
According to the Dynkin diagrams in Figure~\ref{fig:Dynkin},
for each $d\in \sd $ define a symmetric bilinear form
$\sbf{\,}{\,}=\sbf{\,}{\,}_d:V_d\times V_d \rightarrow \ndC $ as follows:
\begin{align*}
\sbf{\al_{i,4}}{\al_{i,4}}=&0\quad \mbox{for $i\in \isv $},&
\sbf{\al_{0,4}}{\al_{3,4}}=&\sbf{\al_{1,4}}{\al_{2,4}}=-1, \\
\sbf{\al_{0,4}}{\al_{1,4}}=&\sbf{\al_{2,4}}{\al_{3,4}}=-x,&
\sbf{\al_{0,4}}{\al_{2,4}}=&\sbf{\al_{1,4}}{\al_{3,4}}=x+1,
\end{align*}
\begin{align*}
\sbf{\al_{0,0}}{\al_{0,0}}=&0,& 
\sbf{\al_{i,0}}{\al_{j,0}}=&0\quad \text{for $i,j\in \{1,2,3\}$, $i\not=j$,}\\
\sbf{\al_{1,0}}{\al_{1,0}}=&-2x,& \sbf{\al_{1,0}}{\al_{0,0}}=&x, \\
\sbf{\al_{2,0}}{\al_{2,0}}=&2(x+1),& \sbf{\al_{2,0}}{\al_{0,0}}=&-x-1, \\
\sbf{\al_{3,0}}{\al_{3,0}}=&-2,& \sbf{\al_{3,0}}{\al_{0,0}}=&1,
\end{align*}
and
\begin{align*}
\sbf{\al_{i,d}}{\al_{j,d}}=&\sbf{\al_{f_d(i),0}}{\al_{f_d(j),0}}
\end{align*}
for $d\in\{1,2,3\}$. 
For $d\in \sd $ define a $\mathbb{Z}$-module map
$p=p_d:\mathbb{Z}\Pi_d\rightarrow\mathbb{Z}$
by 
\begin{align*}
p(\al_{i,d}):=
\begin{cases}
0 & \text{if $\sbf{\al _{i,d}}{\al _{i,d}}\not=0$,}\\
1 & \text{if $\sbf{\al _{i,d}}{\al _{i,d}}=0$}
\end{cases}
\end{align*}
  and call $p(\alpha )$ the \textit{parity} of $\alpha \in \mathbb{Z}\Pi_d$.
The next lemma follows immediately from Eq.~\eqref{eq:mijd} and the
definition of $\sbf{\,}{\,}$.

\begin{lemma}
	One has
	\begin{align*}
		m_{i,j;d}=&2 \quad \text{if $(\al _{i,d},\al _{j,d})=0$,}&
		m_{i,j;d}=&3 \quad \text{if $(\al _{i,d},\al _{j,d})\not=0$.}
	\end{align*}
	for all $d\in \sd $ and $i,j\in \isv $ with $i\not=j$.
	\label{le:mij}
\end{lemma}

In the following lemma we give a  representation of $\Wext $ which is
compatible with the symmetric bilinear form defined above.

\begin{lemma}\label{le:Wextrep}
Let ${\bf V}:=\oplus_{d=0}^4 V_d$.
Then there exists a unique semigroup
homomorphism
${\bf t}:\Wext \rightarrow \End _\ndC ({\bf V})$
such that
\begin{align*}
&{\bf t}(0)=0,\quad
{\bf t}(e_d)={\rm{id}}_{V_d},& 
{\bf t}(\tau _{f,d})(\al_{i,d})=&\al_{f(i),f\tr d},\\
&{\bf t}(s_{i,d})(\al_{i,d})=-\al_{i,n_i\tr d},&
{\bf t}(s_{i,d})(\al_{j,d})=&\al_{j,n_i\tr d}
+(m_{i,j;d}-2)\al_{i,n_i\tr d}
\end{align*}
for all $i,j\in \isv $ and $d\in \sd $ with $i\not=j$.
Moreover, for $w\in \Wext $ with $w e_d=w$
and for $v$, $v'\in V_d$ and $\mu\in\mathbb{Z}\Pi_d$
we have
\begin{align}
\sbf{{\bf t}(w)(v)}{{\bf t}(w)(v')}=\sbf{v}{v'}\quad \mbox{and}\quad 
(-1)^{p({\bf t}(w)(\mu))}=(-1)^{p(\mu)}.
\label{eq:compsbf}
\end{align}
Further, if $p(\al _{i,d})=0$ then $n_i\tr d=d$, $(\al _{i,d},\al _{i,d})\not=0$, and
\begin{align*}
	{\bf t}(s_{i,d})(\al_{j,d})=\al_{j,d}-
	\frac{2\sbf{\al _{i,d}}{\al _{j,d}}}{\sbf{\al _{i,d}}{\al _{i,d}}}\al _{i,d}
\end{align*}
for all $i,j\in \isv $ and $d\in \sd $.
\end{lemma}

\begin{proof}
  One has to check that the definition of $\bf t$ is compatible
with the relations \eqref{eq:rel1}--\eqref{eq:rel9}. This is trivial
for all relations different from \eqref{eq:rel5} and \eqref{eq:rel6},
but also easy for the latter. For example, if $i,j,k\in \isv $ are pairwise
distinct then one gets
\begin{equation}
\begin{aligned}
  s_{i,i}s_{j,i}s_{i,4}(\al _{k,4})=
  s_{i,i}s_{j,i}(\al _{k,i}+\al _{i,i})=&
  s_{i,i}(\al _{k,i}+\al _{i,i}+\al _{j,i})\\
  =&\al _{k,4}+\al _{i,4}+\al _{j,4}=
  s_{j,j}s_{i,j}s_{j,4}(\al _{k,4}).
\end{aligned}
\label{eq:simplecomp}
\end{equation}
Further, Eq.~\eqref{eq:compsbf}
has to be checked for generators $w$ of $\Wext $. Again all calculations
are easily done.
\end{proof}

For affine Lie (super)algebras there exists a distinguished root
$\iisor $ (see Section~\ref{sec:sectiontwo}), which
should be considered here.
For all $d\in \sd $ set
\begin{align}\label{eq:isoroot}
\isor _4:=&\sum _{i=0}^3\al _{i,4},&
\isor _d:=\al _{d,d}+\sum_{i=0}^3\al_{i,d} \quad \mbox{for $d\in \{0,1,2,3\}$}.
\end{align}
Then these are the elements corresponding to 
$\iisor $. 
Note that for all $d\in \sd $ one has
\begin{align*}
\ndC \isor_d =\{\lambda\in V_d\,|\,(\lambda,\mu)=0
\,\mbox{for all $\mu\in V_d$} \}.
\end{align*}
Further, using Eq.~\eqref{eq:omega} and Lemma~\ref{le:Wextrep}
one gets by computations similar to the one in Eq.~\eqref{eq:simplecomp}
the following formulas:
\begin{gather}
	\label{eq:repomega}
	{\bf t}(\omega^\vee_{i,d})(\al_{i,d})=\al_{i,d}-\isor_d,\quad
	{\bf t}(\omega^\vee_{i,d})(\al_{j,d})=\al_{j,d},\quad
	{\bf t}(w)(\isor _d)=\isor _{d'}
\end{gather}
for all $i,j\in \isv \setminus \{0\}$ with $i\not=j$,
and all $w\in \Wext $ and $d,d'\in \sd \setminus \{0\}$ such that
$e_{d'}we_d=w$.

\section{Quantum affine superalgebras of $\boldsymbol{D^{(1)}(2,1;x)}$}

Drinfeld \cite{Drinfeld} gave a second realization of quantum affine algebras
$U_q(\hat{\mfg })$.
He identified the generators of his algebra as loop-like generators in
$U_q(\hat{\mfg })$. We follow Beck's method \cite{beck1}
to define the analog of Drinfeld's second realization for $D^{(1)}(2,1;x)$.
First we introduce the quantum affine superalgebra of $D^{(1)}(2,1;x)$
for any Dynkin diagram. Then in Section~\ref{sec:brgraction}
we will give Lusztig type isomorphisms between these algebras, and
observe that these isomorphisms satisfy the relations of the braid
groupoid.

\subsection{Quantum affine superalgebras $\boldsymbol{U^\prime_d}$}

Fix $\hslash \in\mathbb{C}\setminus \ndZ \pi \sqrt{-1}$.
For any $u\in \ndC $ let
\begin{align*}
q^u:=&\exp(u\hslash )
=\sum_{n=0}^\infty{\frac {(u\hslash )^n} {n!}},&
q:=&q^1,&
[u]_q:=&\frac {q^u-q^{-u}} {q-q^{-1}}.
\end{align*}
In this paper we assume that
\begin{align}
q^{ku}\ne 1\quad \mbox{for all $u\in\{1,x,x+1\}$ 
and $k\in\mathbb{N}$}.
\label{eq:qCond}
\end{align}

Let $d\in \sd$. Let $\cU '_d$
be the $\ndC $-algebra (with $1$)
generated by 
the elements
\begin{align} \label{eq:Ugen}
\sigma _d,\quad K_{i,d}^{{\frac 1 2}},\quad K_{i,d}^{-{\frac 1 2}},\quad
E_{i,d},\quad F_{i,d},\quad \text{where $i\in \isv $}
\end{align}
and defined by the relations
\eqref{eq:cUr1}--\eqref{eq:cUr5} below:
\begin{gather}
  XY=YX\,\,\mbox{for }X,Y\in \{\sigma _d,
  K_{i,d}^{{\frac 1 2}},K_{i,d}^{-{\frac 1 2}}\,|\,i\in \isv \},
  \label{eq:cUr1}\\
  \sigma _d^2=1,\,\,K_{i,d}^{{\frac 1 2}}K_{i,d}^{-{\frac 1 2}}
  =K_{i,d}^{-{\frac 1 2}}K_{i,d}^{{\frac 1 2}}=1,
  \label{eq:cUr2}\\
  \sigma _dE_{i,d}\sigma _d=(-1)^{p(\al_{i,d})}E_{i,d},\qquad \quad
  \sigma _dF_{i,d}\sigma _d=(-1)^{p(\al_{i,d})}F_{i,d},
  \label{eq:cUr3}\\
  K_{i,d}^{{\frac 1 2}}E_{j,d}K_{i,d}^{-{\frac 1 2}}=
  q^{\sbf{\al _{i,d}}{\al _{j,d}}/2}E_{j,d},\quad
  K_{i,d}^{{\frac 1 2}}F_{j,d}K_{i,d}^{-{\frac 1 2}}=
  q^{-\sbf{\al _{i,d}}{\al _{j,d}}/2}F_{j,d},
  \label{eq:cUr4}\\
  E_{i,d}F_{j,d}-(-1)^{p(\al _{i,d})p(\al _{j,d})}F_{j,d}E_{i,d}
  %=\delta _{ij}\frac {K_{i,d}-K_{i,d}^{-1}}{q-q^{-1}} 
  =\delta _{ij}\frac {(K_{i,d}^{\frac 1 2})^2-(K_{i,d}^{-\frac 1 2})^2}{q-q^{-1}} 
  \label{eq:cUr5}
\end{gather}
for all $i,j\in \isv $.
% where $K_{i,d}:=(K_{i,d}^{\frac 1 2})^2$
%and $K_{i,d}^{-1}:=(K_{i,d}^{-{\frac 1 2}})^2$.
Note that for all $i\in I$ the equations $(K_{i,d}^{\frac 1 2})^{-1}=K_{i,d}^{-{\frac 1 2}}$ 
and  $(K_{i,d}^{\frac 1 2})^0=1$ hold.
Later on we will also use the abbreviations
\begin{align}
&K_{i,d}:=(K_{i,d}^{\frac 1 2})^2,\quad K_{i,d}^{\frac m 2}:=(K_{i,d}^{\frac 1 2})^m, 
\quad K_{i,d}^m:=(K_{i,d}^{\frac m 2})^2
\quad
\text{for all $m\in\ndZ$,}	\label{eq:preKmu} \\
&K_{\mu ;d}:=\prod _{i=0}^3K_{i,d}^{\frac {m_i} 2}, \quad
	\text{for all $\mu ={\frac 1 2}\sum _{i=0}^3m_i\al _{i,d}\in
	{\frac 1 2}\ndZ \Pi _d$ with $m_i\in\ndZ$.}\label{eq:Kmu}
\end{align} 

In particular, according to the definition of $\isor _d$
in Eq.~\eqref{eq:isoroot} we have
\begin{align}\label{eq:isoK}
	K_{\isor _d;4}=&\prod _{i=0}^3K_{i,4},&
	K_{\isor _d;d}=K_{d,d}\prod _{i=0}^3K_{i,d}
	\quad \mbox{for $d\in \{0,1,2,3\}$}.
\end{align}

The algebra $\cU '_d$ admits a unique $\ndZ \Pi _d$-grading (see
Section~\ref{sec:sbf} for the definition of $\Pi _d$)
\begin{align*}
  \cU '_d=&\bigoplus _{\lambda \in \ndZ \Pi _d}\cU '_{d,\lambda },&
  1\in &\cU '_{d,0},&
  \cU '_{d,\mu }\cU '_{d,\lambda }\subset &\cU '_{d,\mu +\lambda }
  \quad \text{for all $\mu ,\lambda \in \ndZ \Pi _d$,}
\end{align*}
such that $\sigma _d,K_{i,d}^{\pm \frac 1 2}\in \cU '_{d,0}$,
$E_{i,d}\in \cU '_{d,\al _{i,d}}$, and
$F_{i,d}\in \cU '_{d,-\al _{i,d}}$ for all $i\in \isv $.
Further, there exists a unique algebra automorphism $\Psi _d$ of $\cU '_d$
such that
\begin{align*}
  \Psi _d(\sigma _d)=&\sigma _d,&
  \Psi _d(K_{i,d}^{\pm \frac 1 2})=&K_{i,d}^{\mp \frac 1 2},&
  \Psi _d(E_{i,d})=&(-1)^{p(\al _{i,d})}F_{i,d},&
  \Psi _d(F_{i,d})=&E_{i,d}.
%  \label{eq:Phi}
\end{align*}
  Notice that $\Psi _d^2(X)=\sigma _d X\sigma _d$ for all $x\in \cU '_d$.

\begin{notation}\label{no:bracket}
	(Super-bracket and $q$-super-bracket)
For elements of $\cU '_d$ we use the {\it{super-bracket}} $[\,,\,]$
and the {\it{$q$-super-bracket}} $\wlbr\,,\,\wrbr$
defined as follows. For any $\mu ,\lambda \in \ndZ \Pi _d$, $a\in \ndC $,
and $X_\mu \in \cU '_{d,\mu }$, $Y_\lambda \in \cU '_{d,\lambda }$
we let
\begin{gather}
[X_\mu,Y_\lambda ]_a:=X_\mu Y_\lambda 
-(-1)^{p(\mu )p(\lambda )}aY_\lambda X_\mu ,
\label{eq:genbr}\\
[X_\mu,Y_\lambda ]:=[X_\mu,Y_\lambda ]_1,\qquad
\wlbr X_\mu,Y_\lambda \wrbr :=[X_\mu,Y_\lambda ]_{q^{-(\mu ,\lambda )}}
\label{eq:wbr}
\end{gather}
\end{notation}

Now we define the quantum affine superalgebras $U'_d$ of $D^{(1)}(2,1;x)$
for the Dynkin diagrams labeled by $d\in \sd $.

\begin{defin}\label{definition:defU'd}
For any $d\in \sd $ let $U'_d$ be the
quotient of the algebra $\cU '_d$ by the ideal generated by
the following elements
(see also \cite[Remark~7.1.1]{yamane1}
(or Prop.6.3.1(vii),(viii) in q-alg/9603015)):
\begin{align}
  E_{i,d}^2, \quad &\text{where $i\in \isv $ and $p(\al _{i,d})=1$,}
  \label{eq:Ur6}\\
  \wlbr E_{i,d},E_{j,d}\wrbr , \quad & \text{where $i,j\in \isv $, $i\not=j$,
  and $(\al_{i,d},\al_{j,d})=0$,} 
  \label{eq:Ur7}\\
  \wlbr E_{i,d} , \wlbr E_{i,d} , E_{j,d} \wrbr\wrbr , \quad &
  \text{where $i,j\in \isv $, $i\ne j$, $p(\al_{i,d})=0$,
  and $\sbf{\al _{i,d}}{\al _{j,d}}\not=0$,} 
  \label{eq:Ur8}\\
\label{eq:Ur9}
[(\al_{i,4},\al_{k,4})]_q
\wlbr \wlbr E_{i,4},&E_{j,4}\wrbr ,E_{k,4}\wrbr
-[(\al_{i,4},\al_{j,4})]_q
\wlbr \wlbr E_{i,4},E_{k,4}\wrbr ,E_{j,4}\wrbr ,\\
& \text{if $d=4$,  where $i,j,k\in \isv $ such that $i<j<k$,} \notag \\
[(\al_{i,d}+\al _{d,d},\al_{k,d}+&\al_{d,d})]_q
\wlbr \wlbr \wlbr E_{d,d},E_{i,d}\wrbr ,
\wlbr E_{d,d},E_{j,d}\wrbr \wrbr ,\wlbr E_{d,d},E_{k,d}\wrbr \wrbr -\notag \\
[(\al_{i,d}+\al _{d,d},&\al_{j,d}+\al_{d,d})]_q
\wlbr \wlbr \wlbr E_{d,d},E_{i,d}\wrbr ,
\wlbr E_{d,d},E_{k,d}\wrbr \wrbr ,\wlbr E_{d,d},E_{j,d}\wrbr \wrbr
\label{eq:Ur10}\\
&\mbox{
if $d\not=4$, where $\{i,j,k,d\}=\isv $ and $i<j<k$,} \notag \\
\Psi _d(X), \quad & \text{for all $X$ in \eqref{eq:Ur6}--\eqref{eq:Ur10}.}
\label{eq:Ur11}
\end{align}
\end{defin}

Notice that $\Psi _d$ induces an automorphism of $U'_d$, which will also
be denoted by $\Psi _d$. Further, $U'_d$ inherits a $\ndZ \Pi _d$-grading from
$\cU '_d$, that is
\begin{align*}
  U '_d=&\bigoplus _{\lambda \in \ndZ \Pi _d}U '_{d,\lambda },&
  1\in &U '_{d,0},&
  U '_{d,\mu }U '_{d,\lambda }\subset &U '_{d,\mu +\lambda }
  \quad \text{for all $\mu ,\lambda \in \ndZ \Pi _d$,}
\end{align*}
and hence Notation~\ref{no:bracket} can be applied also
to the elements of $U'_d$.

\begin{lemma}
	\label{le:relbrwbr}
	For any $\mu ,\lambda ,\xi \in \ndZ \Pi _d$ and $X_\mu\in U '_{d,\mu }$,
	$Y_\lambda \in U '_{d,\lambda }$,
	$Z_\xi \in U '_{d,\xi }$,
	(or, alternatively,
	$X_\mu\in \cU '_{d,\mu }$, $Y_\lambda \in \cU '_{d,\lambda }$,
	$Z_\xi \in \cU '_{d,\xi }$,)
	one has the following formulas.
	\begin{gather*}
		\wlbr X_\mu ,K_{\lambda ;d}Y_\lambda \wrbr =
		q^{-(\lambda ,\mu )}K_{\lambda ;d}[X_\mu ,Y_\lambda ],\qquad
		\wlbr K_{\mu ;d}^{-1} X_\mu ,Y_\lambda \wrbr =
		K_{\mu ;d}^{-1}[X_\mu ,Y_\lambda ],\\
		\wlbr \wlbr X_\mu ,Y_\lambda \wrbr ,Z_\xi \wrbr =
		\wlbr X_\mu ,\wlbr Y_\lambda ,Z_\xi \wrbr \wrbr 
		+(-1)^{p(\lambda )p(\xi )}q^{-(\lambda ,\xi)}[\wlbr X_\mu ,
		Z_\xi \wrbr ,Y_\lambda ]_{q^{(\lambda ,\xi -\mu )}},\\
		\wlbr X_\mu ,\wlbr Y_\lambda ,Z_\xi \wrbr \wrbr =
		\wlbr \wlbr X_\mu ,Y_\lambda \wrbr ,Z_\xi \wrbr 
		+(-1)^{p(\mu )p(\lambda )}q^{-(\mu ,\lambda )}[Y_\lambda ,
		\wlbr X_\mu ,Z_\xi \wrbr ]_{q^{(\lambda ,\mu -\xi )}}.
	\end{gather*} 
\end{lemma}

\begin{proof}
	The first two equations follow from
	Eq.s~\eqref{eq:wbr} and \eqref{eq:cUr4}.
	For the other relations one needs only
	the definition of the $q$-super-bracket.
\end{proof}

\begin{lemma}\label{le:Ufu}
Let $d\in \sd $ and $i$, $j\in \isv $ with $i\ne j$.
Then the following equations hold in $\cU '_d$ and in $U'_d$.
\begin{align}
	[\wlbr E_{j,d},E_{i,d}\wrbr , F_{i,d}]
	=&-[(\al _{i,d},\al _{j,d})]_qK_{i,d}^{-1} E_{j,d},
	\label{eq:Ufu1}\\
	[\wlbr E_{j,d},E_{i,d}\wrbr , F_{j,d}]
	=&(-1)^{p(\al _{i,d})p(\al _{j,d})}[(\al _{i,d},\al _{j,d})]_q
	E_{i,d} K_{j,d},
	\label{eq:Ufu2}\\
	[E_{i,d},\wlbr F_{j,d},F_{i,d} \wrbr ]
	=&(-1)^{p(\al _{i,d})p(\al _{j,d})}[(\al _{i,d},\al _{j,d})]_q
	K_{i,d}F_{j,d}, 
	\label{eq:Ufu3}\\
	[E_{j,d}, \wlbr F_{j,d},F_{i,d}\wrbr ]
	=&-[(\al _{i,d},\al _{j,d})]_q F_{i,d}K_{j,d}^{-1}
	\label{eq:Ufu4}\\
	[ \wlbr E_{j,d},E_{i,d}\wrbr , \wlbr F_{j,d},F_{i,d}\wrbr ]
	=&
	\label{eq:Ufu5}\\
	(-1)^{p(\al _{i,d})p(\al _{j,d})}&q^{-(\al _{i,d},\al _{j,d})}
	[(\al _{i,d},\al _{j,d})]_q
	\frac {K _{i,d}K _{j,d}
	 -K _{i,d}^{-1}K _{j,d}^{-1}}{q-q^{-1}} \notag
\end{align}
\end{lemma}

\begin{proof}
	Eq.s~\eqref{eq:Ufu1},\eqref{eq:Ufu2}
	can be checked directly by using
	Eq.s~\eqref{eq:cUr4},\eqref{eq:cUr5}.
	Then
	Eq.s~\eqref{eq:Ufu3},\eqref{eq:Ufu4} can be obtained
	from \eqref{eq:Ufu1},\eqref{eq:Ufu2} by applying
	$\Psi _d$. Finally, Eq.~\eqref{eq:Ufu5} follows from
	formulas~\eqref{eq:Ufu1}--\eqref{eq:Ufu4} by using
	\eqref{eq:cUr1},\eqref{eq:cUr4}.
\end{proof}

Let $\cUpp _d$, $\cUpm _d$, and $\cUpo $
be the subalgebras of $\cU'_d$ generated by the sets
$\{E_{i,d}\,|\,i\in \isv \}$,
$\{F_{i,d}\,|\,i\in \isv \}$, and
$\{\sigma _d ,K_{i,d}^{\pm \frac 1 2}\,|\,i\in \isv \}$,
respectively. Let $\Upp _d$, $\Upm _d$, and $\Upo $ denote the images
of $\cUpp _d$, $\cUpm _d$, and $\cUpo $, respectively, under
the canonical projection $\cU '_d\to U'_d$.
Notice that $\Upm _d=\Psi _d(\Upp _d)$ and $\Psi _d(\Upo _d)=\Upo _d$.

\begin{theorem}\label{th:HopfTri}
{\rm{(1)}}
The $\mathbb{C}$-algebras $\cU'_d$ and $U'_d$
can be regarded as Hopf algebras
$(\cU '_d,\Delta,S,\varepsilon)$ and $(U '_d,\Delta,S,\varepsilon)$
such that
\begin{gather*}
\Delta(X)=X\otimes X,\,\,S(X)=X^{-1},\,\,\varepsilon(X)=1\,\,
\mbox{for $X\in \{\sigma _d, K_{i,d}^{\pm{\frac 1 2}}\,|\,i\in \isv \}$,}
\\
\begin{aligned}
	\Delta(E_{i,d})=&E_{i,d}\otimes 1 +K_{i,d}\sigma _d^{p(\al_{i,d})}
	\otimes E_{i,d},&
	S(E_{i,d})=&-K_{i,d}^{-1}\sigma _d^{p(\al_{i,d})}E_{i,d},\\
	\Delta(F_{i,d})=&F_{i,d}\otimes K_{i,d}^{-1} +\sigma _d^{p(\al _{i,d})}
	\otimes F_{i,d},&
	S(F_{i,d})=&-(-1)^{p(\al_{i,d})}F_{i,d}K_{i,d}\sigma _d^{p(\al_{i,d})},\\
\end{aligned}\\
	\varepsilon(E_{i,d})=\varepsilon(F_{i,d})=0
\end{gather*}
for all $i\in \isv $.

{\rm{(2)}} The algebras $\cU '_d$ and $U' _d$ admit a triangular decomposition.
More precisely, the multiplication maps
\begin{align*}
\cUpp _d\otimes \cUpo \otimes \cUpm _d &\to \cU '_d,&
\cUpm _d\otimes \cUpo \otimes \cUpp _d &\to \cU '_d,\\
\Upp _d\otimes \Upo \otimes \Upm _d &\to U '_d,&
\Upm _d\otimes \Upo \otimes \Upp _d &\to U '_d,
\end{align*}
where $X_1\otimes X_2 \otimes X_3$ is mapped to $X_1X_2X_3$
for all $X_1,X_2,X_3$,
are isomorphisms of $\ndZ \Pi _d$-graded $\ndC $-vector spaces.
Further, the algebra $\cUpp _d$ is the free algebra generated by
the set $\{E_{i,d}\,|\,i\in \isv \}$, and $\Upp _d$
is isomorphic to its quotient by the ideal generated by the elements in
\eqref{eq:Ur6}--\eqref{eq:Ur10}.
The algebras $\cUpo $ and $\Upo $
are both isomorphic to the commutative algebra
generated by the set $\{\sigma _d,K_{i,d}^{\pm \frac 1 2}\,|\,i\in \isv \}$
and defined by the relations~\eqref{eq:cUr2},
so the set 
\begin{align*}
  \Big\{\sigma _d^m\prod_{i=0}^3(K_{i,d}^{{\frac 1 2}})^{n_i}\,|\,
  m\in\{0,1\},n_i\in \mathbb{Z}\text{ for all $i\in \isv $}\Big\}
\end{align*}
is a $\mathbb{C}$-basis of both $\cUpo _d$ and $\Upo _d$.
\end{theorem}

\begin{proof}
	Part (2) of the theorem is standard
	\cite[Prop.\,4.16,\,Thm.\,4.21]{b-Jantzen96}.
	The compatibility of $\Delta $ and $\varepsilon $ with
	the defining relations of $\cU '_d$, and the axioms of
	$S$ can be easily checked. In order to prove that $\Delta $ is
	well-defined on $\Upp _d$ and $\Upm _d$
	we used the computer algebra program
	Mathematica.
\end{proof}

\subsection{Lusztig type isomorphisms}
\label{sec:brgraction}

The main result of this section is Theorem~\ref{th:tzerotone}, which tells
that the elements of the braid semigroup $\tWext $ defined in
Section~\ref{sec:brsgroup} can be represented by isomorphisms between the
quantum affine algebras $U' _d$.

\begin{notation}\label{no:rij}
	For any $d\in \sd $ and $i,j\in \isv $ with $i\not=j$ one has
	\begin{align*}
	  (-1)^{p(\al _{i,d})p(\al _{j,d})}[(\al _{i,d},\al _{j,d})]_q\in 
	  \{0,1,[x]_q,[-x-1]_q\}.
	\end{align*}
	We fix a square root of all of these four complex numbers, and
	write
	\begin{align*}
	  r_{i,j;d}:=1\Big/\sqrt{(-1)^{p(\al _{i,d})p(\al _{j,d})}
	  [(\al _{i,d},\al _{j,d})]_q}
	\end{align*}
	if $(\al _{i,d},\al _{j,d})\not=0$. In this case one also has
	$r_{i,j;n_i\tr d}=r_{i,j;d}=r_{j,i;d}$.
	Further, let
	$\rsign _{i,j;d}=(-1)^{\delta _{ij}+p(\al _{i,d})p(\al _{j,d})}$.
	(As for the symbol $\delta _{ij}$, 
    see the last paragraph in Introduction.)
	%where $\delta $ denotes Kronecker's $\delta $-function.
\end{notation}

\begin{theorem}\label{th:exisT}
	The following statements are valid.

{\rm{(1)}} For all $d\in \sd $ and $i\in \isv $ there exist unique
$\mathbb{C}$-algebra isomorphisms
\begin{align}
	T_{i,d},T^-_{i,d}:U'_d \rightarrow U'_{n_i\tr d} 
\label{eq:Ti0}
\end{align}
satisfying Eq.s~\eqref{eq:Ti1}--\eqref{eq:Ti7} below.
\begin{gather}
	T_{i,d}(\sigma _d)=T^-_{i,d}(\sigma _d)=\sigma _{n_i\tr d},
	\label{eq:Ti1}\\
	\begin{aligned}
	  T_{i,d}(K_{i,d}^{\frac 1 2})=&T^-_{i,d}(K_{i,d}^{\frac 1 2})=
	  K_{i,n_i\tr d}^{-{\frac 1 2}},\\
	  T_{i,d}(K_{j,d}^{\frac 1 2})=&T^-_{i,d}(K_{j,d}^{\frac 1 2})
	  =K_{j,n_i\tr d}^{\frac 1 2}
	  (K_{i,n_i\tr d}^{{\frac 1 2}})^{m_{i,j;d}-2},
	\end{aligned}
        \label{eq:Ti2}\\
	\intertext{where $j\in \isv$ and $j\not=i$,}
	\begin{aligned}
	  T_{i,d}(E_{i,d})=&(-1)^{p(\al _{i,n_i\tr d})}
	  q^{-{\frac {(\al _{i,n_i\tr d},\al _{i,n_i\tr d})} 2}}
	  F_{i,n_i\tr d}K_{i,n_i\tr d},\\
	  T_{i,d}(F_{i,d})=&
	  q^{{\frac {(\al _{i,n_i\tr d},\al _{i,n_i\tr d})} 2}}
	  K^{-1}_{i,n_i\tr d}E_{i,n_i\tr d},
	  \label{eq:Ti3}
	\end{aligned}\\
	\begin{aligned}
	  T^-_{i,d}(E_{i,d})=&
	  q^{-{\frac {(\al _{i,n_i\tr d},\al _{i,n_i\tr d})} 2}}
	  K^{-1}_{i,n_i\tr d}F_{i,n_i\tr d},\\
	  T^-_{i,d}(F_{i,d})=&(-1)^{p(\al _{i,n_i\tr d})}
	  q^{{\frac {(\al _{i,n_i\tr d},\al _{i,n_i\tr d})} 2}}
	  E_{i,n_i\tr d}K_{i,n_i\tr d},
	  \label{eq:Ti4}
	\end{aligned}\\
	\begin{aligned}
	  T_{i,d}(E_{j,d})=&r_{i,j;d}
	  q^{{\frac {(\al_{i,n_i\tr d},\al _{j,n_i\tr d})} 2}}
	  \wlbr E_{j,n_i\tr d} , E_{i,n_i\tr d} \wrbr ,\\
	  T_{i,d}(F_{j,d})=&r_{i,j;d}
	  q^{{\frac {(\al_{i,n_i\tr d},\al _{j,n_i\tr d})} 2}}
	  \wlbr F_{j,n_i\tr d} , F_{i,n_i\tr d} \wrbr ,
	  \label{eq:Ti5}
	\end{aligned}\\
	\begin{aligned}
	  T^-_{i,d}(E_{j,d})=&r_{i,j;d}
	  q^{{\frac {(\al_{i,n_i\tr d},\al _{j,n_i\tr d})} 2}}
	  \wlbr E_{i,n_i\tr d} , E_{j,n_i\tr d} \wrbr ,\\
	  T^-_{i,d}(F_{j,d})=&r_{i,j;d}
	  q^{{\frac {(\al_{i,n_i\tr d},\al _{j,n_i\tr d})} 2}}
	  \wlbr F_{i,n_i\tr d} , F_{j,n_i\tr d} \wrbr ,
	  \label{eq:Ti6}
	\end{aligned}\\
	\intertext{where in Eq.s~\eqref{eq:Ti5},\eqref{eq:Ti6}
	one has $j\in \isv $ with $j\not=i$ and
	$(\al _{i,d},\al _{j,d})\ne 0$,}
	T_{i,d}(E_{j,d})= T^-_{i,d}(E_{j,d})=
	E_{j,n_i\tr d},\quad
	T_{i,d}(F_{j,d})= T^-_{i,d}(F_{j,d})=
	F_{j,n_i\tr d}
	\label{eq:Ti7}
\end{gather}
where $j\in \isv $ with $j\not=i$ and $(\al _{i,d},\al _{j,d})=0$.

{\rm{(2)}} One has $T^-_{i,n_i\tr d}=(T_{i,d})^{-1}$ for all $d\in \sd $
and $i\in \isv $.

{\rm{(3)}} The isomorphisms $T_{i,d}$ satisfy the equations
$\Psi _{n_i\tr d}T_{i,d}=T_{i,d}\Psi _d$.
\end{theorem}

\begin{proof}
Parts (2) and (3) of the theorem are obtained easily
from the definition of $T_{i,d}$, $T^-_{i,d}$, and $\Psi _d$,
and from Lemma~\ref{le:Ufu}.
The uniqueness of $T_{i,d}$ and $T^-_{i,d}$ are obvious,
but to check the compatibility of $T_{i,d}$, $T^-_{i,d}$ with the defining
relations of $U'_d$ we need again the computer algebra program
Mathematica.
\end{proof}

\begin{theorem}\label{th:Tcoxrel}
For any $d\in \sd $ the following statements hold.
%Let $i_1$, $i_2\in\{0,1,2,3\}$
%be such that $i_1\ne i_2$.
%Let $\al:=\al_{i_1,a(j)}$, 
%$\beta:=\al_{i_2,a(j)}$,
%$\al^\prime:=\al_{i_1,n(i_2)\tr a(j)}$,
%$\beta^\prime:=\al_{i_2,n(i_1)\tr a(j)}$,
%$\al^{\prime\prime}:=\al_{i_1,n(i_2)n(i_1)\tr a(j)}$ and
%$\beta^{\prime\prime}:=\al_{i_2,n(i_1)n(i_2)\tr a(j)}$. \par

{\rm{(1)}} If $i,j\in \isv $, where $i\not=j$ and $m_{i,j;d}=2$, then
\begin{align}
	T_{i,d}T_{j,d}=&T_{j,d}T_{i,d},&
	T_{i,d}(E_{j,d})=&E_{j,d},&
	T_{i,d}(F_{j,d})=&F_{j,d}.
	\label{eq:Tbr1}
\end{align}

{\rm{(2)}} If $i,j\in \isv $, where $i\not=j$ and $m_{i,j;d}=3$, then
the following equations hold.
\begin{gather}
	T_{i,n_jn_i\tr d} T_{j,n_i\tr d} T_{i,d} =
	T_{j,n_in_j\tr d} T_{i,n_j\tr d} T_{j,d}
	\label{eq:Tbr2},\\
	T_{i,n_j\tr d} T_{j,d}(E_{i,d})=E_{j,n_in_j\tr d},\qquad
	T_{i,n_j\tr d} T_{j,d}(F_{i,d})=F_{j,n_in_j\tr d}.
	\label{eq:Tbr3}
\end{gather}
\end{theorem}

\begin{proof}
	All statements of the theorem follow from Theorem~\ref{th:exisT},
	the definition of the algebras $U'_d$, and from
	Lemma~\ref{le:Ufu}.
\end{proof}

For all $d\in \sd $ and $f\in \Klf $, see Eq.~\eqref{eq:A4},
let $T_{f,d}:U'_d\to U'_{f\tr d}$ denote the $\mathbb{C}$-algebra
isomorphism satisfying the following equations.
\begin{align*}
	T_{f,d}(\sigma _d)=&\sigma _{f\tr d},&
	T_{f,d}(K_{i,d}^{\pm{\frac 1 2}})=&K_{f(i),f\tr d}^{\pm{\frac 1 2}},&\\
	T_{f,d}(E_{i,d})=&E_{f(i),f\tr d},&
	T_{f,d}(F_{i,d})=&F_{f(i),f\tr d},
\end{align*}
where $i\in \isv $.

\begin{defin}
	Let $\cM _\sd =\cup _{d,d'\in \sd }\mathrm{Hom}\,(U'_d,U'_{d'})
	\cup \{0\}$
	be a disjoint union of sets, where $\mathrm{Hom}\,(U'_d,U'_{d'})$
	denotes the set of unital algebra maps from $U'_d$ to $U'_{d'}$.
	The set $\cM _\sd $ admits a unique semigroup structure with the
	following properties.
%	\begin{align*}
%		0\phi =&\phi 0=0,\\
%		\phi _1\, \phi _2=&
%		\begin{cases}
%			\phi _1\circ \phi _2 & \text{if $d_2=d_3$,}\\
%			0 & \text{otherwise,}
%		\end{cases}
%	\end{align*}
\begin{align*}
0\phi =\phi 0=0,\quad
\phi _1\, \phi _2=
		\begin{cases}
			\phi _1\circ \phi _2 & \text{if $d_2=d_3$,}\\
			0 & \text{otherwise,}
		\end{cases}
\end{align*}
	for all $\phi \in \cM _\sd $ and all
	$\phi _1\in \mathrm{Hom}\,(U'_{d_1},U'_{d_2})$ and
	$\phi _2\in \mathrm{Hom}\,(U'_{d_3},U'_{d_4})$, where
	$d_1,d_2,d_3,d_4\in \sd $.
\end{defin}

Note that the set $\cM _\sd \setminus \{0\}$ in the above definition
can also be considered as the morphisms of a category with objects
$U'_d$, where $d\in \sd $.

Recall the definition of $\tWext $ from Section~\ref{sec:brsgroup}.

\begin{theorem}\label{th:tzerotone}
	There exists a unique semigroup homomorphism
	\begin{align}
		\cT :\tWext \rightarrow \cM _\sd
		\label{eq:T0}
	\end{align}
	such that for all $d\in \sd $, $i\in \isv $, and $f\in \Klf $ one has
	\begin{equation}
			\cT (0)=0,\quad 
			\cT (\te_d)=\id :U'_d\to U'_d,\quad
			\cT (\ts _{i,d})=T_{i,d},\quad
			\cT (\ttau _{f,d})=T_{f,d}.
		\label{eq:T1}
	\end{equation}
	Moreover, for all ${\tilde w}\in \tWext $
	with ${\tilde w}\ne 0$ and $\te_{d'}{\tilde w}\te_d={\tilde w}$ for some
	$d,d'\in \sd $, we have
	\begin{align}
		\cT ({\tilde w})(U'_{d,\mu })
		=&U'_{d',{\bf t}(\mfp ({\tilde w}))(\mu)}
		& &\text{for all $\mu \in \ndZ \Pi _d$},
		\label{eq:T2}\\
		\cT ({\tilde w})(\sigma _d)=&\sigma _{d'},\,\,
		\cT ({\tilde w})(K_{\mu ;d})=
		K_{{\bf t}(\mfp (\tilde w))(\mu );d'}
		& & \text{for all $\mu \in \frac{1}{2}\ndZ \Pi _d$,}
		\label{eq:T3}\\
		\cT (\tilde w)\Psi _d=&\Psi _{d'}\cT (\tilde w).
		\label{eq:T4}
	\end{align}
\end{theorem}

\begin{proof}
	One has to show that Eq.~\eqref{eq:T1} is compatible with
	the defining relations of $\tWext $. This follows from
	Theorem~\ref{th:Tcoxrel} and the definition of the maps $T_{f,d}$,
	where $f\in \Klf $ and $d\in \sd $.
\end{proof}

Note that by Eq.~\eqref{eq:compsbf} also the formula
\begin{align}
	\cT (\tilde w)(\wlbr X_\mu ,Y_\lambda \wrbr )=
	\wlbr \cT (\tilde w)(X_\mu ),
	\cT (\tilde w)(Y_\lambda )\wrbr 
	\label{eq:compTwbr}
\end{align}
holds for all ${\tilde w}\in \tWext $, $X_\mu \in U'_{d,\mu }$,
and $Y_\lambda \in U'_{d,\lambda }$, where
$d\in \sd $, $\mu ,\lambda \in \ndZ \Pi _d$, ${\tilde w}\ne 0$,
and ${\tilde w}e_d={\tilde w}$.

\section{Root vectors associated to imaginary roots}
\label{sec:imagroots}

The aim of this section is to construct root vectors to the
imaginary roots $k\isor _d$ for all $k\in \ndN $ and
$d\in \sd \setminus \{0\}$.
First we will give some technical results on
the isomorphisms $\cT (\tomega ^\vee _{i,d})$, where
$\tomega ^\vee _{i,d}$ and $\cT $ are as in Eq.~\eqref{eq:tomega}
and Theorem~\ref{th:tzerotone} respectively. Then we define
root vectors $\bpsi _{i,k;d}$ of weight $k\isor _d$ for each
$i\in \isv \setminus \{0\}$, $d\in \sd \setminus \{0\}$,
and $k\in \ndN $. Finally we prove that these root vectors
commute with each other if they belong to the same algebra $U'_d$.

\subsection{Preliminary definition of root vectors to imaginary roots}

We start with the calculation of various values of
the isomorphisms $\cT (\tomega ^\vee _{i,d})$.

\begin{lemma} \label{le:propertyTomega} Let $d\in \sd \setminus \{0\}$ and
	$i$, $j\in \isv \setminus \{0\}$ with $i\ne j$.

{\rm{(1)}} The following formulas hold.
\begin{align*}
	\cT (\tomega^\vee_{i,d})(E_{j,d})=&E_{j,d},&
	\cT (\tomega^\vee_{i,d})(F_{j,d})=&F_{j,d},&
	\cT (\tomega^\vee_{i,d})(K^{\pm \frac 1 2}_{j,d})=
	&K_{j,d}^{\pm \frac 1 2}.
\end{align*}

{\rm{(2)}} We have
\begin{equation*}
\begin{aligned}
	\cT (\tomega^\vee_{i,d})(\wlbr E_{i,d},E_{j,d}\wrbr )
	=& \cT(\tomega^\vee_{j,d})(\wlbr E_{j,d},E_{i,d}\wrbr ),\\
	\cT (\tomega^\vee_{i,d})(\wlbr F_{i,d},F_{j,d}\wrbr )
	=& \cT(\tomega^\vee_{j,d})(\wlbr F_{j,d},F_{i,d}\wrbr ).
\end{aligned}
\end{equation*}
\end{lemma}

\begin{proof}
Part~(1) can be checked easily by using Eq.~\eqref{eq:Tbr3}
and Theorem~\ref{th:tzerotone}.

(2) If $(\al_{i,d},\al_{j,d})= 0$ then Lemma~\ref{le:propertyTomega}(2)
holds by Formula~\eqref{eq:Ur7}. Assume now that
$(\al_{i,d},\al_{j,d})\ne 0$.
%$$\kappa:=
%q^{{\frac {(\al_{i,a(j)},\al_{x,a(j)})} 2}}
%\Bigl(\sqrt{(-1)^{p(\al_{i,a(j)})p(\al_{x,a(j)})}
%[(\al_{i,a(j)},\al_{x,a(j)})]}\Bigr)^{-1}.
%$$
Using Theorem~\ref{th:tzerotone} we obtain the following formulas.
\begin{align*}
&q^{{\frac {(\al_{i,d},\al_{j,d})} 2}}r_{i,j;d}
\cT (\tomega^\vee_{i,d})(\wlbr E_{i,d},E_{j,d}\wrbr ) & & \\
&\quad =\cT (\tomega^\vee_{i,d})T_{j,n_j\tr d}(E_{i,n_j\tr d})
& &\mbox{(by Eq.~\eqref{eq:Ti5})}\\
&\quad =T_{j,n_j\tr d}\cT (\tomega^\vee_{i,n_j\tr d})(E_{i,n_j\tr d}) 
& &\mbox{(by Eq.~\eqref{eq:crtomega3})}\\
&\quad = \cT (\tnu _{j,n_j\tr d})\cT (\tomega^\vee _{j,n_j\tr d})
 (E_{i,n_j\tr d})
 & &\mbox{(by Eq.~\eqref{eq:crtomega2},
 Lemma~\ref{le:propertyTomega}(1))}\\
% BCIH <
% &\quad =\cT(\tnu _{j,n_j\tr d})(E_{i,d})
&\quad =\cT(\tnu _{j,n_j\tr d})(E_{i,n_j\tr d})
% > ECIH
& &\mbox{(by Lemma~\ref{le:propertyTomega}(1))}\\
&\quad =\cT(\tnu _{j,n_j\tr d})T_{j,d}T^-_{j,n_j\tr d}(E_{i,n_j\tr d})
& &\mbox{(by Theorem~\ref{th:exisT}(2))}\\
&\quad = 
q^{{\frac {(\al_{i,d},\al_{j,d})} 2}}r_{i,j;d}
\cT(\tomega^\vee_{j,d})(\wlbr E_{j,d},E_{i,d}\wrbr )
& &\mbox{(by Lemma~\ref{le:tnu} and Eq.~\eqref{eq:Ti5}).}
\end{align*}
Hence we have the first equation. The second one can be obtained
from the first one by applying $\Psi _d$ and using Eq.~\eqref{eq:T4}.
\end{proof}

Recall from Notation~\ref{no:rij} that for $d\in\sd$
and $\{i,j,k,l\}= I$ with $p(\al_{i,d})=1$ 
one has the equation
$(r_{i,j;d}r_{i,k;d}r_{i,l;d})^{-2}=[x]_q[-x-1]_q$.

\begin{lemma} \label{le:TomegaFKone}
	Let $d\in \sd \setminus \{0\}$ and $i$, $j$, $k\in \{1,2,3\}$
	such that $\{i,j,k\}=\{1,2,3\}$. 

{\rm{(1)}} The following formulas hold.
\begin{align*}
	\cT (\tomega^\vee_{i,d})(K_{i,d}^{-1}F_{i,d})
	=&q^{{\frac {(\al _{i,d},\al _{i,d})} 2}}
	\cT (\tomega^\vee_{i,d})T^-_{i,n_i\tr d}(E_{i,n_i\tr d})
	= q^{{\frac {(\al_{i,d},\al _{i,d})} 2}}
	\cT(\tnu _{i,n_i\tr d})(E_{i,n_i\tr d})\\
=&
\begin{cases}
-[(\al_{0,4},\al_{k,4})]_q^{-1}\wlbr E_{j,4},\wlbr E_{k,4},
E_{0,4}\wrbr\wrbr & \text{if $d=4$}, \\
r_{i,0;i}r_{i,j;i}r_{i,k;i}\wlbr E_{j,i},\wlbr E_{k,i},
\wlbr E_{i,i} , E_{0,i}\wrbr\wrbr\wrbr
 & \mbox{if $d=i$},
 \\
 r_{j,0;j}r_{j,i;j}r_{j,k;j}\wlbr E_{j,j},\wlbr E_{k,j},
\wlbr E_{j,j} , E_{0,j}\wrbr\wrbr\wrbr
 & \mbox{if $d=j$}.
\end{cases}
\end{align*}

{\rm{(2)}} For all $l,m\in \ndZ $ we have the following equations.
\begin{gather*}
	\wlbr \cT(\tomega^\vee_{i,d})(K_{i,d}^{-1}F_{i,d}),
	K_{i,d}^{-1}F_{i,d}\wrbr =0, \qquad
	\wlbr E_{i,d}, \cT(\tomega^\vee_{i,d})^{-1}(E_{i,d}) \wrbr =0,\\
	\wlbr \cT (\tomega ^\vee _{i,d})^l(E_{i,d}),
	\cT (\tomega ^\vee _{j,d})^m(K_{j,d}^{-1}F_{j,d})\wrbr=0.
\end{gather*}
\end{lemma}

\begin{proof}
	Part (1) can be proven directly by using
~\eqref{eq:Tbr1},\eqref{eq:Tbr3}.
The first equation of part~(2)
can be proven by using part~(1) and \ref{eq:Ufu2}.
The second equation of part~(2) follows from the first one
by applying the algebra map $\cT(\tomega^\vee_{i,d})^{-1}\circ \Psi _d$.

Let now $l,m\in \ndZ $. One gets
\begin{align*}
	& \wlbr \cT (\tomega ^\vee _{i,d})^l(E_{i,d}),
	\cT (\tomega ^\vee _{j,d})^m(K_{j,d}^{-1}F_{j,d})\wrbr\\
	&=
	\cT (\tomega ^\vee _{j,d})^m (
	\wlbr \cT (\tomega ^\vee _{j,d})^{-m}\cT (\tomega ^\vee _{i,d})^l
	(E_{i,d}),
	K_{j,d}^{-1}F_{j,d}\wrbr )&
	&\text{(by Eq.~\eqref{eq:compTwbr}})\\
	&= \cT (\tomega ^\vee _{j,d})^m (
	\wlbr \cT (\tomega ^\vee _{i,d})^l\cT (\tomega ^\vee _{j,d})^{-m}
	(E_{i,d}),K_{j,d}^{-1}F_{j,d}\wrbr )
	& &\text{(by Theorems~\ref{th:main}(1),
	\ref{th:tzerotone}})\\
	&= \cT (\tomega ^\vee _{j,d})^m (
	\wlbr \cT (\tomega ^\vee _{i,d})^l(E_{i,d}),
	K_{j,d}^{-1}F_{j,d}\wrbr )&
	&\text{(by Lemma~\ref{le:propertyTomega}(1))}\\
	&= \cT (\tomega ^\vee _{j,d})^m
	\cT (\tomega ^\vee _{i,d})^l(\wlbr E_{i,d},
	K_{j,d}^{-1}F_{j,d}\wrbr )&
	&\text{(as in the previous steps)}\\
	&= \cT (\tomega ^\vee _{j,d})^m
	\cT (\tomega ^\vee _{i,d})^l(q^{(\al _{i,d},\al _{j,d})}
	K_{j,d}^{-1}[E_{i,d},F_{j,d}]) & &
	\text{(by Lemma~\ref{le:relbrwbr})}\\
	&=0 &
	&\text{(by Eq.~\eqref{eq:cUr5}).}
\end{align*}
\end{proof}

%We have:
%\begin{eqnarray*}
%\lefteqn{\wlbr E_\al , T(\tomega^\vee_\al)^{-1}(
%E_\al) \wrbr } \\
%&=& q^{-(\al,\al)}[E_\al K_\al ,
%T^\prime_{\al^\prime}T(\tomega^\vee_{\al^\prime})
%T^\prime_\al (E_\al ) ] \\
%&=& q^{-(\al,\al)}(-1)^{p(\al^\prime)}q^{-(\al^\prime,\al^\prime)}
%T^\prime_{\al^\prime}([F_{\al^\prime},
%T(\tomega^\vee_{\al^\prime})(K_{\al^\prime}^{-1}F_{\al^\prime})]) \\
%&=& 0.
%\end{eqnarray*}\hfill $\Box$

%\subsection{Formulas 
%of $[{\bar \psi}^{(s)}_{\al k} , T(\tomega^\vee_\beta)^m(E_\al) ]$
%and $[{\bar \psi}^{(s)}_{\al k} , T(\tomega^\vee_\beta)^m(F_\al) ]$
%and notation
%$\stackrel{k}{\equiv}$}

%\subsection{Formulas 
%of $[{\bar \psi}^{(0)}_{\al k} , E_\beta ]$
%and $[{\bar \psi}^{(1-k)}_{\al k} , F_\beta]$}

In order to define root vectors for integer multiples of $\isor _d$,
we now define
for all $i\in \isv \setminus \{0\}$, $k\in \ndN $, and
$d\in \sd \setminus \{0\}$
a family $\{\bpsi ^{(s)}_{i,k;d}\,|\,s\in \ndZ \}$ of elements
of $U'_d$. In Proposition~\ref{pr:allpsis} it will be shown
that the cardinality of each of these families is one,
and their unique elements
will be considered as the root vectors associated to the pairs
$(\al _{i,d},k\isor _d)$.

\begin{defin}\label{de:t0imrv}
Recall the definition of $K_{\isor _d;d}$ from Eq.~\eqref{eq:isoK}.
For $d\in \sd \setminus \{0\}$, $i\in \isv \setminus \{0\}$,
$k\in\mathbb{N}$, and $s\in\mathbb{Z}$
let
\begin{align} \label{eq:bpsis}
\bpsi ^{(s)}_{i,k;d}=&
(-1)^kq^{-(\al _{i,d},\al_{i,d})}
\cT(\tomega^\vee_{i,d})^s(\wlbr E_{i,d},\cT(\tomega^\vee_{i,d})^k(
K_{i,d}^{-1}F_{i,d}) \wrbr ) \\
=&(-1)^k q^{-(\al _{i,d},\al _{i,d})}
\wlbr \cT(\tomega^\vee_{i,d})^s(E_{i,d}),
\cT(\tomega^\vee_{i,d})^{s+k}(
K_{i,d}^{-1}F_{i,d}) \wrbr \notag \\
%=&(-1)^kq^{-(\al _{i,k},\al _{i,k})}\cT(\tomega^\vee_{i,d})^s(E_{i,d})
%\cT(\tomega^\vee_{i,d})^{s+k}(K_{i,d}^{-1}F_{i,d}) \notag \\
%& \quad -(-1)^{p(\al _{i,d})+k}\cT(\tomega^\vee_{i,d})^{s+k}(
%K_{i,d}^{-1}F_{i,d})\cT(\tomega^\vee_{i,d})^s(E_{i,d}) \notag \\
=& (-1)^kK_{\isor _d;d}^{s+k}K_{i,d}^{-1}
[T(\tomega^\vee_{i,d})^s(E_{i,d}),T(\tomega^\vee_{i,d})^{s+k}(F_{i,d})],
\notag
\end{align}
where the last equation follows from Theorem~\ref{th:tzerotone},
Eq.~\eqref{eq:repomega}, and Lemma~\ref{le:relbrwbr}.
Note that $\bpsi ^{(s)}_{i,k;d}\in U'_{d,k\isor _d}$.
\end{defin}

We will also need the element
\begin{align}
	\bpsi ^{(0)}_{i,0;d}:=q^{-(\al _{i,d},\al_{i,d})} \wlbr E_{i,d},
	K_{i,d}^{-1}F_{i,d}\wrbr ={\frac {1-K_{i,d}^{-2}} {q-q^{-1}}}
	\in U'_{d,0}.
	\label{eq:bpsi0}
\end{align}

\begin{lemma}
	\label{le:Tomega(bpsi)}
	Let $d\in \sd \setminus \{0\}$, $i,j\in \isv \setminus \{0\}$,
	$k\in \ndN $, and $s\in \ndZ $. Then we have
	\begin{align*}
		\cT (\tomega ^\vee _{i,d})(\bpsi ^{(s)}_{j,k;d})=
		\begin{cases}
			\bpsi ^{(s+1)}_{j,k;d} & \text{if $i=j$,}\\
			\bpsi ^{(s)}_{j,k;d} & \text{if $i\not=j$.}
		\end{cases}
	\end{align*}
\end{lemma}

\begin{proof}
	The statement in the case $i=j$ follows from the definition of
	$\bpsi ^{(s)}_{j,k;d}$. If $i\not=j$ then one obtains the claim
	by Eq.~\eqref{eq:crtomega1} and Lemma~\ref{le:propertyTomega}(1).
\end{proof}

\begin{lemma} \label{le:psiTomegaEF} Let $d\in \sd \setminus \{0\}$
	and $i$, $j\in \isv \setminus \{0\}$. Then for all $k\in\mathbb{N}$
	the following equations hold.
	\begin{align*}
		[\bpsi ^{(-1)}_{i,k;d}, E_{j,d}]=&
		\rsign _{i,j;d}q^{(\al _{i,d},\al _{j,d})}
		[\bpsi ^{(0)}_{i,k-1;d},\cT(\tomega^\vee_{j,d})^{-1}(E_{j,d})]
		_{q^{-2(\al _{i,d},\al _{j,d})}},\\
		[\bpsi ^{(1-k)}_{i,k;d},K_{j,d}^{-1}F_{j,d}]=&
		\rsign _{i,j;d}q^{-(\al _{i,d},\al _{j,d})}
		[\bpsi ^{(1-k)}_{i,k-1;d},
		\cT(\tomega^\vee_{j,d})(K_{j,d}^{-1}F_{j,d})
		]_{q^{2(\al _{i,d},\al _{j,d})}}.
	\end{align*}
%\begin{eqnarray*}
%\lefteqn{[{\bar \psi}^{(1-k)}_{\al k} , K_\beta^{-1}F_\beta ]} \\
%&=& (-1)^{\delta_{\al\beta}+p(\al)p(\beta)}
%(q^{-(\al,\beta)}{\bar \psi}^{(1-k)}_{\al k-1}T(\tomega^\vee_\beta) ( K_\beta^{-1}F_\beta) 
%-q^{(\al,\beta)}T(\tomega^\vee_\beta) 
%( K_\beta^{-1}F_\beta){\bar \psi}^{(1-k)}_{\al k-1} ).
%\end{eqnarray*}
In particular, one has the formulas
\begin{align*}
	[\bpsi ^{(-1)}_{i,1;d} , E_{j,d}]=&\rsign _{i,j;d}
	[(\al_{i,d},\al _{j,d})]_q\cT (\tomega^\vee_{j,d})^{-1}(E_{j,d}),\\
	[\bpsi ^{(0)}_{i,1;d} , K_{j,d}^{-1}F_{j,d}]=&-\rsign _{i,j;d}
	[(\al_{i,d},\al _{j,d})]_q
	\cT (\tomega^\vee_{j,d})(K_{j,d}^{-1}F_{j,d}).
\end{align*}
\end{lemma}

\begin{proof}
Assume first that $i=j$.
By Eq.~\eqref{eq:bpsis},
Lemma~\ref{le:TomegaFKone}(2), and Lemma~\ref{le:relbrwbr}
we have the following.
\begin{align*}
	&[\bpsi ^{(-1)}_{i,k;d}, E_{i,d}]=-\wlbr E_{i,d},\bpsi ^{(-1)}_{i,k;d}
	\wrbr \\
	&\quad =(-1)^{k+1}q^{-(\al _{i,d},\al _{i,d})}
	\wlbr E_{i,d},\wlbr \cT(\tomega^\vee_{i,d})^{-1}(E_{i,d}),
	\cT(\tomega^\vee_{i,d})^{k-1}(K_{i,d}^{-1}F_{i,d}) \wrbr \wrbr\\
	&\quad =(-1)^{k+1}q^{-(\al _{i,d},\al _{i,d})}
	(-1)^{p(\al _{i,d})+1}q^{(\al _{i,d},\al _{i,d})}\\
	&\qquad \quad [ \wlbr
	E_{i,d},
	\cT(\tomega^\vee_{i,d})^{k-1}(K_{i,d}^{-1}F_{i,d})\wrbr ,
	\cT(\tomega^\vee_{i,d})^{-1}(E_{i,d})
	]_{q^{-2(\al _{i,d},\al _{i,d})}} \\
	&\quad =\rsign _{i,i;d}q^{(\al _{i,d},\al _{i,d})}
	[\bpsi ^{(0)}_{i,k-1;d},\cT(\tomega^\vee_{i,d})^{-1}(E_{i,d})]
	_{q^{-2(\al _{i,d},\al _{i,d})}},\\
	&[\bpsi ^{(1-k)}_{i,k;d},K_{i,d}^{-1}F_{i,d}]=
	\wlbr \bpsi ^{(1-k)}_{i,k;d},K_{i,d}^{-1}F_{i,d}\wrbr \\
	&\quad =(-1)^kq^{-(\al _{i,d},\al _{i,d})}
	\wlbr \wlbr 
	\cT(\tomega^\vee_{i,d})^{1-k}(E_{i,d}),
	\cT(\tomega^\vee_{i,d})(K_{i,d}^{-1}F_{i,d}) \wrbr ,
	K_{i,d}^{-1}F_{i,d} \wrbr\\
	&\quad =(-1)^kq^{-(\al _{i,d},\al _{i,d})}
	(-1)^{p(\al _{i,d})}q^{-(\al _{i,d},\al _{i,d})}\\
	&\qquad \quad [ \wlbr
	\cT(\tomega^\vee_{i,d})^{1-k}(E_{i,d}) ,
	K_{i,d}^{-1}F_{i,d} \wrbr ,
	\cT(\tomega^\vee_{i,d})(K_{i,d}^{-1}F_{i,d})
	]_{q^{2(\al _{i,d},\al _{i,d})}} \\
	&\quad =\rsign _{i,i;d}q^{-(\al _{i,d},\al _{i,d})}
	[\bpsi ^{(1-k)}_{i,k-1;d},\cT(\tomega^\vee_{i,d})(K_{i,d}^{-1}F_{i,d})]
	_{q^{2(\al _{i,d},\al _{i,d})}}.
\end{align*}
Assume now that $i\not=j$.
First notice that
because of  
%By 
Lemmata~\ref{le:TomegaFKone} and \ref{le:propertyTomega}
we have the following.
\begin{align*}
	&[\bpsi ^{(-1)}_{i,k;d},E_{j,d}]=
	-\wlbr E_{j,d},\bpsi ^{(-1)}_{i,k;d}\wrbr \\
	&\quad =(-1)^{k+1}q^{-(\al_{i,d},\al_{i,d})}\wlbr E_{j,d},
	\wlbr \cT(\tomega^\vee_{i,d})^{-1}(E_{i,d}),
	\cT(\tomega^\vee_{i,d})^{k-1}(K_{i,d}^{-1}F_{i,d})\wrbr \wrbr \\ 
	&\quad =(-1)^{k+1}q^{-(\al_{i,d},\al_{i,d})}\wlbr \wlbr E_{j,d},
	\cT(\tomega^\vee_{i,d})^{-1}(E_{i,d})\wrbr ,
	\cT(\tomega^\vee_{i,d})^{k-1}(K_{i,d}^{-1}F_{i,d})\wrbr \\ 
	&\quad =(-1)^{k+1}q^{-(\al_{i,d},\al_{i,d})}\wlbr
	\cT(\tomega^\vee_{i,d})^{-1}( \wlbr E_{j,d},
	E_{i,d}\wrbr ),
	\cT(\tomega^\vee_{i,d})^{k-1}(K_{i,d}^{-1}F_{i,d})\wrbr \\ 
	&\quad =(-1)^{k+1}q^{-(\al_{i,d},\al_{i,d})}\wlbr
	\cT(\tomega^\vee_{j,d})^{-1}( \wlbr E_{i,d},
	E_{j,d}\wrbr ),
	\cT(\tomega^\vee_{i,d})^{k-1}(K_{i,d}^{-1}F_{i,d})\wrbr \\ 
	&\quad =(-1)^{k+1}q^{-(\al_{i,d},\al_{i,d})}\wlbr
	\wlbr E_{i,d},\cT(\tomega^\vee_{j,d})^{-1}( E_{j,d})\wrbr ,
	\cT(\tomega^\vee_{i,d})^{k-1}(K_{i,d}^{-1}F_{i,d})\wrbr \\ 
	&\quad =(-1)^{k+1}q^{-(\al_{i,d},\al_{i,d})}
	(-1)^{p(\al _{i,d})p(\al _{j,d})}q^{(\al _{j,d},\al _{i,d})}\\
	&\qquad \quad [\wlbr E_{i,d},
	\cT(\tomega^\vee_{i,d})^{k-1}(K_{i,d}^{-1}F_{i,d})\wrbr , 
	\cT(\tomega^\vee_{j,d})^{-1}( E_{j,d})]_{q^{-2(\al _{j,d},\al _{i,d})}}
	\\
	&\quad =\rsign _{i,j;d}q^{(\al _{i,d},\al _{j,d})}
	[\bpsi _{i,k-1;d}^{(0)}, 
	\cT(\tomega^\vee_{j,d})^{-1}( E_{j,d})]_{q^{-2(\al _{i,d},\al _{j,d})}}.
\end{align*}
The remaining equation can be proved similarly.
\end{proof}

\subsection{Definition of type one imaginary root vectors}

The main result of this subsection is the following statement.

\begin{proposition} \label{pr:allpsis}
	Let $d\in \sd \setminus \{0\}$.

{\rm{(1)}}
One has $\bpsi ^{(s)}_{i,k;d}=\bpsi ^{(0)}_{i,k;d}$
for all $i\in \isv \setminus \{0\}$,
$k\in\mathbb{N}$, and $s\in\mathbb{Z}$.

{\rm{(2)}}
One has
$[\bpsi ^{(0)}_{i,k;d},\bpsi ^{(0)}_{j,r;d}]=0$
for all $i$, $j\in \isv \setminus \{0\}$ and
$k$, $r\in\mathbb{N}$. 
\end{proposition}

\begin{defin}\label{de:t1imrv}
	Let $\bpsi _{i,k;d}:=\bpsi ^{(0)}_{i,k;d}$ for all
	$i\in \isv \setminus \{0\}$, $k\in \ndN $, and
	$d\in \sd \setminus \{0\}$.
	For $k\in -\ndN $ and $i\in \isv \setminus \{0\}$,
	$d\in \sd \setminus \{0\}$
	set $\bpsi _{i,k;d}:=\Psi _d(\bpsi _{i,-k;d})$.
	The elements  $\bpsi _{i,k;d}\in U'_{d,k\isor _d}$, where
	$k\in \ndZ \setminus \{0\}$, are called
	\textit{type one imaginary root vectors.}
\end{defin}

In order to prove the above proposition we need a technical lemma.
We use the notation $\{k;c\}:=\sum _{j=0}^{k-1}c^{2j-k+1}$
for all $k\in \ndN $ and $c\in \ndC $.

\begin{lemma} \label{le:specalg}
	Let $d\in \sd $, and let $\epsilon ,\eta \in
	\{1,-1\}\subset \ndC $, $a$, $b$, $c\in \ndC \setminus \{0\}$,
	and $m,n\in\mathbb{N}$ with $m\le n$. Let $(X_u)_{0\le u\le n}$,
	$(Y_u)_{0\le u\le n}$, $(Z^+_v)_{v\in \ndZ }$, and
	$(Z^-_v)_{v\in \ndZ }$ be families
	of $\ndZ \Pi _d$-homogeneous elements of $U'_d$, with the following
	properties.
	\begin{enumerate}
		\item The parity of the $\ndZ \Pi _d$-degree of $X_u$
			is even for all $u\in \{0,1,\ldots ,n\}$.
		\item The families $(X_u)_{0\le u\le n}$,
			$(Y_u)_{0\le u\le n}$, $(Z^+_v)_{v\in \ndZ }$, and
			$(Z^-_v)_{v\in \ndZ }$
			satisfy the following equations.
\begin{equation}
	\begin{split}
		[X_u,Z^\pm _v]=&\epsilon (c^{\pm 1}X_{u-1}Z^\pm _{v\mp 1}
		-c^{\mp 1}Z^\pm _{v\mp 1}X_{u-1})\qquad
		\text{if $1\le u\le m$, $v\in \ndZ $,}\\
		Y_0=&aZ^+_0Z^-_0-bZ^-_0Z^+_0,\\
		Y_u=&\eta ^u(aZ^+_0Z^-_u-bZ^-_iZ^+_0)
        \quad\qquad \mbox{if $1\le u\le n-1$,}\\
		Y_u=&\eta^u(aZ^+_{-1}Z^-_{u-1}-bZ^-_{u-1}Z^+_{-1})
		\qquad \mbox{if $1\le u\le n-1$,}\\
		Y_n=&\eta ^n(aZ^+_{-1}Z^-_{n-1}-bZ^-_{n-1}Z^+_{-1}).
	\end{split}
	\label{eq:specalgrel0}
\end{equation}
	\item One has $[X_0,Y_n]=0$.
	\end{enumerate}
Then for all $u$ with $2\le u\le m$ we have
\begin{align}
	[X_u,Y_{n-u}]=(\eta \epsilon)^{u-1}\{u;c\}[X_1,Y_{n-1}].
	\label{eq:specalgrel1}
\end{align}

%In particular,
%$[{\cal H}^{(1)}_0,{\cal H}^{(2)}_n]=0$. \par
Moreover if there exists $r\in \ndC $ such that the equations
\begin{align}
[X_1,Z^+_0]=&rZ^+_{-1},&
[X_1,Z^-_{n-1}]=&-rZ^-_n
	\label{eq:specalgrel2}
\end{align}
hold, then one also has
\begin{align}
rY_n-\eta [X_1,Y_{n-1}]=
r\eta ^n (aZ^+_0Z^-_n-bZ^-_nZ^+_0).
\label{eq:specalgrel3}
\end{align}
\end{lemma}

\begin{proof}
	The strategy of the proof of the first part of the lemma
	is the following. First we prove that
	\begin{align}\label{eq:partialeq}
		[X_u,Y_{n-u}]=\eta \epsilon (c+c^{-1})[X_{u-1},Y_{n-u+1}]
		-[X_{u-2},Y_{n-u+2}]
	\end{align}
	for all $u\in \{2,\ldots ,m\}$. Then assumption 3 implies that
	$$[X_2,Y_{n-2}]=\eta \epsilon (c+c^{-1})[X_1,Y_{n-1}],$$
	and induction on $u$ using \eqref{eq:partialeq} implies
	Eq.~\eqref{eq:specalgrel1}.

	Now we prove Eq.~\eqref{eq:partialeq}.
	Using ~\eqref{eq:specalgrel0}
	and Jacobi identity for $[\,,\,]$ one obtains for all
$u\in\{2,\ldots ,m\}$ the following.

\begin{align*}
	[X_u,Y_{n-u}]
%\\&\quad 
%&
&=[X_u,\eta ^{n-u}aZ^+_0Z^-_{n-u}-\eta ^{n-u}bZ^-_{n-u}Z^+_0]\\
%	&\quad 
&=\eta ^{n-u}a\Big\{\epsilon (cX_{u-1}Z^+_{-1}
	-c^{-1}Z^+_{-1}X_{u-1})Z^-_{n-u}\\
	&\quad \qquad
	+Z^+_0\epsilon (c^{-1}X_{u-1}Z^-_{n-u+1}-cZ^-_{n-u+1}X_{u-1})\Big\}\\
	&\quad \qquad -\eta ^{n-u}b\Big\{
	\epsilon (c^{-1}X_{u-1}Z^-_{n-u+1}-cZ^-_{n-u+1}X_{u-1})Z^+_0\\
	&\quad \qquad
	+Z^-_{n-u}\epsilon (cX_{u-1}Z^+_{-1}-c^{-1}Z^+_{-1}X_{u-1})\Big\}.
\end{align*}
We calculate the first four summands of the last expression separately
by using Eq.s~\eqref{eq:specalgrel0}.
\begin{align*}
	&\eta ^{n-u}\epsilon acX_{u-1}Z^+_{-1}Z^-_{n-u}
	=\eta ^{n-u}\epsilon ac
	X_{u-1}\Big(\frac{Y_{n-u+1}}{\eta ^{n-u+1}a}+\frac{b}{a}
	Z^-_{n-u}Z^+_{-1}\Big)\\
	&\quad =\eta \epsilon cX_{u-1}Y_{n-u+1}+\eta ^{n-u}\epsilon bc
	(Z^-_{n-u}X_{u-1}
	+\epsilon c^{-1}X_{u-2}Z^-_{n-u+1}
	-\epsilon cZ^-_{n-u+1}X_{u-2})Z^+_{-1},\\
	&-\eta ^{n-u}\epsilon ac^{-1}Z^+_{-1}X_{u-1}Z^-_{n-u}\\
	&\quad =-\eta ^{n-u}\epsilon ac^{-1}Z^+_{-1}(Z^-_{n-u}X_{u-1}
	+\epsilon c^{-1}X_{u-2}Z^-_{n-u+1}
	-\epsilon cZ^-_{n-u+1}X_{u-2})\\
	&\quad =-\eta ^{n-u}\epsilon ac^{-1}
	\Big( \Big(\frac{Y_{n-u+1}}{\eta ^{n-u+1}a}+\frac{b}{a}
	Z^-_{n-u}Z^+_{-1}\Big)X_{u-1}\\
	&\quad \qquad +\epsilon c^{-1}Z^+_{-1}X_{u-2}Z^-_{n-u+1}
	-\epsilon cZ^+_{-1}Z^-_{n-u+1}X_{u-2}\Big),\\
	&\eta ^{n-u}\epsilon ac^{-1}Z^+_0X_{u-1}Z^-_{n-u+1}\\
	&\quad =\eta ^{n-u}\epsilon ac^{-1}(X_{u-1}Z^+_0-
	\epsilon cX_{u-2}Z^+_{-1}
	+\epsilon c^{-1}Z^+_{-1}X_{u-2})Z^-_{n-u+1},\\
	&\quad =\eta ^{n-u}\epsilon ac^{-1}X_{u-1}Z^+_0Z^-_{n-u+1}
	+\eta ^{n-u}ac^{-2}Z^+_{-1}X_{u-2}Z^-_{n-u+1}\\
	&\qquad \quad -X_{u-2}(
	Y_{n-u+2}+\eta ^{n-u}bZ^-_{n-u+1}Z^+_{-1}),\\
	&-\eta ^{n-u}\epsilon acZ^+_0Z^-_{n-u+1}X_{u-1}
	=-\eta \epsilon c(Y_{n-u+1}+\eta ^{n-u+1}bZ^-_{n-u+1}Z^+_0)
	X_{u-1}\\
	&\quad =-\eta \epsilon cY_{n-u+1}X_{u-1}+\eta ^{n-u}\epsilon bc
	Z^-_{n-u+1}(-X_{u-1}Z^+_0+\epsilon cX_{u-2}Z^+_{-1}-\epsilon c^{-1}
	Z^+_{-1}X_{u-2}).
\end{align*}
Comparison of the latter formulas gives that
\begin{align*}
	&[X_u,Y_{n-u}]\\
	&\quad =\eta \epsilon cX_{u-1}Y_{n-u+1}
	-\eta \epsilon c^{-1}Y_{n-u+1}X_{u-1}\\
	&\qquad \quad
	+\eta ^{n-u}aZ^+_{-1}Z^-_{n-u+1}X_{u-2}
	+\eta ^{n-u}\epsilon ac^{-1}X_{u-1}Z^+_0Z^-_{n-u+1}
	-X_{u-2}Y_{n-u+2}\\
	&\qquad \quad -\eta \epsilon cY_{n-u+1}X_{u-1}
	-\eta ^{n-u}bZ^-_{n-u+1}Z^+_{-1}X_{u-2}
	-\eta ^{n-u}\epsilon bc^{-1}X_{u-1}Z^-_{n-u+1}Z^+_0\\
	&\quad =\eta \epsilon c[X_{u-1},Y_{n-u+1}]
	-\eta \epsilon c^{-1}Y_{n-u+1}X_{u-1}\\
	&\qquad \quad +Y_{n-u+2}X_{u-2}
	+\eta \epsilon c^{-1}X_{u-1}Y_{n-u+1}-X_{u-2}Y_{n-u+2}\\
	&\quad =\eta \epsilon (c+c^{-1})[X_{u-1},Y_{n-u+1}]
	-[X_{u-2},Y_{n-u+2}].
\end{align*}
This proves Eq.~\eqref{eq:partialeq}.

Now we prove Eq.~\eqref{eq:specalgrel3}. We have
\begin{align*}
	rY_n=&r\eta ^naZ^+_{-1}Z^-_{n-1}-r\eta ^nbZ^-_{n-1}Z^+_{-1}\\
	=&\eta ^na(X_1Z^+_0-Z^+_0X_1)Z^-_{n-1}
	-\eta ^nbZ^-_{n-1}[X_1,Z^+_0]\\
	=& \eta X_1(Y_{n-1}+\eta ^{n-1}bZ^-_{n-1}Z^+_0)
	-\eta ^naZ^+_0(Z^-_{n-1}X_1-rZ^-_n)
	-\eta ^nbZ^-_{n-1}[X_1,Z^+_0]\\
	=& \eta X_1Y_{n-1}+\eta ^nb[X_1,Z^-_{n-1}]Z^+_0
	+\eta ^narZ^+_0Z^-_n-\eta ^n(aZ^+_0Z^-_{n-1}-bZ^-_{n-1}Z^+_0)X_1\\
	=&\eta X_1Y_{n-1}-\eta ^nbrZ^-_nZ^+_0
	+\eta ^narZ^+_0Z^-_n-\eta Y_{n-1}X_1,
\end{align*}
which gives \eqref{eq:specalgrel3}.
\end{proof}

\begin{proof}[Proof of Proposition~\ref{pr:allpsis}]
We prove both parts of the proposition simultaneously by induction on $k+r$,
where $r=0$ in part~(1).
Let $n\in\mathbb{N}$ and assume that we have proved
(1) for all $k\in \ndN $ with $k<n$ and (2) for all $k,r\in \ndN $ with
$k+r<n$.

First we assume that $i\not=j$. We want to apply Lemma~\ref{le:specalg}
with $m=n$ and the following setting.
\begin{align*}
	\epsilon :=&\epsilon _{i,j;d},&
	\eta :=&-1,& 
	r:=&\epsilon _{i,j;d}[(\al _{i,d},\al _{j,d})]_q,\\
	a:=&q^{-(\alpha _{i,d},\al _{i,d})},&
	b:=&(-1)^{p(\al _{j,d})},&
	c:=&q^{(\al _{i,d},\al _{j,d})},\\
X_u:=&\bpsi ^{(0)}_{i,u;d} \quad (0\le u\le n),&
Y_u:=&\bpsi ^{(0)}_{j,u;d} \quad (0\le u\le n-1),&
Y_n:=&\bpsi ^{(-1)}_{j,n;d},\\
Z^+_v&:=\cT (\tomega ^\vee _{j,d})^v(E_{j,d}),&
Z^-_v&:=\cT (\tomega ^\vee _{j,d})^v(K_{j,d}^{-1}F_{j,d}),
\end{align*}
for all $v\in \ndZ $. Note that the first formula in
Eq.~\eqref{eq:specalgrel0} holds by 
Lemmata~\ref{le:propertyTomega}(1), \ref{le:Tomega(bpsi)}
and \ref{le:psiTomegaEF}
and by the induction hypothesis for part (1). The fourth line
in \eqref{eq:specalgrel0} can be shown by applying
$\cT (\tomega _{j,d}^\vee)$ to the definition of $\bpsi ^{(0)}_{j,u;d}$
and using the induction hypothesis for part (1). All other formulas in
\eqref{eq:specalgrel0} hold essentially by the definition of
$Y_u$. Eq.s~\eqref{eq:specalgrel2} follow essentially from
Lemmata~\ref{le:propertyTomega}(1), \ref{le:Tomega(bpsi)} and \ref{le:psiTomegaEF}.
Eq.~\eqref{eq:bpsi0} implies assumption~3.
Thus Lemma~\ref{le:specalg} gives that
\begin{align}
[\bpsi^{(0)}_{i,u;d},\bpsi ^{(0)}_{j,n-u;d}]=(-\epsilon _{i,j;d})^{u-1}
\{u;q^{(\al _{i,d},\al _{j,d})}\}[\bpsi ^{(0)}_{i,1;d},\bpsi ^{(0)}_{j,n-1;d}],
	\label{eq:inotj-2}
\end{align}
for $2\leq u\leq n$. Let $u=n$ in the last formula. Then from the definition
of $\bpsi ^{(0)}_{j,0;d}$ one gets
\begin{align*}
	0=[\bpsi ^{(0)}_{i,n;d},\bpsi ^{(0)}_{j,0;d}]=
	(-\epsilon _{i,j;d})^{n-1}\{n;q^{(\al _{i,d},\al _{j,d})}\}
	[\bpsi ^{(0)}_{i,1;d},\bpsi ^{(0)}_{j,n-1;d}].
\end{align*}
Therefore the assumption in \eqref{eq:qCond} implies that
\begin{align}
[\bpsi ^{(0)}_{i,1;d},\bpsi ^{(0)}_{j,n-1;d}]=0.
\label{eq:comm1-n-1}
\end{align}
This with Eq.~\eqref{eq:inotj-2}
gives part (2) of the proposition for $i\not=j$ and $k+r=n$.
In order to prove part~(1) of the proposition, assume that
$i\not=j$ and $(\al _{i,d},\al _{j,d})\not=0$.
The second part of Lemma~\ref{le:specalg} gives that
\begin{align}
\label{eq:psijn}
\epsilon _{i,j;d}[ ( \al _{i,d},\al _{j,d}) ]_q
\bpsi ^{(-1)}_{j,n;d}=-[\bpsi ^{(0)}_{i,1;d},\bpsi ^{(0)}_{j,n-1;d}]
+\epsilon _{i,j;d}[ ( \al _{i,d},\al _{j,d}) ]_q
\cT (\tomega ^\vee _{j,d})(\bpsi ^{(-1)}_{j,n;d}).
\end{align}
and hence  part (1) of the proposition follows from
\eqref{eq:comm1-n-1} and Lemma~\ref{le:Tomega(bpsi)} together with
assumption~\eqref{eq:qCond} and relation
$[(\al _{i,d},\al _{j,d})]_q\not=0$.

Second we assume that $i=j$. 
We can apply the same argument as above,
together with the equation
%Since part (1) with $k\le n$ has already been
%proved, one has
\begin{align*}
  \bpsi ^{(-1)}_{i,k;d}=\cT (\tomega ^\vee _{i,d})(\bpsi ^{(0)}_{i,k;d})
  =\bpsi ^{(0)}_{i,k;d}
\end{align*}
for $1\le k\le n$,
since part (1) with $k\le n$ has already been
proved (cf. Lemma~\ref{le:Tomega(bpsi)}).
%Hence we can apply
%the argument obtained from the above one by letting $j=i$.
Then one has the same equations as in \eqref{eq:inotj-2}--\eqref{eq:comm1-n-1}
with $j=i$. As a result one gets part (2) for $i=j$.
\end{proof}

\subsection{Definition of type two imaginary root vectors}

We start with the definition.

\begin{defin}\label{de:t3imrv}
  Let $d\in \sd $ and $z$ a formal parameter.
  The coefficients $\bh_{i,k;d}\in U'_{d,k\isor _d}$ of the formal power series
  \begin{align}
    \sum _{k=1}^\infty \bh_{i,k;d}z^k:=\frac{1}{q-q^{-1}}\log \Big(
    1+\sum _{r=1}^\infty (q-q^{-1})\bpsi _{i,r;d}z^r\Big)
    \label{eq:t2irv}
  \end{align}
  together with the elements $\bh_{i,-k;d}:=-\Psi _d(\bh_{i,k;d})$,
  where $k\in \ndN $, are called \textit{type two imaginary root vectors}.
  
\end{defin}

Note that the above definition is equivalent to the more common formulas
  \begin{align*}
    \exp \Big( (q-q^{-1})\sum _{k=1}^\infty \bh_{i,k;d}z^k\Big)=&
    1+\sum _{r=1}^\infty (q-q^{-1})\bpsi _{i,r;d}z^r,\\
    \exp \Big(-(q-q^{-1})\sum _{k=1}^\infty \bh_{i,-k;d}z^k\Big)=&
    1+\sum _{r=1}^\infty (q-q^{-1})\bpsi _{i,-r;d}z^r.
  \end{align*}

In order to determine commutation relations between root vectors to real
roots and type two imaginary root vectors, we need two technical lemmata.
The first of them is standard in a more general setting.

\begin{lemma}\label{le:XZseries}
Let $d\in \sd $ and $\epsilon \in\ndC $, $c\in \ndC \setminus \{0\}$.
Let $(X_u)_{u\in \ndN }$ and $(Z_v)_{v\in \ndZ }$ be families of
$\ndZ \Pi _d$-homogeneous elements of $U'_d$ such that
the parities of the degrees of the elements $X_u$ are even for all
$u\in \ndN $.  Then the three conditions
{\rm{(i)}}, {\rm{(ii)}} and {\rm{(iii)}} below are equivalent.

{\rm{(i)}} For all $u\ge 2$ one has
$$ 
[X_u,Z_v]=\epsilon (cX_{u-1}Z_{v-1}-c^{-1}Z_{v-1}X_{u-1}).
$$

{\rm{(ii)}} For all $u\ge 1$ one has
$$ 
[X_u,Z_v]=(\epsilon c)^{u-1}[X_1,Z_{v-u+1}]
+\epsilon (c-c^{-1})\sum_{r=1}^{u-1}(\epsilon c)^{u-1-r}Z_{v-u+r}X_r.
$$

{\rm{(iii)}} For all $u\ge 1$ one has
$$
[X_u,Z_v]=(\epsilon c^{-1})^{u-1}[X_1,Z_{v-u+1}]
+\epsilon (c-c^{-1})\sum _{r=1}^{u-1}(\epsilon c^{-1})^{u-1-r}X_rZ_{v-u+r}.
$$
\end{lemma}

\begin{lemma} \label{le:L}
	Let $d\in \sd $, and let $\epsilon \in
	\{1,-1\}\subset \ndC $, and $c\in \ndC \setminus \{0\}$.
	Let $(X_u)_{u\in \ndN }$, $(Z^ +_v)_{v\in \ndZ }$ and
	$(Z^-_v)_{v\in \ndZ }$ be families
	of $\ndZ \Pi _d$-homogeneous elements of $U'_d$, with the following
	properties.
	\begin{enumerate}
		\item The parity of the $\ndZ \Pi _d$-degree of $X_u$
			is even for all $u\in \ndN $.
		      \item For all $u,u'\in \ndN $ one has $[X_u,X_{u'}]=0$.
		\item The families $(X_u)_{u\in \ndN }$,
			$(Z^+_v)_{v\in \ndZ }$, and
			$(Z^-_v)_{v\in \ndZ }$
			satisfy the following equations.
	\begin{equation}
	  \label{eq:XZcomm}
	\begin{split}
		[X_u,Z^\pm _v]=&\epsilon (c^{\pm 1}X_{u-1}Z^\pm _{v\mp 1}
		-c^{\mp 1}Z^\pm _{v\mp 1}X_{u-1})\quad
		\text{for all $u\in \ndN $, $v\in \ndZ $,}\\
		[X_1,Z^+_v]=&rZ^+_{v-1},\quad [X_1,Z^-_v]=-rZ^-_{v+1}
		\quad \text{for all $v\in \ndZ $.}
	\end{split}
	\end{equation}
	\end{enumerate}
	Let $b\in \ndC \setminus \{0\}$ such that $rb=\epsilon (c-c^{-1})$.
	For all $u\in \ndN $ define
	${\cal L}_u\in U'_d$ by the following generating function in $z$.
	\begin{align}
	  \exp\Big(b\sum_{u=1}^\infty {\cal L}_uz^u\Big)=
	  1+b\sum_{u=1}^\infty X_uz^u.
	  \label{eq:L1}
	\end{align}
	Then for all $u\in \ndN $ and $v\in \ndZ $ we have
	\begin{align}
	  [{\cal L}_u,Z^\pm_v]=\epsilon^u{\frac {c^{\pm u}-c^{\mp u}} {ub}}
	  Z^\pm_{v\mp u}.
	  \label{eq:L2}
	\end{align}
\end{lemma}

\begin{proof}
  We proceed by induction on $u$. For $u=1$ one has ${\cal L}_1=X_1$,
  and hence the second line in \eqref{eq:XZcomm} together with the
  definition of $b$ implies the claim.

  Let $n\in \ndN $ and assume that Eq.~\eqref{eq:L2} holds for $u<n$.
  Let $\cL (z)=\sum _{u=1}^\infty \cL _uz^u$
  and $\cZ ^\pm (y)=\sum _{v=-\infty}^\infty Z^\pm _vy^v$,
  where $y$ is a formal parameter, and let
  \begin{align*}
   \cP (y):=&[{\cal L}_n,\cZ ^\pm (y)]-
   \epsilon^n{\frac {c^{\pm n}-c^{\mp n}}{nb}}y^{\pm n}\cZ ^\pm (y)\\
   =&\sum _{v=-\infty }^\infty \big([\cL _n,Z^\pm _v]
   -\epsilon^n{\frac {c^{\pm n}-c^{\mp n}}{nb}}Z^\pm _{v\mp u}\big)y^v.
  \end{align*}
  
In the remainder of this proof
we treat equations in the algebra $U'_d\wlbr z\wrbr /(z^{n+1})$. 
We have the following:
\begin{align*}
  [{\cal L}_u,{\cal Z}^\pm (y)]=&\sum _{v=-\infty }^\infty \epsilon ^u
  \frac{c^{\pm u}-c^{\mp u}}{ub}Z^\pm _{v\mp u}y^v
  =y^{\pm u}\epsilon ^u\frac{c^{\pm u}-c^{\mp u}}{ub}{\cal Z}^\pm (y),\\
  \intertext{for all $u\in \ndN $ with $u<n$ by induction hypothesis,
  and hence}
  [b\cL (z),\cZ ^\pm (y)]=&b\cP (y)z^n
  +\sum _{u=1}^\infty \epsilon ^u\frac{c^{\pm u}-c^{\mp u}}{u}y^{\pm u}z^u
  \cZ ^\pm (y).
\end{align*}
This gives
\begin{align*}
  &(\exp b\cL (z))\cZ ^\pm (y)(\exp b\cL (z))^{-1}=\exp (\mathrm{ad}\,
  b\cL (z))(\cZ ^\pm (y))\\
  &\quad =\cZ ^\pm (y)+b\cP (y)z^n+\sum _{u=1}^\infty
  \epsilon ^u\frac{c^{\pm u}-c^{\mp u}}{u}y^{\pm u}z^u\cZ ^\pm (y)\\
  &\qquad \quad +\sum _{v=2}^\infty \frac{1}{v!}\Big(
  \sum _{u=1}^\infty \epsilon ^u\frac{c^{\pm u}-c^{\mp u}}{u}
  y^{\pm u}z^u\Big)^v\cZ ^\pm (y)\\
  &\quad =b\cP (y)z^n+\exp \Big(
  \sum _{u=1}^\infty \frac{1}{u}\big((\epsilon c^{\pm 1}y^{\pm 1}z)^u-
  (\epsilon c^{\mp 1}y^{\pm 1}z)^u\big)\Big)\cZ ^\pm (y)\\
  &\quad =b\cP (y)z^n+\exp \big(-\log (1-\epsilon c^{\pm 1}y^{\pm 1}z)
  +\log (1-\epsilon c^{\mp 1}y^{\pm 1}z)\big)
  \cZ ^\pm (y)\\
  &\quad =b\cP (y)z^n+(1-\epsilon c^{\mp 1}y^{\pm 1}z)
  \sum _{v=0}^\infty (\epsilon c^{\pm 1}y^{\pm 1}z)^v \cZ ^\pm (y)\\
  &\quad =b\cP (y)z^n+\cZ ^\pm (y)+
  (1-c^{\mp 2})\sum _{v=1}^\infty (\epsilon c^{\pm 1}y^{\pm 1})^vz^v
  \cZ ^\pm (y).
\end{align*}

Next we compute $[\exp b\cL (z),\cZ ^\pm (y)]$ in two different ways.
By Lemma~\ref{le:XZseries} we obtain that
\begin{align*}
  [X_u,\cZ^\pm (y)]=&\sum _{v=-\infty }^\infty [X_u,Z^\pm _v]y^v
  =\sum _{v=-\infty}^\infty \Big((\epsilon c^{\pm 1})^{u-1}
  [X_1,Z^\pm _{v\mp (u-1)}]\\
  &+\epsilon (c^{\pm 1}-c^{\mp 1})\sum _{m=1}^{u-1}
  (\epsilon c^{\pm 1})^{u-1-m}Z^\pm _{v\mp (u-m)}X_m\Big)y^v,
  \intertext{and using the second line in Equation~\eqref{eq:XZcomm}
  we get}
  [X_u,\cZ^\pm (y)]=&\pm r(\epsilon c^{\pm 1})^{u-1}y^{\pm u}\cZ^\pm (y)\\
  &+\epsilon (c^{\pm 1}
  -c^{\mp 1})\sum _{m=1}^{u-1}(\epsilon c^{\pm 1})^{u-1-m}
  y^{\pm (u-m)}\cZ ^\pm (y)X_m.
\end{align*}
Thus we have
\begin{align*}
  &[\exp b\cL (z),\cZ ^\pm (y)]=b\sum _{u=1}^\infty [X_u,\cZ ^\pm (y)]z^u\\
  &
\quad 
=\pm br\sum _{u=1}^\infty (\epsilon c^{\pm 1})^{u-1}y^{\pm u}z^u
  \cZ ^\pm (y)
\\&\qquad \quad 
  +b\epsilon (c^{\pm 1}-c^{\mp 1})\sum _{u=1}^\infty \sum _{m=1}^{u-1}
  (\epsilon c^{\pm 1})^{u-1-m}y^{\pm (u-m)}z^u\cZ^\pm (y)X_m,
  \intertext{and also}
  &[\exp b\cL (z),\cZ ^\pm (y)]=(\exp b\cL (z))\cZ^\pm (y)(\exp b\cL (z))^{-1}
  \exp b\cL (z)
%\\&\qquad \quad 
-\cZ ^\pm (y)\exp b\cL (z)\\
  &\quad =(b\cP (y)z^n+\cZ ^\pm (y)+
  (1-c^{\mp 2})\sum _{v=1}^\infty (\epsilon c^{\pm 1}y^{\pm 1})^vz^v
  \cZ ^\pm (y))\exp b\cL (z)\\
  &\qquad \quad -\cZ ^\pm (y)\exp b\cL (z)\\
  &\quad =b\cP (y)z^n+(1-c^{\mp 2})
  \sum _{v=1}^\infty (\epsilon c^{\pm 1}y^{\pm 1})^vz^v
  \cZ ^\pm (y)(1+b\sum _{u=1}^\infty X_uz^u)
\end{align*}
Comparison of the two expressions for $[\exp b\cL (z),\cZ ^\pm (y)]$
using the formula $\pm br=\epsilon (c^{\pm 1}-c^{\mp 1})$ gives that
$b\cP (y)z^n=0$, and hence $\cP (y)=0$.
This gives the statement of the lemma.
\end{proof}

Next we calculate commutation relations between root vectors to real
roots and type two imaginary root vectors.
Let
\begin{align}
  \Theta (k)=
  \begin{cases}
    0 & \text{if $k\le 0$,}\\
    1 & \text{if $k>0$,}
  \end{cases}
  \label{eq:Heavi}
\end{align}
denote the Heaviside function.

\begin{lemma} \label{le:hcomm}
  Let $d\in \sd \setminus \{0\}$.

  {\rm{(1)}} For all $i,j\in \isv \setminus \{0\}$ and
  $k\in \ndZ \setminus \{0\}$, $m\in \ndZ $ one has
  \begin{align*}
    [\bh_{i,k;d},\cT(\tomega^\vee_{j,d})^m(E_{j,d})] &=\epsilon _{i,j;d}^k
    \frac{[k(\al _{i,d},\al _{j,d})]_q}{k}K_{\isor _d;d}^{\Theta( -k) k}
    \cT(\tomega^\vee_{j,d})^{m-k}(E_{j,d}), \\
    [\bh_{i,k;d},\cT(\tomega^\vee_{j,d})^m(F_{j,d})] &=-\epsilon _{i,j;d}^k
    \frac{[k(\al _{i,d},\al _{j,d})]_q}{k}K_{\isor _d;d}^{\Theta(k)k}
    \cT(\tomega^\vee_{j,d})^{m+k}(F_{j,d}).
  \end{align*}

  {\rm{(2)}} For all $i,j\in \isv \setminus \{0\}$ and $k,l\in
  \ndZ \setminus \{0\}$ one has
\begin{align*}
[\bh_{i,k;d},\bh_{j,l;d}]=\delta_{k,-l}(-1)^k\epsilon _{i,j;d}^k
{\frac {[k(\al _{i,d},\al _{j,d})]_q} k}
{\frac {K_{\isor _d;d}^k-K_{\isor _d;d}^{-k}} {q-q^{-1}}}
\end{align*}
\end{lemma}

\begin{proof}
Part~(1) follows immediately from Lemmata~\ref{le:Tomega(bpsi)},
\ref{le:psiTomegaEF}, 
%\ref{pr:allpsis} and
\ref{le:L} 
and Proposition~\ref{pr:allpsis}. 

To part (2). If $k>l>0$ or $k<l<0$, then the 
equation $[\bh_{i,k;d},\bh_{j,l;d}]=0$
holds by Preposition~\ref{pr:allpsis}(2).
Assume now that $k>0>l$ and let $m=-l$.
By 
Definitions~\ref{de:t1imrv}, \ref{de:t0imrv} and 
Preposition~\ref{pr:allpsis}(1),
we have for all $s\in \ndZ $ the equation
\begin{align*}
  \bpsi _{j,l;d}=\Psi(\bpsi _{j,m;d}^{(s)})=(-1)^m K_{\isor _d;d}
  ^{-(s+m)}K_{j,d}(-1)^{p(\al _{j,d})}
  [\cT(\tomega^\vee_{j,d})^s(F_{j,d}),\cT(\tomega^\vee_{j,d})^{s+m}(E_{j,d})].
\end{align*}
By part (1) of the lemma and by Proposition~\ref{pr:allpsis}(1) we have
\begin{eqnarray*}
\lefteqn{[\bh_{i,k;d}, {\bpsi}_{j,l;d}]
=[\bh_{i,k;d}, (-1)^lK_{\isor _d;d}^lK_{j,d}
(-1)^{p(\al _{j,d})}
[F_{j,d},\cT(\tomega^\vee_{j,d})^m(E_{j,d})]]} \\
&=&(-1)^{m+p(\al _{j,d})}K_{\isor _d;d}^lK_{j,d}
\epsilon _{i,j;d}^k{\frac{[k(\al _{i,d},\al _{j,d})]_q} k}\big\{ \\
& & 
[-K_{\isor _d;d}^k\cT(\tomega^\vee_{j,d})^k(F_{j,d}),\cT(\tomega^\vee_{j,d})^m
(E_{j,d})] 
+[F_{j,d},\cT(\tomega^\vee_{j,d})^{m-k}(E_{j,d})]\big\}
\\ 
&=& 
\left\{\begin{array}{ll}
0 & \mbox{if $m<k$,} \\
\displaystyle
(-1)^{k+1}\epsilon _{i,j;d}^k{\frac
{[k(\al _{i,d},\al _{j,d})]_q} k}
{\frac {K_{\isor _d;d}^k-K_{\isor _d;d}^{-k}} {q-q^{-1}}} & \mbox{if $m=k$,}
\\
\displaystyle
(-1)^{k+1}\epsilon _{i,j;d}^k{\frac
{[k(\al _{i,d},\al _{j,d})]_q} k}
(K_{\isor _d;d}^k-K_{\isor _d;d}^{-k})
{\bpsi}_{d,k-m;d} & \mbox{if $m>k$.}
\end{array}\right.
\end{eqnarray*}
Then part (2) of the lemma can be shown
along the lines of the last part of the proof of Lemma~\ref{le:L}.
\end{proof}

The commutation relations in the following two lemmata 
will also be needed in the
last section.

\begin{lemma}\label{le:tcomm1} Let $d\in \sd \setminus \{0\}$ and
$i,j\in \isv \setminus \{0\}$. 
For all $k$, $l\in \ndZ $, one has
\begin{align*}
&\wlbr \cT(\tomega^\vee_{i,d})^k(E_{i,d}),\cT(\tomega^\vee_{j,d})^l(E_{j,d})\wrbr 
= (-1)^{\delta_{i,j}}\wlbr \cT(\tomega^\vee_{j,d})^{l+1}(E_{j,d}),
\cT(\tomega^\vee_{i,d})^{k-1}(E_{i,d})\wrbr, \\
&\wlbr \cT(\tomega^\vee_{i,d})^k(F_{i,d}),\cT(\tomega^\vee_{j,d})^l(F_{j,d})\wrbr 
= (-1)^{\delta_{i,j}}\wlbr \cT(\tomega^\vee_{j,d})^{l+1}
(F_{j,d}),
%(\tomega^\vee_{j,d})(F_{j,d}),
\cT(\tomega^\vee_{i,d})^{k-1}(F_{i,d})\wrbr .
\end{align*}
\end{lemma} 

\begin{proof}
First, treat the first equation. If $i\not=j$ then the statement
follows from Lemma~\ref{le:propertyTomega},
Theorem~\ref{th:tzerotone},
and Theorem~\ref{th:main}(1).
Assume now that $i=j$. Without loss of generality,
by Lemma~\ref{le:propertyTomega}(1),
Theorem~\ref{th:tzerotone},
and Theorem~\ref{th:main}(1),
let $k\ge 1$
and $l=0$.
If $k=1$ then the statement of the lemma follows from
Lemma~\ref{le:TomegaFKone}(2).
Let $l\in \isv \setminus \{0\}$ be such that
$(\al _{i,d},\al _{l,d})\ne 0$. We proceed by induction on $k$.
By Lemma~\ref{le:hcomm}(1) one has
\begin{eqnarray*}
0&=&
{\frac {\epsilon _{i,l;d}} {[(\al _{i,d},\al _{l,d})]_q}}
K_{\isor _d;d}^{-1}[\bh_{l,-1;d},\wlbr \cT(\tomega^\vee_{i,d})^k(E_{i,d}),
E_{i,d}\wrbr \\
& & +\wlbr \cT(\tomega^\vee_{i,d})(E_{i,d}),\cT(\tomega^\vee_{i,d})^{k-1}
(E_{i,d})\wrbr ] \\
&=& \wlbr T(\tomega^\vee_{i,d})^{k+1}(E_{i,d}),E_{i,d}\wrbr
+\cT(\tomega^\vee_{i,d})(
\wlbr \cT(\tomega^\vee_{i,d})^{k-1}(E_{i,d}),E_{i,d}\wrbr \\
& &+\wlbr \cT(\tomega^\vee_{i,d})(E_{i,d}),
\cT(\tomega^\vee_{i,d})^{k-2}(E_{i,d})\wrbr )
+\wlbr \cT(\tomega^\vee_{i,d})(E_{i,d}),T(\tomega^\vee_{i,d})^k(E_{i,d})\wrbr \\
&=&\wlbr \cT(\tomega^\vee_{i,d})^{k+1}(E_{i,d}),E_{i,d}\wrbr
+\wlbr \cT(\tomega^\vee_{i,d})(E_{i,d}),
\cT(\tomega^\vee_{i,d})^k(E_{i,d})\wrbr ,
\end{eqnarray*} as desired.
The second equation is obtained from the first one by applying
$\Psi _d$. 
\end{proof}

\begin{lemma} \label{le:tcomm2}  Let $d\in \sd \setminus \{0\}$ and
$i,j\in \isv \setminus \{0\}$. 
 
{\rm{(1)}} If $d\in\{1,2,3\}$ and $i\not=d$, then for all $k$, $r$, $l\in \ndZ $,
one has
\begin{align*}
&\wlbr \cT (\tomega^\vee_{i,d})^k(E_{i,d}), \wlbr \cT (\tomega^\vee_{i,d})^r(E_{i,d}),
\cT (\tomega^\vee_{j,d})^l(E_{j,d}) \wrbr\wrbr
+ \\
&\quad \wlbr \cT (\tomega^\vee_{i,d})^r(E_{i,d}), \wlbr \cT (\tomega^\vee_{i,d})^k(E_{i,d}),
 \cT (\tomega^\vee_{j,d})^l(E_{j,d}) \wrbr\wrbr 
=0, \\ 
&\wlbr \cT (\tomega^\vee_{i,d})^k(F_{i,d}), \wlbr \cT (\tomega^\vee_{i,d})^r(F_{i,d}),
\cT (\tomega^\vee_{j,d})^l(F_{j,d}) \wrbr\wrbr
+ \\
&\quad \wlbr \cT (\tomega^\vee_{i,d})^r(F_{i,d}), \wlbr \cT (\tomega^\vee_{i,d})^k(F_{i,d}),
\cT (\tomega^\vee_{j,d})^l(F_{j,d}) \wrbr\wrbr 
=0.
\end{align*}

{\rm{(2)}} If $(\al_{i,d},\al_{j,d})=0$, then for all $k\in \ndZ $,
one has
\begin{align*}
&\wlbr \cT (\tomega^\vee_{i,d})^k(E_{i,d}),\cT (\tomega^\vee_{j,d})^l(E_{j,d})\wrbr 
=0, \\
&\wlbr \cT (\tomega^\vee_{i,d})^k(F_{i,d}),\cT (\tomega^\vee_{j,d})^l(F_{j,d})\wrbr 
=0.
\end{align*}

{\rm{(3)}} If $d=4$, then for all $k$, $r$, $l\in \ndZ $,
one has
\begin{align*}
&[(\al_{1,4},\al_{3,4})]_q \wlbr\wlbr 
\cT (\tomega^\vee_{1,d})^k(E_{1,4}),\cT (\tomega^\vee_{2,d})^r(E_{2,4})
\wrbr ,
\cT (\tomega^\vee_{3,d})^l(E_{3,4}) \wrbr
+ \\
&\quad [(\al_{1,4},\al_{2,4})]_q \wlbr\wlbr 
\cT (\tomega^\vee_{1,4})^r(E_{1,4}), \cT (\tomega^\vee_{3,4})^k(E_{3,4})
\wrbr ,
\cT (\tomega^\vee_{2,4})^l(E_{2,4}) \wrbr 
=0, \\ 
&[(\al_{1,4},\al_{3,4})]_q \wlbr\wlbr 
\cT (\tomega^\vee_{1,4})^k(F_{1,4}),\cT (\tomega^\vee_{2,4})^r(F_{2,4})
\wrbr ,
\cT (\tomega^\vee_{3,4})^l(F_{3,4}) \wrbr
+ \\
&\quad [(\al_{1,4},\al_{2,4})]_q \wlbr\wlbr 
\cT (\tomega^\vee_{1,4})^r(F_{1,4}), \cT (\tomega^\vee_{3,4})^k(F_{3,4})
\wrbr ,
\cT (\tomega^\vee_{2,4})^l(F_{2,4}) \wrbr 
=0.
\end{align*} 
\end{lemma}
\begin{proof} These equations are obtained
from the equations $X=0$ in $U'_d$, where 
$X$ 
%are the same elements as those 
are the elements having the same expressions
as those
in \eqref{eq:Ur6}-\eqref{eq:Ur11},
by applying $\cT (\tomega^\vee_{u,d})^m$
and the $\ndC$-linear map ${\rm{ad}}\,\bh_{u,m;d}$
defined by ${\rm{ad}}\,\bh_{u,m;d}(Y)=[\bh_{u,m;d},Y]$ and using
Lemma~\ref{le:propertyTomega}(1),
Theorem~\ref{th:tzerotone},
and Theorem~\ref{th:main}(1).
\end{proof}

\section{Second realization of the quantum affine superalgebras}
\label{sec:secreal}
\subsection{Main theorem for $U'_d$}
Here for each $d\in\sd \setminus\{0\}
=\{1,2,3,4\}$, we introduce Drinfeld second realization associated with
$D^{(1)}(2,1;x)$
and prove that it is isomorphic to $U'_d$ as a $\ndC$-algebra.
We first give a modified version of the Drinfeld second realization of
$U'_d$.
Then via the version we give 
the Drinfeld second realization
of $U'_d$.

\begin{defin} \label{definition:defUD}
Let $d\in\sd \setminus\{0\}=\{1,2,3,4\}$.
Let $DU'_d=\bigoplus _{\mu \in \ndZ \Pi _d}DU'_{d,\mu }$
be the $\ndZ \Pi _d$-graded $\mathbb{C}$-algebra generated by
the elements
\begin{gather}
\sigma _d,\,\,
K_{i,d}^{{\frac 1 2}},\,\,K_{i,d}^{-{\frac 1 2}} \in DU'_{d,0},\,\,
(i\in \isv ) 
  \label{eq:DG1}\\
  x^\pm _{i,k;d}\in DU'_{d,\pm \alpha _{i,d}+k\isor _d},  
(i\in \isv \setminus \{0\},\,k\in\mathbb{Z}) \label{eq:DG2dash}\\
\psi _{i,r;d},\,\,h_{i,r;d}\in DU'_{d,r\isor _d},\,\,
(i\in \isv \setminus \{0\},\,r\in\mathbb{Z}\setminus \{0\}) 
\label{eq:DG2}
\end{gather}
and defined by the relations below, where the elements $K_{i;d}$
and $K_{\isor _d;d}$
are defined as in Eq.~\eqref{eq:isoK}.
\begin{gather}
\sigma _d^2=1,\,\,
K_{i,d}^{{\frac 1 2}}K_{i,d}^{-{\frac 1
2}}=K_{i,d}^{-{\frac 1 2}}
K_{i,d}^{{\frac 1 2}}=1, \label{eq:DRL1} \\
XY=YX \quad \mbox{for all $X$, $Y$ in Eq.~\eqref{eq:DG1}} 
\label{eq:DRL2} \\
  \sigma _d X \sigma _d=(-1)^{p(\mu )}X,\,\,
K_{i,d}^{{\frac 1 2}}XK_{i,d}^{-{\frac 1 2}}
=q^{{\frac {(\al _{i,d},\mu )} 2}}X \,\,\mbox{for all $X\in DU'_{d,\mu }$,
$\mu \in \ndZ \Pi _d$,} \label{eq:DRL3}
\\
[x^+_{i,k;d},x^-_{j,l;d}]  =
\begin{cases}
  0 & \text{if $i\not=j$,}\\
  K_{\isor _d;d}^{-l}K_{i,d}\psi_{i,k+l;d} & \text{if $i=j$ and $k+l>0$,}\\
  \displaystyle
{\frac {K_{\isor _d;d}^kK_{i,d}
-K_{\isor _d;d}^{-k}K_{i,d}^{-1}} {q-q^{-1}}} & \text{if $i=j$ and $k+l=0$,}
\\
-K_{\isor _d;d}^{-k}K_{i,d}^{-1} \psi_{i,k+l;d} & \text{if $i=j$ and $k+l<0$,}
\end{cases} \label{eq:xpxm} \\
\exp\Big(\pm (q-q^{-1})\sum_{k>0}z^kh_{i,\pm k;d}\Big)
=1+(q-q^{-1})\sum_{k>0}z^k\psi_{i,\pm k;d} \label{eq:expexp}
\end{gather} (as equations of generating functions in $z$)
\begin{align}
[h_{i,k;d},h_{j,l;d}] & =
\delta_{k,-l}{\frac {[k(\al_{i,d},\al_{j,d})]_q} k}
{\frac {K_{\isor _d;d}^k-K_{\isor _d;d}^{-k}}{q-q^{-1}}} \label{eq:hikdhjkd}
\\
[h_{i,k;d},x^\pm_{j,l;d}] & =
\pm{\frac {[k(\al_{i,d},\al_{j,d})]_q} k}
K_{\isor _d;d}^{\Theta(\mp k)k}x^\pm_{j,k+l;d} \label{eq:hikdxjkd}
\\
\wlbr x^\pm_{i,k;d}, x^\pm_{j,l;d} \wrbr & =
\begin{cases}
0 & \quad \text{if $(\al _{i,d},\al _{j,d})=0$,} \\ 
-\wlbr x^\pm_{j,l\mp 1;d}, x^\pm_{i,k\pm 1;d} \wrbr &
 \quad \text{if $d\in \{1,2,3\}$,} \\
 \wlbr x^\pm_{j,l\mp 1;4}, x^\pm_{i,k\pm 1;4} \wrbr &
 \quad \text{if $d=4$} 
\end{cases} \label{eq:prc1}
\\ \wlbr x^\pm_{i,r;d}, \wlbr x^\pm_{i,k;d},
 x^\pm_{j,l;d}\wrbr\wrbr
& +
\wlbr x^\pm_{i,k;d}, \wlbr x^\pm_{i,r;d},
 x^\pm_{j,l;d}\wrbr\wrbr
=0 \,\,\text{if $d\in\{1,2,3\}$ and $i\not=d$,} 
\label{eq:prc2}
\end{align}\begin{align}
&[(\al_{1,4},\al_{3,4})]_q
\wlbr \wlbr x^\pm_{1,r;4}, x^\pm_{2,k;4}\wrbr,
x^\pm_{3,l;4}\wrbr
-[(\al_{1,4},\al_{2,4})]_q
\wlbr \wlbr x^\pm_{1,r;4}, x^\pm_{3,k;4}\wrbr,
 x^\pm_{2,l;4}\wrbr=0\,\text{if 
$d=4$.} \label{eq:prc3}
\end{align}

Note that if $i\in \isv \setminus \{0\}$
then in 
%$DU'_{d,0},$
$DU'_d$,
one has 
$\wlbr x^\pm _{i,k;d},x^\pm _{i,k^\pm 1;d}\wrbr =0$
%$\wlbr x^\pm _{i,k;d},x^\pm _{i,k^\mp 1;d}\wrbr =0$ 
for all 
$k\in \ndZ$.
Note that if $i\in \isv \setminus \{0\}$ and $p(\al _{i,d})=1$ 
then in 
%$DU'_{d,0},$
$DU'_d$,
one has $(x^\pm _{i,k;d})^2=\frac{1}{2}
\wlbr x^\pm _{i,k;d},x^\pm _{i,k;d}\wrbr =0$
for all 
$k\in \ndZ$.
\end{defin}

\begin{proposition}\label{pr:exFd}
Let $d\in\sd \setminus\{0\}$. Then
there exists a unique 
$\mathbb{C}$-algebra homomorphism
\begin{align}
\cF_d: DU'_d\to U'_d
		\label{eq:cF}
\end{align} such that for all $u\in I$,
$i\in I\setminus\{0\}$,
$k\in\mathbb{Z}$, and $r\in\mathbb{Z}\setminus\{0\}$,
one has
\begin{align}
			&\cF_d (\sigma _d)=\sigma _d, \quad
			\cF_d  (K_{u,d}^{\pm{\frac 1 2}})=K_{u,d}^{\pm{\frac 1 2}},
\label{eq:cFg1} \\
&\cF_d (x^+ _{i,k;d})=
(-\epsilon _{i,i;d})^k\cT (\tomega^\vee_{i,d})^{-k}(E_{i,d}), 
\label{eq:cFg2} \\
&\cF_d (x^- _{i,k;d})=
(-\epsilon _{i,i;d})^k\cT (\tomega^\vee_{i,d})^k(F_{i,d}), 
\label{eq:cFg3} \\
&\cF_d (h _{i,r;d})=\epsilon _{i,i;d}^r\bh_{i,r;d}, \quad
\cF_d (\psi _{i,r;d})=\epsilon _{i,i;d}^r\bpsi _{i,r;d}.
\label{eq:cFg4}
\end{align}
\end{proposition}
\begin{proof}
This follows immediately
from Definitions~\ref{de:t0imrv}, \ref{de:t1imrv}, \ref{de:t3imrv},
Proposition~\ref{pr:allpsis}, and Lemmata~\ref{le:hcomm},
\ref{le:tcomm1}, and \ref{le:tcomm2}.
\end{proof}
Let $d\in\sd \setminus\{0\}$.
Define $D\Psi_d$ to be the automorphism of $DU'_d$
such that $D\Psi_d(\sigma _d)=\sigma _d$,
$D\Psi_d(K_{u,d}^{\pm{\frac 1 2}})=K_{u,d}^{\mp{\frac 1 2}}$,
$D\Psi_d(x^\pm _{i,k;d})=(\mp 1)^{p(\al_{i,d})}x^\mp _{i,-k;d}$,
$D\Psi_d(\psi _{i,r;d})=\psi _{i,-r;d}$, and 
$D\Psi_d(h _{i,r;d})=-h _{i,-r;d}$.
Note that \begin{align}\label{eq:PDP=DPP}
\cF_dD\Psi_d=\Psi_d\cF_d. 
\end{align}
Define the two elements $\ndX _d^{\pm}\in DU'_{d,\pm\al_{0,d}}$ by 
\begin{align}\label{eq:defndX}
\ndX _d^{\pm}=
\begin{cases}
-r_{d,j;d}r_{d,k;d}r_{d,0;d}K_{0,d}^{\pm 1}
\wlbr \wlbr \wlbr x^\mp _{d,\pm 1;d}, x^\mp _{j,0;d}\wrbr, x^\mp _{k,0;d}\wrbr, x^\mp _{d,0;d}\wrbr \\
\quad \text{if $d\in \{1,2,3\}$, where $\{j,k,d\}=\{1,2,3\}$ and $j<k$,} \\
\pm [(\al_{1,4},\al_{2,4})]_q^{-1}
K_{0,4}^{\pm 1}\wlbr \wlbr x^\mp _{1,\pm 1;4}, x^\mp _{2,0;4}\wrbr, x^\mp _{3,0;4}\wrbr,
\quad  \text{if $d=4$.}
\end{cases}
\end{align}
\begin{lemma}\label{le:isofu1}
Let $d\in\sd \setminus\{0\}$.
One has
$\cF_ d (\ndX _d^+)=E_{0,d}$ and 
%$\cF_ d (\mathbb{X}_d^+)=F_{0,d}$.
$\cF_ d (\mathbb{X}_d^-)=F_{0,d}$.
\end{lemma}
\begin{proof} Assume that $d=4$. 
By Lemma~\ref{le:relbrwbr}, for all $i$, $j\in I=\{0,1,2,3\}$ and all
$Y_\lambda
\in U'_{4,\lambda}$, one has
$\wlbr E_{i,4},F_{j,4}K_{j,4}^{-1}\wrbr
=\delta_{ij}{\frac {1-K_{i,4}^{-2}} {q-q^{-1}}}$, and
\begin{align}\label{eq:le531}
\wlbr \wlbr E_{i,4},Y_\lambda\wrbr , F_{j,4}K_{j,4}^{-1}\wrbr
&=\wlbr E_{i,4},\wlbr Y_\lambda ,F_{j,4}K_{j,4}^{-1}\wrbr\wrbr
+\delta_{ij}(-1)^{p(\lambda)p(\al_{i,4})}[(\lambda,\al_{i,4})]_qY_\lambda .
\end{align}
Note that for $X_\mu
\in U'_{4,\mu}$ and $Y_\lambda
\in U'_{4,\lambda}$, one has 
\begin{align}\label{eq:le532}
\wlbr X_\mu K_{\mu;4} ,Y_\lambda K_{\lambda;4} \wrbr =q^{(\al_{\mu,4},\al_{\lambda,4})}
\wlbr X_\mu ,Y_\lambda \wrbr K_{\mu + \lambda;4} .
\end{align}
By Eq.s \eqref{eq:isoroot}, \eqref{eq:repomega} and \eqref{eq:T3}, one has
\begin{align}\label{eq:le533}
\cT (\tomega^\vee_{1,4})(K_{1,4}^{-1})=K_{0,4}K_{2,4}K_{3,4} .
\end{align}
Then one has
\begin{align*}
\cF_4 (\ndX _4^+)&={\frac {K_{0,4} 
\wlbr \wlbr - \cT (\tomega^\vee_{1,d})(F_{1,4}), 
F_{2,4}\wrbr, F_{3,4}\wrbr} {[(\al_{1,4},\al_{2,4})]_q}} 
\quad \text{(by Eq.s~\eqref{eq:cFg3} and \eqref{eq:defndX})} \\
&=-[(\al_{1,4},\al_{2,4})]_q^{-1}\wlbr \wlbr \cT (\tomega^\vee_{1,4})(F_{1,4}K_{1,4}^{-1}), 
F_{2,4}K_{2,4}^{-1}\wrbr, F_{3,4}K_{3,4}^{-1}\wrbr \\
&\quad \text{(by Eq.s~\eqref{eq:le532}, \eqref{eq:le533}, and \eqref{eq:cUr4})} \\
&={\frac {\wlbr \wlbr 
\wlbr E_{3,4},\wlbr E_{2,4}, E_{0,4} 
\wrbr \wrbr, 
F_{2,4}K_{2,4}^{-1}\wrbr, F_{3,4}K_{3,4}^{-1}\wrbr} {[(\al_{1,4},\al_{2,4})]_q[(\al_{0,4},\al_{2,4})]_q}}
\quad \text{(by Lemma~\ref{le:TomegaFKone}(1))} \\
&=-[(\al_{1,4},\al_{2,4})]_q^{-1}\wlbr 
\wlbr E_{3,4},E_{0,4} \wrbr,
F_{3,4}K_{3,4}^{-1}\wrbr\quad \text{(by Eq.~\eqref{eq:le531})} \\
&=E_{0,4}\quad \text{(by Eq.~\eqref{eq:le531})}.
\end{align*}
Thus one has $\cF_4 (\ndX _4^+)=E_{0,4}$.
By this and Eq.~\eqref{eq:PDP=DPP}, 
one obtains 
\begin{align*}
\cF_ 4 (\ndX _4^-)=\cF_ 4 (-D\Psi_d(\ndX _4^+))
=-\Psi_d(\cF_ 4 (\ndX _4^+))=-\Psi_d(E_{0,4})=F_{0,4},
\end{align*} as desired.

%By Eq.~\eqref{eq:defndX} and
%Proposition~\ref{pr:exFd}, one has
%\begin{align*}
%\cF_d (\ndX _d^+)
%=
%-[(\al_{1,4},\al_{2,4})]_q^{-1}\wlbr \wlbr \cT (\tomega^\vee_{1,d})(F_{1,d}K_{1,d}^{-1}), 
%F_{2,d}K_{2,d}^{-1}\wrbr, F_{3,d}K_{3,d}^{-1}\wrbr.
%\end{align*} 
%Lemma~\ref{le:relbrwbr} implies that for all $i$, $j\in I=\{0,1,2,3\}$ and all
%$Y_\lambda
%\in U'_{4,\lambda}$ one has 
%\begin{align*}
%\wlbr E_{i,4},F_{j,4}K_{j,4}^{-1}\wrbr
%&=\delta_{ij}(q-q^{-1})^{-1}(1-K_{i,4}^{-2}), \\
%\wlbr \wlbr E_{i,4},Y_\lambda\wrbr , F_{j,4}K_{j,4}^{-1}\wrbr
%&=\wlbr E_{i,4},\wlbr Y_\lambda ,F_{j,4}K_{j,4}^{-1}\wrbr\wrbr
%+\delta_{ij}(-1)^{p(\lambda)p(\al_{i,4})}[(\lambda,\al_{i,4})]_qY_\lambda .
%\end{align*}
%Lemma~\ref{le:TomegaFKone}(1) and these three equations imply that 
%$\cF_ 4 (\ndX _4^+)=E_{0,4}$. 
%Applying $\Psi_4$ to the last equation, by
%using \eqref{eq:PDP=DPP}, one obtains $\cF_ 4 (\ndX _4^-)=F_{0,4}$.

The lemma for $d\in \{1,2,3\}$
can be proven in an entirely similar way.
\end{proof}

For $i\in I\setminus\{0\}$ and $d\in\sd \setminus\{0\}$,
let $z^\pm _{i,k;d}:=x^\pm  _{i,k;d}K_{i,d}^{\pm 1}K_{\isor _d;d}^k$. 
%Then one has
%\begin{align}
%\wlbr z^\pm _{i,k;d},z^\pm _{j,r;d}\wrbr
%& =q^{(\al_{i,d},\al_{j,d})}
%\wlbr x^\pm _{i,k;d},x^\pm _{j,r;d}\wrbr K_{i,d}^{\pm 1}K_{j,d}^{\pm 1}K_{\isor _d;d}^{k+r},
%\label{eq:zzxx} \\
%[z^+ _{i,k;d},z^- _{j,r;d}] & = 
%q^{-(\al_{i,d},\al_{j,d})}[x^+ _{i,k;d},x^- _{j,r;d}]K_{i,d}K_{j,d}^{- 1}K_{\isor _d;d}^{k+r} 
%\label{eq:zzxx2} \\
%& = \begin{cases}
%q^{-(\al_{i,d},\al_{j,d})}
%{\frac {K_{i,d}K_{\isor _d;d}^k-K_{i,d}^{-1}K_{\isor _d;d}^{-k}} {q-q^{-1}}} 
%\quad\text{if $i=j$ and $r=-k$,}\\
%0 \quad\text{if $i\not=j$}
%\end{cases}
%\nonumber
%\end{align} for all $d\in \sd\setminus\{0\}$, $i,j\in I\setminus\{0\}$,
%and $k,r\in\ndZ$.

To prove the existence of the inverse of $\cF$,
we need the following lemma, which is similar to
Lemma~\ref{le:relbrwbr}.

\begin{lemma}\label{le:DUfms}
Let $\mu ,\lambda ,\xi \in \ndZ \Pi _d$,
$i,j\in I\setminus\{0\}$, and $k,r\in\ndZ$, and let
$X_\mu\in DU '_{d,\mu }$,
$Y_\lambda \in DU '_{d,\lambda }$, and
$Z_\xi \in DU '_{d,\xi }$.

{\rm{(1)}}
One has 
%the same equations as the ones in Lemma~\ref{le:relbrwbr}
%with $DU '_{d,\mu }$, $DU '_{d,\lambda }$, 
%and $DU '_{d,\xi }$ in place of 
%$U '_{d,\mu }$, $U '_{d,\lambda }$, 
%and $U '_{d,\xi }$ respectively. 
	\begin{gather}
		\wlbr X_\mu ,K_{\lambda ;d}Y_\lambda \wrbr =
		q^{-(\lambda ,\mu )}K_{\lambda ;d}[X_\mu ,Y_\lambda ],\qquad
		\wlbr K_{\mu ;d}^{-1} X_\mu ,Y_\lambda \wrbr =
		K_{\mu ;d}^{-1}[X_\mu ,Y_\lambda ],
\label{eq:usfpr1}\\
		\wlbr \wlbr X_\mu ,Y_\lambda \wrbr ,Z_\xi \wrbr =
		\wlbr X_\mu ,\wlbr Y_\lambda ,Z_\xi \wrbr \wrbr 
		+(-1)^{p(\lambda )p(\xi )}q^{-(\lambda ,\xi)}[\wlbr X_\mu ,
		Z_\xi \wrbr ,Y_\lambda ]_{q^{(\lambda ,\xi -\mu )}},
\label{eq:usfpr2}\\
		\wlbr X_\mu ,\wlbr Y_\lambda ,Z_\xi \wrbr \wrbr =
		\wlbr \wlbr X_\mu ,Y_\lambda \wrbr ,Z_\xi \wrbr 
		+(-1)^{p(\mu )p(\lambda )}q^{-(\mu ,\lambda )}[Y_\lambda ,
		\wlbr X_\mu ,Z_\xi \wrbr ]_{q^{(\lambda ,\mu -\xi )}}.
\label{eq:usfpr3}	\end{gather} Further, one has
\begin{align}
%& \wlbr X_\mu ,K_{\lambda ;d}Y_\lambda \wrbr =
%		q^{-(\lambda ,\mu )}K_{\lambda ;d}[X_\mu ,Y_\lambda ], \\
%&\wlbr \wlbr X_\mu ,Y_\lambda \wrbr ,Z_\xi \wrbr =
%		\wlbr X_\mu ,\wlbr Y_\lambda ,Z_\xi \wrbr \wrbr 
%		+(-1)^{p(\lambda )p(\xi )}q^{-(\lambda ,\xi)}[\wlbr X_\mu ,
%		Z_\xi \wrbr ,Y_\lambda ]_{q^{(\lambda ,\xi -\mu )}}, \label{eq:usf0}
% \\
&\wlbr x^+ _{i,k;d},\wlbr Y_\lambda,  
z^- _{j,-k;d} \wrbr \wrbr =\wlbr \wlbr x^+ _{i,k;d},Y_\lambda \wrbr ,
z^- _{j,-k;d} \wrbr - \delta_{ij}(-1)^{p(\al_{j,d})p(\lambda )}
[(\al_{i,d},\lambda )]_q Y_\lambda \label{eq:usf} \\ 
&\wlbr z^\pm _{i,k;d},z^\pm _{j,r;d}\wrbr
=q^{(\al_{i,d},\al_{j,d})}
\wlbr x^\pm _{i,k;d},x^\pm _{j,r;d}\wrbr K_{i,d}^{\pm 1}K_{j,d}^{\pm 1}K_{\isor _d;d}^{k+r},
\label{eq:zzxx} \\
&[z^+ _{i,k;d},z^- _{j,r;d}]  = 
q^{-(\al_{i,d},\al_{j,d})}[x^+ _{i,k;d},x^- _{j,r;d}]K_{i,d}K_{j,d}^{- 1}K_{\isor _d;d}^{k+r} 
\label{eq:zzxx2} \\
& = \begin{cases}
q^{-(\al_{i,d},\al_{j,d})}
{\frac {K_{i,d}K_{\isor _d;d}^k-K_{i,d}^{-1}K_{\isor _d;d}^{-k}} {q-q^{-1}}} 
& \text{if $i=j$ and $r=-k$,}\\
0 & \text{if $i\not=j$}.
\end{cases}
\nonumber
\end{align}

{\rm{(2)}} 
If $p(\lambda)=1$ and $(\lambda,\lambda )=0$,
then one has 
\begin{align}
& \wlbr \wlbr X_\mu ,Y_\lambda \wrbr ,Y_\lambda \wrbr =
\wlbr Y_\lambda , \wlbr Y_\lambda, X_\mu \wrbr  \wrbr =0 \quad
\text{if $Y_\lambda^2=0$}, \label{eq:usf2} \\
& \wlbr X_\mu ,Y_\lambda \wrbr^2=0
 \quad \text{if $Y_\lambda^2=0$, $\wlbr X_\mu ,\wlbr X_\mu ,Y_\lambda \wrbr
\wrbr =0$,
and $q^{(\mu ,\mu )}+(-1)^{p(\mu )}\not= 0$}, \label{eq:usf3} \\
& \wlbr Y_\lambda, X_\mu \wrbr^2=0
 \quad \text{if $Y_\lambda^2=0$, $\wlbr \wlbr 
Y_\lambda, X_\mu \wrbr , X_\mu  \wrbr
=0$,
and $q^{(\mu ,\mu )}+(-1)^{p(\mu )}\not= 0$}. \label{eq:usf4}
\end{align}

{\rm{(3)}} 
If $(\mu,\lambda)+(\lambda,\xi)+(\xi,\mu)=0$
	and $p(\mu)=p(\lambda)=p(\xi)=1$, then one has
	\begin{align*}
	&[(\mu ,\xi )]_q\wlbr \wlbr X_\mu ,Y_\lambda \wrbr ,Z_\xi \wrbr
	-[(\mu ,\lambda )]_q\wlbr \wlbr X_\mu ,Z_\xi \wrbr , Y_\lambda \wrbr \\
	&=[(\mu ,\xi )]_q(X_\mu Y_\lambda Z_\xi -Z_\xi Y_\lambda X_\mu )
    +[(\lambda ,\mu )]_q(Y_\lambda Z_\xi X_\mu - X_\mu Z_\xi Y_\lambda ) \\
    &\quad +[(\xi ,\lambda )]_q(Z_\xi X_\mu Y_\lambda - Y_\lambda X_\mu Z_\xi ) \\
    &=-[(\mu ,\xi )]_q\wlbr Z_\xi, \wlbr Y_\lambda , X_\mu \wrbr\wrbr
	+[(\mu ,\lambda )]_q\wlbr Y_\lambda, \wlbr Z_\xi , X_\mu \wrbr\wrbr .
	\end{align*}
	
{\rm{(4)}} 
%{\rm{(3)}} 
Let $X_{-\mu}\in DU '_{d,-\mu }$ and
$Y_{-\lambda} \in DU '_{d,-\lambda }$.
Assume that $[X_\mu,X_{-\mu}]=
{\frac {K_{\mu;d}-K_{-\mu;d}} {q-q^{-1}}}$,
$[Y_\lambda,Y_{-\lambda}]=
{\frac {K_{\lambda;d}-K_{-\lambda;d}} {q-q^{-1}}}$,
and $[X_{\pm\mu},Y_{\mp\lambda}]=0$.
Then one has 
\begin{align*}
[\wlbr X_\mu ,Y_\lambda\wrbr,\wlbr X_{-\mu} ,Y_{-\lambda}\wrbr ]
=(-1)^{p(\lambda ) p(\mu )}q^{-(\lambda ,\mu )}
[(\lambda ,\mu )]_q{\frac {K_{\mu +\lambda ;d}-K_{-\mu -\lambda ;d}} {q-q^{-1}}},
\end{align*} where $K_{\mu;d}$ and $K_{\lambda;d}$
are defined as in \eqref{eq:Kmu}.

\end{lemma}
\begin{proof}
One needs only the definition of $q$-super-bracket
and Eq.~\eqref{eq:xpxm}.
\end{proof}

We also need the following lemma.

\begin{lemma}\label{le:caiso} 
Let $\{d,i,j\}=\{1,2,3\}$ and $l,r,m\in \ndZ $.
Let $\ndY^\pm_1$ and $\ndY^\pm_2$ be the elements of $DU'_d$ defined by
$\ndY^\pm_1:= \wlbr \wlbr z^\pm _{d,\mp l;d},z^\pm _{j,\mp r;d}\wrbr , 
z^\pm _{k,\mp m;d}\wrbr$
and $\ndY^\pm_2:= \wlbr \ndY^\pm_1 , z^\pm _{d,\mp (l-1);d}\wrbr$. 

{\rm{(1)}} One has the following:
\begin{gather}
(z^- _{d,l;d})^2=0,\quad \wlbr z^- _{j,r;d},z^- _{k,m;d}\wrbr=
\wlbr z^- _{k,m;d}, z^- _{j,r;d}\wrbr=0,\label{eq:caiso1} \\
\wlbr z^- _{d,l;d},z^- _{j,r;d}\wrbr=-\wlbr z^- _{j,r+1;d}, z^- _{d,l-1;d}\wrbr ,
\quad \wlbr \wlbr z^- _{d,l;d},z^- _{k,m;d}\wrbr , z^- _{k,m;d}\wrbr=0,
\label{eq:caiso2} \\
\ndY^-_1= \wlbr \wlbr z^- _{d,l;d},z^- _{k,m;d}\wrbr , z^- _{j,r;d}\wrbr ,\quad
\wlbr \wlbr z^- _{j,r;d}, z^- _{d,l;d}\wrbr , z^- _{d,l;d}\wrbr =0,\,\, 
%\wlbr z^- _{j,r;d}, z^- _{d,l;d}\wrbr^2 =0,
\label{eq:caiso3} \\
\wlbr z^- _{j,r;d}, z^- _{d,l;d}\wrbr^2 =0,\quad 
\wlbr z^- _{d,l;d}, z^- _{j,r;d}\wrbr^2 =0,
\label{eq:caiso3b}  \\
\wlbr \ndY^-_1, z^- _{j,r;d}\wrbr =0, \quad
\wlbr \ndY^-_1, z^- _{k,m;d}\wrbr =0, \quad (\ndY^-_1)^2=0, \quad
\wlbr \ndY^-_1, \ndY^-_2\wrbr =0,
%\quad \wlbr \ndY^-_2, z^- _{d,l-1;d}\wrbr =0, 
\label{eq:caiso4}  \\
\ndY^-_1=- \wlbr \wlbr z^- _{k,m+1;d}, z^- _{d,l-1;d} \wrbr, z^- _{j,r;d} \wrbr
=- \wlbr z^- _{k,m+1;d}, \wlbr z^- _{d,l-1;d}, z^- _{j,r;d} \wrbr\wrbr ,
\label{eq:caiso5aa} \\
%\wlbr \ndY^-_1, z^- _{j,r;d}\wrbr =0,
%\label{eq:caiso5a}
%\\
%\wlbr \ndY^-_2, z^- _{j,r;d}\wrbr =0,  \quad
%\wlbr \ndY^-_2, z^- _{k,m;d}\wrbr =0, \label{eq:caiso5} 
%\\
[\ndY^-_2, x^- _{j,r;d}] =0,  \quad
[\ndY^-_2, x^- _{k,m;d}] =0. \label{eq:caiso5} 
\end{gather}

{\rm{(2)}} One has
\begin{gather}
\wlbr x^+ _{d,-l+1;d},\ndY^-_1 \wrbr =0, \quad
\wlbr x^+ _{d,-l+1;d},\ndY^-_2 \wrbr =-[(\al_{j,d}+\al_{k,d},\al_{d,d})]_q\ndY^-_1 , 
\label{eq:caiso6}\\
\wlbr\wlbr x^+ _{d,-l+1;d},\ndY^-_2 \wrbr ,\ndY^-_2 \wrbr=0,
\quad
\wlbr x^+ _{j,-r;d}, \ndY^-_2 \wrbr=0, 
\quad \wlbr x^+ _{k,-m;d}, \ndY^-_2 \wrbr=0. \label{eq:caiso7}
\end{gather}

{\rm{(3)}} If $l=1$ and $r=m=0$ then one has
\begin{gather}
[\ndY^- _2, \ndY^+ _2]=
[(\al_{d,d},\al_{j,d})]_q[(\al_{d,d},\al_{k,d})]_q[(\al_{d,d},\al_{0,d})]_q
{\frac {K_{0,d}-K_{0,d}^{-1}} {q-q^{-1}}}.\label{eq:caiso8}
\end{gather}
\end{lemma}
\begin{proof} To Part~(1): 
Eq.s~\eqref{eq:caiso1},\eqref{eq:caiso2} hold just by
\eqref{eq:prc1}-\eqref{eq:prc2} and \eqref{eq:zzxx}.
Each equation in \eqref{eq:caiso3}-\eqref{eq:caiso5aa}
follows from the ones in \eqref{eq:caiso1}-\eqref{eq:caiso5aa}
before it
and Lemma~\ref{le:DUfms}(1),(2).
%\begin{align*}
%\ndY_1 & = \wlbr \wlbr z^- _{d,l;d},z^- _{k,m;d}\wrbr , z^- _{j,r;d}\wrbr
%=-\wlbr \wlbr z^- _{k,m+1;d},z^- _{d,l-1;d}\wrbr ,z^- _{j,r;d}\wrbr \\
%&=-\wlbr z^- _{k,m+1;d},\wlbr z^- _{d,l-1;d},z^- _{j,r;d}\wrbr\wrbr , \\
%\wlbr \ndY^- _1, z^- _{j,r;d}\wrbr &=
%- \wlbr z^- _{k,m+1;d}, \wlbr\wlbr z^- _{d,l-1;d}, z^- _{j,r;d}\wrbr , z^- _{j,r;d} \wrbr\wrbr \\
%&\,\quad\text{(by Equation~\eqref{eq:usfpr2}} \\
%&\quad\quad 
%\text{and the second equations of Equation~\eqref{eq:caiso1},\eqref{eq:caiso5a})}\\
%&=0.\,\,\text{(by the second equation of Equation~\eqref{eq:caiso2})}
%\end{align*} 
%Similarly 
One obtains the first equation in
\eqref{eq:caiso5} in the following way.
\begin{align*}
&[\ndY^- _2, x^- _{j,r;d}]K_{j,d}^{-1}K_{\isor _d;d}^r 
= \wlbr \ndY^- _2, z^- _{j,r;d}\wrbr\quad 
\text{(by the equation $[K_{j,d}^{-1}K_{\isor _d;d}^r,\ndY^- _2]=0$)} \\
\quad &=\wlbr\wlbr  \ndY^- _1, z^- _{d,l-1;d}\wrbr, z^- _{j,r;d}\wrbr 
\quad \text{(by the definition of $\ndY^- _2$)}
\\
\quad &=\wlbr \ndY^- _1 , \wlbr z^- _{d,l-1;d}, z^- _{j,r;d} \wrbr\wrbr
%\\ &\, 
\quad \text{(by Eq. \eqref{eq:usfpr2} and
the first equation in \eqref{eq:caiso4})} \\
\quad &=-\wlbr \wlbr z^- _{k,m+1;d},\wlbr z^- _{d,l-1;d},z^- _{j,r;d}\wrbr\wrbr , 
\wlbr z^- _{d,l-1;d}, z^- _{j,r;d} \wrbr\wrbr 
\quad \text{(by using \eqref{eq:caiso5aa})} \\
\quad &=0 \quad \text{(by the second equation in \eqref{eq:caiso3b}
and by using \eqref{eq:usf2})}.
\end{align*} 
Similarly one has the second equation in
\eqref{eq:caiso5}.

To Part~(2): Applying Eq.~\eqref{eq:usf} twice,
one has the first equation in \eqref{eq:caiso6}
by using the first one in \eqref{eq:caiso5aa}
and the second one in \eqref{eq:caiso1}.
Similarly, the second equation in \eqref{eq:caiso6}
holds by the first one in \eqref{eq:caiso6}.
The first equation in \eqref{eq:caiso7} follows from 
the second one in \eqref{eq:caiso6}
and 
%the third one
the fourth one
in \eqref{eq:caiso4}.
As for the second equation in \eqref{eq:caiso7}, one has
\begin{align*}
\wlbr x^+ _{j,-r;d}, \ndY^-_2 \wrbr
%& =[(\al_{d,d},\al_{j,d})]_q\wlbr\wlbr z^- _{d,l;d},z^- _{k,m;d}\wrbr,
%z^- _{d,l-1;d} \wrbr \\
%& \quad\quad\text{(by the definition of $\ndY^-_1$ and Eq.~\eqref{eq:usf})} \\
& =-[(\al_{d,d},\al_{j,d})]_q
\wlbr\wlbr z^- _{k,m+1;d},z^- _{d,l-1;d}\wrbr,
z^- _{d,l-1;d} \wrbr  \\
& \,\quad\text{(by using \eqref{eq:usf} and 
the first equation in \eqref{eq:caiso5aa})}  \\
&= 0\,\,\text{(by the second equation in \eqref{eq:caiso3})}.
\end{align*} Similarly one has the third equation in \eqref{eq:caiso7}.

To Part~(3): Using the second equation in \eqref{eq:usfpr1}
and applying $D\Psi_d$ to the first one in \eqref{eq:caiso6}
one has $[z^\pm _{d,-l+1;d}, Y^\mp _2]=0$. Then 
applying Lemma~\ref{le:DUfms}~(4) and 
Eq.~\eqref{eq:zzxx2} repeatedly one obtains \eqref{eq:caiso8}.
\end{proof}

We first give a modified version
of the Drinfeld second realization of $U'_d$.

\begin{theorem}\label{theorem:iso}
For each $d\in\sd \setminus\{0\}$ the map 
$\cF_d:DU'_d\to U'_d$ is a
$\ndZ\Pi_d$-graded $\ndC$-algebra isomorphism.
\end{theorem} 

\begin{proof}
We show
the existence of the inverse map $\cF_d^{-1}$
directly.

%Assume that $d\in\{1,2,3\}$.
%For $i\in I\setminus\{0\}$
%let $z^- _{i,k;d}:=x^-  _{i,k;d}K_{i,d}^{- 1}K_{\isor _d;d}^k$.
%
%
%For the following equations we often use the formulas
%of Lemma~\ref{le:DUfms} 
%and the equations in Definition~\ref{definition:defUD},
%especailly Equations~\eqref{eq:prc1} and \eqref{eq:prc2},
%without refering to them.

Assume that $d\in \{1,2,3\}$,
$\{d,j,k\}=\{1,2,3\}$, and $j<k$.
Let $\tndX_d^+=-(r_ {d,j;d}r_ {d,k;d}r_ {d,0;d})^{-1}\ndX _d^+$.
One has 
\begin{align}\label{eq:mathbbX}
\tndX_d^+ & =\wlbr \wlbr \wlbr 
z^- _{d,1;d}, z^- _{j,0;d}\wrbr, z^- _{k,0;d}\wrbr, z^- _{d,0;d}\wrbr
%=-\wlbr \wlbr \wlbr 
%z^- _{j,1;d}, z^- _{d,0;d}\wrbr, z^- _{k,0;d}\wrbr, z^- _{d,0;d}\wrbr \\
% & 
=\wlbr \wlbr \wlbr 
z^- _{d,1;d}, z^- _{k,0;d}\wrbr, z^- _{j,0;d}\wrbr, z^- _{d,0;d}\wrbr,
%=-\wlbr \wlbr \wlbr 
%z^- _{k,1;d}, z^- _{d,0;d}\wrbr, z^- _{j,0;d}\wrbr, z^- _{d,0;d}\wrbr .
%\nonumber
\end{align} where we used the first part of 
\eqref{eq:caiso3} for the second equation.
By the first equation in \eqref{eq:mathbbX} and
Lemmas~\ref{le:DUfms}(1), \ref{le:caiso}(1),(2), one has 
\begin{equation}
	\begin{aligned}
& \wlbr \wlbr x^+ _{d,0;d},x^+ _{j,0;d}\wrbr, \wlbr x^+ _{d,0;d},
\tndX _d^+ \wrbr\wrbr 
=\wlbr x^+ _{d,0;d}, 
\wlbr x^+ _{j,0;d}, \wlbr x^+ _{d,0;d},
\tndX _d^+ \wrbr\wrbr\wrbr 
\\&\quad \quad \text{(by Eq.s~\eqref{eq:usfpr2},\eqref{eq:usf2}
and the equation $(x^+ _{j,0;d})^2=0$)}
\\
&\quad =[(\al_{d,d},\al_{0,d} )]_q
\wlbr x^+ _{d,0;d}, 
\wlbr x^+ _{j,0;d},
\wlbr \wlbr 
z^- _{d,1;d}, z^- _{j,0;d}\wrbr, z^- _{k,0;d}\wrbr\wrbr\wrbr \\
&\quad \quad \text{(by the second equation
in Eq.s~\eqref{eq:caiso6})} \\
&\quad =-[(\al_{d,d},\al_{0,d} )]_q
\wlbr x^+ _{d,0;d}, 
\wlbr x^+ _{j,0;d},
\wlbr \wlbr 
z^- _{k,1;d}, z^- _{d,0;d}\wrbr, z^- _{j,0;d}\wrbr\wrbr\wrbr
\\
&\quad \quad \text{(by the first equation
in \eqref{eq:caiso5aa})} \\
&\quad =-[(\al_{d,d},\al_{0,d} )]_q[(\al_{d,d},\al_{j,d} )]_q[(\al_{d,d},\al_{k,d} )]_q
z^- _{k,1;d}\,\,
\text{(by using \eqref{eq:usf})}.  
	\end{aligned}\label{eq:DUeq8}
\end{equation}
By Eq.s~\eqref{eq:prc1},\eqref{eq:usf}
one has
\begin{align}
& \wlbr \wlbr x^+ _{d,0;d},x^+ _{k,0;d}\wrbr, z^- _{k,1;d}\wrbr
=-\wlbr \wlbr x^+ _{k,-1;d},x^+ _{d,1;d}\wrbr, z^- _{k,1;d}\wrbr
=-[(\al_{d,d},\al_{k,d} )]_qx^+ _{d,1;d},  \label{eq:DUeq9} 
\end{align}
By the second equation in \eqref{eq:mathbbX},
using the same argument as above, one also has
the equations obtained from 
\eqref{eq:DUeq8},\eqref{eq:DUeq9}
by changing $j$ and $k$. Hence one conclude
that the following equation holds:
%By Equations~\eqref{eq:DUeq8} and \eqref{eq:DUeq9}
%and the fact that the definition of the element $\tndX _d^+$
%is symmetric with respect to $j$ and $k$, one has 
\begin{equation}
\begin{aligned}
&[(\al_{d,d},\al_{j,d} )]_q\wlbr \wlbr x^+ _{d,0;d},x^+ _{k,0;d}\wrbr, 
\wlbr \wlbr x^+ _{d,0;d},x^+ _{j,0;d}\wrbr, \wlbr x^+ _{d,0;d},
\tndX _d^+ \wrbr\wrbr\wrbr 
\\
&\quad 
=[(\al_{d,d},\al_{0,d} )]_q[(\al_{d,d},\al_{j,d} )]_q^2[(\al_{d,d},\al_{k,d} )]_q^2
x^+ _{d,1;d}, 
\\
&\quad 
=[(\al_{d,d},\al_{k,d} )]_q
\wlbr \wlbr x^+ _{d,0;d},x^+ _{j,0;d}\wrbr, 
\wlbr \wlbr x^+ _{d,0;d},x^+ _{k,0;d}\wrbr, \wlbr x^+ _{d,0;d},
\tndX _d^+ \wrbr\wrbr\wrbr .
\end{aligned} \label{eq:llstUr10}
\end{equation}
By \eqref{eq:llstUr10} and Lemma~\ref{le:DUfms}(3), 
one obtains
\begin{equation}
\begin{aligned}
&[(\al_{d,d}+\al_{0,d}, \al_{j,d}+\al_{d,d})]_q \wlbr 
\wlbr\wlbr x^+ _{d,0;d},\ndX _d^+\wrbr,\wlbr x^+ _{d,0;d},x^+ _{k,0;d}\wrbr\wrbr, 
\wlbr x^+ _{d,0;d},x^+ _{j,0;d}\wrbr\wrbr \\
&\, =
[(\al_{d,d}+\al_{0,d}, \al_{k,d}+\al_{d,d})]_q \wlbr 
\wlbr\wlbr x^+ _{d,0;d},\ndX _d^+\wrbr,\wlbr x^+ _{d,0;d},x^+ _{j,0;d}\wrbr\wrbr, 
\wlbr x^+ _{d,0;d},x^+ _{k,0;d}\wrbr\wrbr 
\end{aligned} \label{eq:comlstUr10} 
\end{equation}
%By Lemma~\ref{le:DUfms}, one also has
%\begin{align}
%& [x^+ _{j,0;d},\ndX _d^+]=0,\quad[x^+ _{k,0;d},\ndX _d^+]=0 
%\quad (\text{by Equation~\eqref{eq:mathbbX}}), 
%\label{eq:DUaxeq1} \\
%&\wlbr z^- _{d,1;d}, z^- _{j,0;d}\wrbr ^2=0,
%\label{eq:DUaxeq2}
%\\ & 
%\wlbr \wlbr \wlbr 
%z^- _{j,1;d}, z^- _{d,0;d}\wrbr, z^- _{k,0;d}\wrbr, z^- _{k,0;d}\wrbr
%=\wlbr z^- _{j,1;d}, \wlbr\wlbr 
%z^- _{d,0;d}, z^- _{k,0;d}\wrbr, z^- _{k,0;d}\wrbr \wrbr
%=0, \label{eq:DUaxeq3}
%\\
%&
%\wlbr x^+ _{d,0;d},\ndX _d^+\wrbr^2=0
%\quad (\text{by Equations~\eqref{eq:DUeq3} and \eqref{eq:DUaxeq2}, 
%\eqref{eq:DUaxeq2}}), \label{eq:DUaxeq4}
%\\
%&\wlbr\wlbr x^+ _{d,0;d},\ndX _d^+\wrbr ,\ndX _d^+\wrbr=0
%\quad (\text{by Equations~\eqref{eq:DUeq3} and \eqref{eq:DUaxeq4}}), 
%\label{eq:DUaxeq5}
%\\
%& \wlbr \tndX _d^+, z^- _{k,0;d}\wrbr
%= -\wlbr \wlbr\wlbr 
%z^- _{j,1;d}, \wlbr
%z^- _{d,0;d}, z^- _{k,0;d}\wrbr\wrbr, \wlbr z^- _{d,0;d}, z^- _{k,0;d}\wrbr\wrbr
%=0, \label{eq:DUaxeq6}
% \\
%& \wlbr \tndX _d^+, z^- _{j,0;d}\wrbr=0 , \label{eq:DUaxeq7}
%\\ & \wlbr \tndX _d^+, z^- _{d,0;d}\wrbr=0 
%\quad (\text{By Lemma~\ref{le:DUfms}(2)} ). \label{eq:DUaxeq8}
%\end{align}
%By Equation~\eqref{eq:DUeq2} and Lemma~\ref{le:DUfms}(1), 
%using $D\Psi_d$, one has
%\begin{align}
%&[ \wlbr \wlbr 
%x^\pm _{d,1;d}, x^\pm _{j,0;d}\wrbr, x^\pm _{k,0;d}\wrbr \wrbr , 
% x^\mp _{d,0;d} ] =0, \label{DUbxeq1} \\
%& [\ndX _d^+, \ndX _d^-]={\frac {K_{0,d}-K_{0,d}^{-1}} {q-q^{-1}}}
%\quad (\text{by Equation~\eqref{DUbxeq1} 
%and Lemma~\ref{le:DUfms}(4)}). \label{eq:DUbxeq2}
%\end{align}

Denote by $(U'_d)^{cl}$ 
the $\ndC$-subalgebra of $U'_d$ generated by the elements
$\sigma _d$, $K_{l,d}^{\pm{\frac 1 2}}$
($l\in\isv $) and $E_{r,d}$, $F_{r,d}$
($r\in\isv \setminus\{0\}$). By Theorem~\ref{th:HopfTri}(2),
one conclude that $(U'_d)^{cl}$ admits the presentation
with these generators and
the relations formed by Eq.s~\eqref{eq:cUr1}-\eqref{eq:cUr5}
and the relations $X=0$ for all elements $X$ in
\eqref{eq:Ur6}-\eqref{eq:Ur11}. 
Note that one also has the same fact
with $U'_d$ and $\isv$ in place of 
$(U'_d)^{cl}$ and $\isv \setminus\{0\}$
respectively.
Clearly, by Definition~\ref{definition:defUD}, 
one has a unique $\ndC$-algebra homomorphism
$(\cF^\prime_d)^{cl}:(U'_d)^{cl}\to DU'_d$ such that
$(\cF^\prime_d)^{cl}(\sigma _d)=\sigma _d$, 
$(\cF^\prime_d)^{cl}(K_{l,d}^{\pm{\frac 1 2}})=K_{l,d}^{\pm{\frac 1 2}}$
($l\in\isv $),  
$(\cF^\prime_d)^{cl}(E_{r,d})=x^+ _{r,0;d}$, 
$(\cF^\prime_d)^{cl}(F_{r,d})=x^- _{r,0;d}$
($r\in\isv \setminus\{0\}$). 
Then,
by Eq.s~\eqref{eq:comlstUr10}, 
\eqref{eq:caiso5},\eqref{eq:caiso7},\eqref{eq:caiso8},
using $D\Psi_d$, one has a unique
$\ndC$-algebra homomorphism
$\cF^\prime_d:U'_d\to DU'_d$ such that
%$\cF^\prime_d(\sigma _d)=\sigma _d$, 
%$\cF^\prime_d(K_{l,d}^{\pm{\frac 1 2}})=K_{l,d}^{\pm{\frac 1 2}}$
%($l\in\{0,1,2,3\}$),  
%$\cF^\prime_d(E_{r,d})=x^+ _{r,0;d}$, 
%$\cF^\prime_d(F_{r,d})=x^- _{r,0;d}$
%($r\in\{1,2,3\}$)  
$\cF^\prime_d(Y)=(\cF^\prime_d)^{cl}(Y)$
($Y\in (U'_d)^{cl}$)
and 
$\cF^\prime_d(E_{0,d})=\ndX _d^+$, 
$\cF^\prime_d(F_{0,d})=\ndX _d^-$.
By the equations in Definition~\ref{definition:defUD},
one concludes that as a $\ndC$-algebra,
$DU'_d$ is generated by the elements
$\sigma_d$, $K_{l,d}^{\pm{\frac 1 2}}$ ($l\in\{0,1,2,3\}$),
$x^\pm _{r,0;d}$ ($r\in\{1,2,3\}$) and $x^\pm _{k,\mp 1;d}$.
Hence Eq.s~\eqref{eq:DUeq8} and 
\eqref{eq:PDP=DPP} imply that the homomorphism 
$\cF^\prime_d$ is surjective. By Lemma~\ref{le:isofu1},
one obtains the equation $\cF_d\cF^\prime_d={\rm{id}}_{DU'_d}$.
Hence $\cF^\prime_d$ is injective. Thus one 
gets this theorem for $d\in\{1,2,3\}$.

The theorem for $d=4$ can be proved similarly, or more easily.
\end{proof}

The following lemma implies that 
Eq.s~\eqref{eq:expexp}-\eqref{eq:hikdxjkd} are
equivalent to
the ones of the original Drinfeld second realization.

\begin{lemma} \label{le:forhat}
Let $d\in\sd \setminus\{0\}$. Let
\begin{gather}
\hps^\pm _{i,\pm k;d}:=K_{i,d}^{\pm 1}
(\delta_{k0}+\Theta (k) (q-q^{-1})K_{\isor _d;d}^{\pm {\frac k 2}}\psi _{i,\pm k;d}),\,\,
(i\in \isv \setminus \{0\},\,k\in\mathbb{Z}) \label{eq:hatDG1} \\
\hh _{i,r;d}:=K_{\isor _d;d}^{- {\frac r 2}}h _{i,r;d}.\,\,(i\in \isv \setminus \{0\},\,r\in\mathbb{Z}\setminus \{0\})
\label{eq:hatDG2}
\end{gather} Then one has
\begin{gather}
K_{i,d}^{\pm 1}\exp(\pm (q-q^{-1})\sum_{r=1}^\infty z^r\hh_{i,\pm r;d})
=1+(q-q^{-1})\sum_{k=-\infty}^\infty z^k\hps_{i,\pm k;d}, \label{eq:hatexpexp}
\end{gather} (as equations of generating functions in $z$)
\begin{align}
[x^+ _{i,k;d},x^- _{j,l;d}] &=\delta_{ij}
{\frac {
K_{\isor ;d}^{{\frac {k-l} {2}}}\hps^+ _{i,k+l;d}
-K_{\isor _d;d}^{{\frac {l-k} {2}}}\hps^- _{i,k+l;d}
} 
{q-q^{-1}}}
\label{eq:hatxpxm}  \\
[\hh_{i,k;d},\hh_{j,l;d}] & =
\delta_{k,-l}{\frac {[k(\al_{i,d},\al_{j,d})]_q} k}
{\frac {K_{\isor _d;d}^k-K_{\isor _d;d}^{-k}}{q-q^{-1}}} \label{eq:hathikdhjkd}
\\
[\hh_{i,k;d},x^\pm_{j,l;d}] & =
\pm{\frac {[k(\al_{i,d},\al_{j,d})]_q} k}
K_{\isor _d;d}^{\mp {\frac {|k|} 2}}x^\pm_{j,k+l;d} \label{eq:hathikdxjkd}
\end{align}
\end{lemma} 
\begin{proof} This follows from \eqref{eq:expexp}-\eqref{eq:hikdxjkd}
by using  
the definition of 
the elements in \eqref{eq:hatDG1}-\eqref{eq:hatDG2}.
\end{proof} 

Now we give the Drinfeld second realization of $U'_d$.

\begin{theorem} \label{theorem:prDs}
Let $d\in \sd \setminus\{0\}=\{1,2,3,4\}$.

{\rm{(1)}} As a $\ndZ\Pi_d$-graded $\ndC$-algebra,
$DU'_d$ admits the presentation with the generators
$\sigma _d$,
$K_{u,d}^{\pm {\frac 1 2}}$, 
$x^\pm _{i,k;d}\in DU'_{d,\pm \alpha _{i,d}+k\isor _d}$, and
$\hh_{i,r;d}\in DU'_{d,r\isor _d}$,
where $u\in \isv $, $i\in \isv \setminus \{0\}$,
$k\in\mathbb{Z}$ and $r\in\mathbb{Z}\setminus \{0\}$,
and the defining relations obtained from
Eq.s~\eqref{eq:DRL1}, \eqref{eq:DRL2}, \eqref{eq:DRL3},
\eqref{eq:prc1}, \eqref{eq:prc2}, \eqref{eq:prc3},
\eqref{eq:hatxpxm}, 
\eqref{eq:hathikdhjkd}, and \eqref{eq:hathikdxjkd}
by determining the elements
$\hps^\pm_{i,k;d}$ for all $i\in I\setminus\{0\}$ and
all $k\in\ndZ$ by Eq.s~\eqref{eq:hatexpexp}.

{\rm{(2)}} There exists
a unique $\ndZ\Pi_d$-graded $\ndC$-algebra
isomorphism ${\widehat \cF}_d:DU'_d\to U'_d$ 
satisfying
the equations obtained from
Eq.s~\eqref{eq:cFg1}, \eqref{eq:cFg2},
\eqref{eq:cFg3}
by replacing $\cF _d$ with 
${\widehat \cF}_d$, and the equations
${\widehat \cF}_d (\hh _{i,r;d})=\epsilon _{i,i;d}^r
K_{\isor _d;d}^{- {\frac r 2}}
\bh_{i,r;d}$ for all
$i\in I\setminus\{0\}$ and
all $r\in\mathbb{Z}\setminus \{0\}$.
Further, ${\widehat \cF}_d$ coincides with $\cF _d$ as a map
and one has
${\widehat \cF}_d (\hps^\pm _{i,\pm l;d})=
\epsilon _{i,i;d}^l (q-q^{-1})K_{i,d}^{\pm 1}
K_{\isor _d;d}^{\pm {\frac l 2}}\bpsi _{i,\pm l;d}$
for all $i\in I\setminus\{0\}$ and
all $l\in\mathbb{N}$. 
\end{theorem}
\begin{proof}
This theorem holds by Definition~\ref{definition:defUD},
Theorem~\ref{theorem:iso}, 
and Lemma~\ref{le:forhat}.
\end{proof}

\subsection{Extension of the $\ndC$-algebras}

Recall the definitions of $\hat{\mfg }$ and $D^{(1)}(2,1;x)$ from
Section~\ref{sec:sectiontwo}.
Assume $\hat{\mfg }$ to be $D^{(1)}(2,1;x)$.
%Let $\hat{\mfg }:=D^{(1)}(2,1;x)$.
%(See Section~\ref{sec:sectiontwo}.)
Strictly $U'_d$ is the quantum superalgebra of
$[\hat{\mfg },\hat{\mfg }]$. 
Here we treat the quantum superalgebra
of $\hat{\mfg }$.
\begin{defin} \label{definition:extUd}
{\rm{(1)}}
Let $d\in\sd $. 
Define the additive group map $\chi_d:\ndZ\Pi_d\to\ndZ$ 
by $\chi_d(\al _{i,d})=\delta _{i0}$ 
for $i\in I$. The group ring of $\ndZ $ is the
commutative and cocommutative Hopf algebra
$\ndC [K_{\Lambda_0;d}^{\frac 1 2},K_{\Lambda_0;d}^{-\frac 1 2}]$, where
$$
K_{\Lambda_0;d}^{\frac 1 2}K_{\Lambda_0;d}^{-\frac 1 2}=1,\quad
\Delta (K_{\Lambda_0;d}^{\frac 1 2})=
K_{\Lambda_0;d}^{\frac 1 2}\otimes K_{\Lambda_0;d}^{\frac 1 2},\quad
\Delta (K_{\Lambda_0;d}^{-\frac 1 2})=
K_{\Lambda_0;d}^{-\frac 1 2}\otimes K_{\Lambda_0;d}^{-\frac 1 2}.
$$
Then $U'_d$ is a left
$\ndC [K_{\Lambda_0;d}^{\frac 1 2},K_{\Lambda_0;d}^{-\frac 1 2}]$-module
algebra \cite[Sect.\,4.1]{b-Montg93} with left action $\lact :
\ndC [K_{\Lambda_0;d}^{\frac 1 2},K_{\Lambda_0;d}^{-\frac 1 2}]\times
U'_d\to U'_d$ defined by
$$K_{\Lambda_0;d}^{{\frac 1 2}}\lact X_\mu = 
q^{{\frac {\chi_d(\mu)} 2}}X_\mu ,\qquad \mu\in \ndZ\Pi_d, 
X_\mu\in U' _{d,\mu}.$$
Let $U _d$ be the smash product algebra
\cite[Def.\,4.1.3]{b-Montg93}
$U_d:=U'_d\#\ndC [K_{\Lambda_0;d}^{\frac 1 2},K_{\Lambda_0;d}^{-\frac 1 2}]$.

\rm{(2)} Let $d\in\sd\setminus\{0\}$. Similarly to the construction in
Part (1) define the smash product algebra
$DU_d:=
DU'_d\#\ndC [K_{\Lambda_0;d}^{\frac 1 2},K_{\Lambda_0;d}^{-\frac 1 2}]$.
\end{defin}

We extend Theorem~\ref{theorem:prDs} to that for $U _d$.
 
\begin{theorem}\label{theorem:fin}
The map ${\widehat \cF} _d$ 
can be extended to a $\ndC$-algebra isomorphism
from $DU _d$ to $U _d$ by letting 
${\widehat \cF} _d (K_{\Lambda_0,d}^{{\frac m 2}})= 
K_{\Lambda_0,d}^{{\frac m 2}}$
for all $m\in\ndZ$.
\end{theorem}

\begin{proof}
This theorem follows from Theorem~\ref{theorem:prDs} and 
Definition~\ref{definition:extUd}.
\end{proof}

\vspace{1cm}
\begin{center}
{\bf{Acknowledgments}}
\end{center} The authors would like to thank the referee
for the careful reading and valuable comments.
This work is supported in part by funds provided by the U.S. Department 
of Energy (D.O.E.) under cooperative research agreement DEFG02-05ER41360.
A.\,Torrielli thanks INFN (Istituto Nazionale di Fisica Nucleare) for 
supporting him through a ``Bruno Rossi'' postdoctoral fellowship.

\newpage
{\small \sc Istv\'an Heckenberger, 
Mathematisches Institut,
Universit\"at Leipzig,
Pf 100920, D-04009 Leipzig, Germany
%Ludwig-Maximilians-Universit\"at,
%Mathematisches Institut,
%Theresienstr. 39,
%D-80333, 
%M\"unchen,
%Germany
}

{\small \textit{E-mail address:} \texttt{Istvan.Heckenberger@math.uni-leipzig.de}}

\vspace{.5\baselineskip}

{\small \sc Fabian Spill, 
Humboldt--Universit\"at zu Berlin,
Institut f\"ur Physik,
Newtonstra\ss e 15,
D-12489 Berlin, Germany

%Theoretical Physics,
%Imperial College London,
%South Kensington campus,
%Prince Consort Road,
%London,
%SW7 2AZ,
%UK
}

{\small \textit{E-mail address:} 
\texttt{spill@physik.hu-berlin.de}}
%\texttt{fabian.spill@imperial.ac.uk}
\vspace{.5\baselineskip}

{\small \sc Alessandro Torrielli,
Center for Theoretical Physics,
Laboratory for Nuclear Sciences,
and
Department of Physics,
Massachusetts Institute of Technology,
Cambridge, Massachusetts 02139, USA}

{\small \textit{E-mail address:} \texttt{torriell@mit.edu}}

\vspace{.5\baselineskip}

{\small \sc Hiroyuki Yamane,
Department of Pure and Applied Mathematics,
Graduate School of Information Science
and Technology, Osaka University, Toyonaka 560-0043,
Japan}

{\small \textit{E-mail address:} \texttt{yamane@ist.osaka-u.ac.jp}}

\end{document}